\documentclass{article}
\usepackage[a4paper,left=3cm,top=3.5cm]{geometry}
\usepackage{verbatim}
\usepackage[all,cmtip]{xy}

\usepackage{mathtools}

\usepackage{CJK,CJKspace,CJKnumb}
\usepackage{latexsym,bm,easymat}
\usepackage{amsmath,amscd,amssymb,amsfonts,mathdots,ntheorem}
\usepackage{indentfirst}
\usepackage[colorlinks,linkcolor=blue,citecolor=red]{hyperref}
\usepackage{xcolor,graphicx,blkarray}
\usepackage{latexsym}
\usepackage{longtable}
\usepackage{colortab}
\usepackage{cases}
\usepackage{slashed}

\usepackage{shorttoc}
\begin{document}
\makeatletter
\DeclareRobustCommand\widecheck[1]{{\mathpalette\@widecheck{#1}}}
\def\@widecheck#1#2{%
    \setbox\z@\hbox{\m@th$#1#2$}%
    \setbox\tw@\hbox{\m@th$#1%
       \widehat{%
          \vrule\@width\z@\@height\ht\z@
          \vrule\@height\z@\@width\wd\z@}$}%
    \dp\tw@-\ht\z@
    \@tempdima\ht\z@ \advance\@tempdima2\ht\tw@ \divide\@tempdima\thr@@
    \setbox\tw@\hbox{%
       \raise\@tempdima\hbox{\scalebox{1}[-1]{\lower\@tempdima\box
\tw@}}}%
    {\ooalign{\box\tw@ \cr \box\z@}}}
\makeatother

%%%%%%%%%%%define of the proof=============
\def\comp{\ensuremath\mathop{\scalebox{.6}{$\circ$}}}
\def\QEDclosed{\mbox{\rule[0pt]{1.3ex}{1.3ex}}} %
\def\QEDopen{{\setlength{\fboxsep}{0pt}\setlength{\fboxrule}{0.2pt}\fbox{\rule[0pt]{0pt}{1.3ex}\rule[0pt]{1.3ex}{0pt}}}}
\def\QED{\QEDopen} %
\def\pf{\noindent{\bf Proof}} %
\def\endpf{\hspace*{\fill}~\QED\par\endtrivlist\unskip \hfill}
\def\Dirac{\slashed{D}}
%%%%%%%%%%%%%%%% Define the symbol and notation =====
\def\Hom{\mbox{Hom}}
\def\ch{\mbox{ch}}
\def\Ch{\mbox{Ch}}
\def\Tr{\mbox{Tr}}
\def\End{\mbox{End}}
\def\Re{\mbox{Re}}
\def\cs{\mbox{Cs}}
\def\spin{\mbox{spin}}
\def\Ind{\mbox{Ind}}
\def\Hom{\mbox{Hom}}
\def\ch{\mbox{ch}}
\def\cs{\mbox{cs}}
\def\nTbe{\nabla^{T,\beta,\epsilon}}
\def\nbe{\nabla^{\beta,\epsilon}}
\def\mfa{\mathfrak{a}}
\def\mfb{\mathfrak{b}}
\def\mfc{\mathfrak{c}}
\def\mcM{\mathcal{M}}
\def\mcC{\mathcal{C}}
\def\mcB{\mathcal{B}}
\def\mcG{\mathcal{G}}
\def\mcE{\mathcal{E}}
\def\mfs{\mathfrak s}
\def\mcWs{\mathcal W_{\mathfrak s}}
\def\umcF{\underline{\mathcal F}}
\def\mce{\mathfrak e}
\def\mbT{\mathbf T}
\def\mbc{\mathbf c}
\def\spinc{$\mbox{spin}^c~$}

%%%%%%%%%%
%email-address
\newcommand{\at}{\makeatletter @\makeatother}
\renewcommand{\contentsname}{\center{Content}}
\newtheorem{defi}{Definition}[section]
\newtheorem{thm}[defi]{\textbf{Theorem}}
\newtheorem{conj}[defi]{\textbf{Conjecture}}
\newtheorem{cor}[defi]{\textbf{Corollary}}
\newtheorem{lemma}[defi]{\textbf{Lemma}}
\newtheorem{exa}[defi]{\textbf{Example}}
\newtheorem{prop}[defi]{Proposition}
\newtheorem{assump}[defi]{Assumption}
\newtheorem{hypo}[defi]{Hypothesis}
%%%%%%%%-----------------title

\title{Monopole Floer homology for codimension-3 Riemannian foliations}
\date{}
\author{Dexie Lin}
\maketitle

%\begin{center}  Dedicated to my parents.\end{center}

%%%%%%%%%%%%%%%%%%%%%%%%%%%%%%%%%%%%%
\begin{abstract}
  In this paper, we give a systematic study of Seiberg-Witten theory on a closed
   oriented manifold $M$ with a codimension-$3$ oriented Riemannian foliation $F$. Under a  certain topological condition,  we   construct the basic Seiberg-Witten invariant and the monopole Floer homologies $\overline{HM}(M,F,\mfs;\Gamma),~\widehat{HM}(M,F,\mfs;\Gamma),
   ~\widecheck{HM}(M,F,\mfs;\Gamma)$,  for each transverse \spinc structure $\mfs$, where $\Gamma$ is a complete local system. We will show that these homologies are independent of the bundle-like metrics and generic perturbations. To define the basic Seiberg-Witten invariant for manifolds with codimension 3 Riemannian foliations, we do not need the taut condition. The main difference between the  basic monopole Floer homologies and the ones on manifolds is the necessity to use the Novikov ring  on basic monopole Floer homologies. %Using the similar idea and methods, we also  construct  the monopole Seiberg-Witten-Floer homology groups for closed oriented $3$-orbifold, and show that they are independent of the metric and generic perturbation.
\end{abstract}

\tableofcontents

%%%%%%%%%%%%%%%

 \section{Introduction}
The interaction between geometry and partial differential equation in dimension $4$ is a theme which runs through a great deal of   work by many mathematicians on gauge theory over the past decades. In particular, the Seiberg--Witten theory's development  is one of the main motivation for the study of differential
topology and low dimensional manifolds. Since the foundational paper \cite{W} by Witten, a lot of work has
been done to apply this theory to various aspects of $3$ and $4$-dimensional
manifolds. Seiberg-Witten theory can   be generalized for studying the orbifolds(see Baldridge's work \cite{Bald} for the extension to $3$-orbifolds). This article lays the groundwork for the case in which the higher-dimensional
manifold admits a Riemannian foliation of codimension $3$. From the viewpoint of analysis, gauge theory is closely related to the study of (nonlinear)Fredholm operator and the index of its linearized operator.
 To extending the framework of gauge theory to the manifold with Riemannian foliation, a natural idea is to study the transverse (nonlinear)elliptic operator on  Riemannian foliation. For instance, the compactness of the basic Seiberg-Witten moduli space for manifolds with  codimension $4$ Riemannian foliations
 is showed by  Kordyukov,  Lejmi and Weber \cite{KLW}.  The author gives a construction of basic cohomotopy Seiberg-Witten invariant for codimension 4 Riemannian foliation \cite{Lin}. In the same paper, the author gives an application to the basic index of the basic Dirac operator.

  The
theme of this article is to generalize the well-known constructions  of Seiberg--Witten theory in $3$-manifolds   to the manifolds with codimension $3$ Riemannian foliations.
It is worth to note that in codimension 3 case, we do not need  the taut condition to define the basic Seiberg-Witten invariant, which is quite different to the codimension 4 case. %In particular, it is a widely generalized theory for orbifolds, which can often be realized as a particular class of foliations.
For any closed oriented $3$-manifold $M$ with a \spinc structure $\mfs$, Kroheimer and Mrowka gave the construction of monopole  Floer homologies  $\overline{HM}(M,\mfs),~\widehat{HM}(M,\mfs),~\widecheck{HM}(M,\mfs)$  in their celebrated book \cite{KM}. The main purpose of this paper is to construct the monopole   Floer homologies for the manifold with a codimension $3$ Riemannian foliation $(M,F)$ satisfying a certain condition.
The idea is to apply  the  arguments of the non-exact perturbation of the monopole Floer homologies to the case of Riemannian foliation. The following theorem could summarize the result of this paper.

\begin{thm}
  Let $(M,F)$ be an oriented closed manifold with codimension 3 oriented taut Riemannian foliation $F$ and admit a transverse \spinc structure $\mathfrak s$. Suppose  that  $H^1_b(M)\cap H^1(M,\mathbb Z)\subset H^1(M)$ is a lattice of $H^1_b(M)$. Then, using a bundle-like metric $g$, a generic perturbation $\eta$ and the Novikov ring $\Gamma$,   we construct the basic monopole  Floer homologies $$\overline{HM}(M,F,\mathfrak s, g,\eta;\Gamma),~\widehat{HM}(M,F,\mathfrak s, g,\eta;\Gamma),~\widecheck{HM}(M,F,\mathfrak s, g,\eta;\Gamma),$$  Moreover,  these homologies are independent of the  bundle-like metrics  and the generic perturbations , which are denoted by
         \[\overline{HM}(M,F,\mathfrak s ;\Gamma),~\widehat{HM}(M,F,\mathfrak s;\Gamma ),~\widecheck{HM}(M,F,\mathfrak s ;\Gamma).\]
\end{thm}

Some notations and terms will be defined in Section 7.
The necessity of using the Novikov ring to construct the homologies $\overline{HM}(M,F,\mathfrak s ),~\widehat{HM}(M,F,\mathfrak s ),~\widecheck{HM}(M,F,\mathfrak s )$   is reflected in defining the partial operator of the Floer complex.% Explicitly speaking,  there might be infinitely many terms in general.

The structure of this paper is as follows: in Section 2, we review some notions and necessary properties about the Riemannian foliation;
in Section 3,  some analysis properties for some transverse
equations will be given, which are necessary for the later sections; in Section 4, we construct the basic Seiberg-Witten invariant for manifold with codimension $3$ Riemannian foliation; in Section 5, we construct the basic Chern-Simons-Dirac functional and give some properties of it; in Section 6,  the gluing theorem for the basic moduli spaces will be proved,  which is essential to  construct the basic monopole  Floer homologies; in Section 7, the proof of the above theorem will be given; in Section 8, we construct the   monopole Floer homologies for a certain kind of $3$ orbifolds, and we give   some examples and a method to construct the Riemannian foliation satisfying the assumption of the above theorem.

\vspace{3mm}

{\bf Acknowledgement:} The author  warmly thanks Mikio Furuta for his long time invaluable help in both mathematics and life. The author is  grateful to Kim. A. Fr\o yshov for his helpful discussion on Floer homology and Ken. Richardson for the discussion on Riemannian foliation. The research is partially sponsored by the FMSP of The University of Tokyo and The Fundamental Research Funds for the Central Universities No. 2021CDJQY-009..

 \section{Preliminary}

 In this section, we review some results about the geometry of Riemannian foliations and foliated bundle. %In the first subsection, we review the classical results of the geometric foliation; in the second subsection, we review the Seiberg-Witten theory for codimension $4$ foliation.
% \subsection{Geometry of Foliation}
  Let $M$ be a closed oriented $n$ dimensional manifold with rank $p$ foliation $F$, and let $Q=TM/F$ be the quotient bundle. We  denote the codimension of this foliation by $q=n-p$. For more details of this section, we give a reference \cite{Ton2}.
 \begin{defi}
   A Riemannian metric $g_Q$ on $Q$ is said to be bundle-like, if it holds that 
   \[L_Xg_Q\equiv0,\]
   for any $X\in \Gamma(F)$. We say $(M,F)$ is a Riemannian foliation, if $Q$ admits a bundle-like metric.
 \end{defi}

 Given a metric $g$ on $TM$, $Q$ is identified with the orthogonal complement to $F^\perp$ by $g$. In turn, $Q$ inherits a metric $g_{F^\perp}$, where $g_{F^\perp}=g|_{F^\perp}$. We have the following equivalence,
 \[\mbox{a metric }g\mbox{ of }TM\mbox{ corresponds a triple }(g_F,\pi_F, g_Q),\]
 where $g_F=g|_{F}$ and $\pi_F$ is the projection $TM\to F$.

\noindent
A Riemannian metric $g$ on $TM$ is said to be \emph{bundle-like}, if the induced metric $g_{F^\perp}$ is bundle-like.
 By the work of Reinhart $\cite{Reinhart}$, it is known  that the bundle-like metric can be locally written as $g=\sum_{i,j}g_{ij}(x,y)\omega^i\otimes \omega^j+\sum_{k,l}g_{k,l}(y)dy^k\otimes dy^l$, where $(x,y)$ is in the foliated chart of $M$ and $\omega^i=dx^i+a^i_\alpha(x,y) dy^\alpha$.
In this paper, we  assume that $(M,F)$ is a manifold with Riemannian foliation without specific mention.
Let $\pi_Q:TM\to Q$ be the canonical projection.  We define a connection $\nabla^{T}$ on $Q$, by
$$\nabla^{T}_Xs=\begin{cases}
   \pi_Q([X,Z_s])& X\in \Gamma(F),\\
   \pi_Q(\nabla^g_X Z_s)& X\in \Gamma(F^\perp),
 \end{cases}$$
 for any section $s\in\Gamma(Q)$,
 where $Z_s\in \Gamma(TM)$ is a lift of $s$, i.e. $\pi_Q(Z_s)=s$ and $\nabla^g$ denotes the Levi-Civita connection of $g$.  We call $\nabla^T$ transverse Levi-Civita connection. If $(M,F)$ is a Riemannian foliation, then by the Koszul-formula \cite[Theorem 5.9]{Ton2},  it is clear that $\nabla^T$ is  uniquely determined by $g_{Q}$. Moreover, one can verify that it is  torsion free and metric-compatible, whose leafwise restriction coincides with the Bott-connection. We set $R^T$ as the curvature of this connection.
 Similarly, we define the transverse Ricci curvature and scalar curvature by
 \[Ric^{T}(Y)=\sum^{q}_{i=1}R^{T}(Y,e_i)e_i,~
 Scal^{T}=\sum^{q}_{i=1}g_Q(Ric^{T}(e_i),e_i),\]
 where $\{e_i\}$ is a local  orthonormal frame of $Q$.
 % One has that $R^T$ satisfies the   condition \[\iota_XR^T=0,\]for any $X\in\Gamma(F)$.
  We define the basic forms as follows:
 \[\Omega^r_b(M)=\{\omega\in\Omega^r(M)\big|~\iota_X(\omega)=0,~L_X(\omega)=0,
 \mbox{ for all } X\in\Gamma(F)\}.\]

 By the work of Alvarez L\`opez \cite{AL},  the following $L^2$ orthogonal decomposition holds  for the forms on $M$, i.e.
 \[\Omega(M)=\Omega_b(M)\oplus\Omega^\perp_b(M),\]
 with respect to the $C^\infty$-Fr\'echet topology.

 Choosing  a local orthonormal basis $\{e_i\}_{1\leq i\leq p}$ of $F$, we define the character form $\chi_F$ of the foliation by,
 $\chi_F(Y_1,\cdots, Y_p)=\det(g_F(e_i,Y_j))_{1\leq i,j\leq p}$,
 for any section $Y_1,\cdots, Y_p\in \Gamma(TM)$. By the metric $g_Q$($g_{F^\perp}$),  we  have the basic Hodge-star operator, $$\bar*:\bigwedge^rQ^*\to \bigwedge^{q-r}Q^*.$$
 The basic Hodge-star operator is  related to the usual Hodge-star operator by the formula $\bar*\alpha=(-1)^{(q-r)\dim(F)}*(\alpha\wedge\chi_F)$. Moreover, we have $\bar *:\Omega^r_b(M)\to \Omega^{q-r}_b(M)$ and the volume density formula,
 $dvol_M=dvol_Q\wedge\chi_F$.
 For a section $\alpha\in \Omega^r_b(M)$, we define its $L^2$ norm by
\[\|\alpha\|^2_{L^2}=\int_M\alpha\wedge\bar*\alpha\wedge\chi_F.\]
We set $d_b$ as the restriction of $d$ to the basic forms, the complex $d_b:\Omega^r_b(M)\to\Omega^{r+1}_b(M)$ is a subcomplex of the deRham complex, whose cohomology is called  basic cohomology, and  denoted by $H^r_b(M)$. It is known that $H^1_b(M)\subset H^1(M)$.
 We denote by $b^r_b=\dim H^r_b(M)$. Before going on, we introduce the mean curvature field and mean curvature form.
 %The main difference on analysis between the Riemannian manifolds and Riemannian foliations is the
% For any section $\alpha\in \Gamma(\bigwedge^rQ^*)$, we have that \[\int_M\alpha\wedge\bar*\alpha\wedge\chi_F=\int_M\alpha\wedge*\alpha.\]
%For a bundle-like metric, we  have $\bar*:\Omega^r_b(M)\to\Omega^{q-r}_b(M).$
\begin{defi}
  The mean curvature vector field is defined by $\tau=\sum^{\dim F}_{i=1}\pi_Q(\nabla^g_{\xi_i}\xi_i)\in\Gamma(Q)$, where $\{\xi_i\}$ is a local orthonormal basis of $F$. Let $\kappa\in \Gamma(Q^*)$ be the dual to $\tau$ via the metric  $g_Q$.
\end{defi}

The mean curvature form is related to the volume of the foliation, by the following proposition.

\begin{prop}[Rummler \cite{Rummler}]\label{formula-Rummler}
For any metric $g$ on $TM$,
  it holds that
  \[d\chi_F=-\kappa\wedge\chi_F+\phi_0,\]
  where $\phi_0$ belongs to $F^2\Omega^{p+1}=\{\omega\in\Omega^{p+1}(M)\big|\iota_{X_1}\cdots\iota_{X_p}\omega=0,\mbox{ for any }X_1,\cdots,X_p\in \Gamma(F)\}$. This implies that $w\wedge\phi_0=0$ for any $w\in \Gamma(\bigwedge^{q-1}Q^*)$.
\end{prop}

% By direct calculation \cite[Chpater 7]{Ton2},    we have that the $L^2$ formal adjoint  of $d$ is $\delta=(-1)^{m(*+1)+1}\bar*(d-\kappa\wedge)\bar*$.
% \end{prop}
% \begin{pf}   By direct calculation, we have for any sections $\alpha,~\beta\in \Omega^*(Q^*)$, the following identity holds:   \begin{eqnarray*}     d(\alpha\wedge\bar*\beta\wedge\chi_F)&=&((d\alpha)\wedge\bar*\beta)\wedge\chi_F+   (-1)^{r-1}(\alpha\wedge(d\bar*\beta))\wedge\chi_F\\   &&+(-1)^{m-1}(\alpha\wedge\bar*\beta)\wedge   (-\kappa\wedge\chi_F)+(\alpha\wedge\bar*\beta)\wedge\phi_0.   \end{eqnarray*}   When $\alpha\in\Omega^{r-1}(Q^*)$ and $\beta\in\Omega^r(Q^*)$, one can deduce that $(\alpha\wedge\bar*\beta)\wedge\phi_0\equiv0$. Taking the integral of the above formula, we have   \begin{eqnarray*}     (d\alpha,\beta)_{L^2}&=&(-1)^r\int_M\alpha\wedge(d\bar*\beta-(-1)^{m-r}    \bar*\beta\wedge\kappa)\wedge\chi_F\\     &=&(-1)^{m(r+1)+1}\int_M\alpha\wedge\bar*(\bar*(d-\kappa\wedge)\bar*\beta)\wedge     \chi_F.   \end{eqnarray*} \end{pf}

By the decomposition, we have that  $\kappa=\kappa_b+\kappa_0$ for a bundle-like metric $g$, where $\kappa_b\in \Omega^1_b(M)$ and  $(\kappa_0,\omega_b)_{L^2}=0$ for any basic one form $\omega_b$.  Dominguez \cite{D} shows that
 any Riemannian foliation $F$ carries a tense bundle-like metrics, i.e. having basic mean curvature form $\kappa=\kappa_b$.
The form $\kappa_b$ is called the basic mean curvature form.  It is known that
 $d\kappa_b=0$,
 and the cohomology class $[\kappa_b]$ is independent of any bundle-like metric \cite{AL}.

\begin{defi}
  We say a foliation is \emph{taut}, if  there is a metric on $M$ such that $\kappa=0$, i.e. all leaves are minimal submanifolds.
\end{defi}
In this paper,
a bundle-like metric is called to be taut, if the induced mean curvature form vanishes.
For a fixed Riemannian foliation $F$, the taut condition has a topological obstruction.

\begin{prop}[Alvarez L\`opez \cite{AL}]\label{prop-taut}
  Let $F$ be a  Riemannian foliation on a closed manifold. Then, $F$ is taut if and only if the class $[\kappa_b]$ is trivial. Furthermore, when $F$ is transversely oriented  the foliation is taut if and only if $H^q_b(M)\neq0$.
\end{prop}

By \cite[Page 99]{Ton2},   the Poincare duality holds for the basic cohomologies  under the  taut  condition, i.e. $H^r_b(M)\cong H^{q-r}_b(M)$.
%\begin{prop}[ Dominguez \cite{D}]\end{prop}\textbf{Remark}: It is known that any bundle-like metric can be deformed inthe leaf directions leaving the transverse part unchanged in such a way that the mean curvature form becomes basic.In this article, we always let  $\kappa$  be basic, i.e. $\kappa=\kappa_b$.

 \begin{prop}[Tondeur {\cite[Theorem 7.18]{Ton2}}]
   Let $d_b$ denote the restriction of $d$ on the basic forms. %For any basic form $\alpha$ we define the $L^2$-norm by \[\|\alpha\|^2_{L^2}=\int_M\alpha\wedge\bar*\alpha\wedge\chi_F.\]
   Then, the   $L^2$-formal adjoint   of $d_b$ is $\delta_b=(-1)^{q(*+1)+1}\bar*(d_b-\kappa_b\wedge)\bar*$.
 \end{prop}
 We define the basic Laplacian operator by $\Delta_b=d_b\delta_b+\delta_bd_b$.

\noindent
Now,
we review the definitions of foliated vector bundle and basic connections.

 \begin{defi}
   A principal bundle $P\to M$ is called foliated, if it is equipped with a lifted foliation $F_P$ invariant under the structure group action, such that it is transversal to  the tangent space to the fiber and $F_P$ projects isomorphically onto $F$. We say a vector bundle $E\to M$ is foliated, if its principal bundle $P_E$ is foliated.
 \end{defi}

 \begin{defi}
   A connection $\omega$ of the foliated principal bundle $P$ is called adapted, if the horizontal distribution associated to this connection   contains the foliation $F_P$. A covariant connection on a foliated vector bundle is called adapted, if its associated connection on the principal bundle is. We say an adapted connection $\omega$ is called basic, if it is a Lie algebra valued basic form. Similarly, an adapted covariant connection is called basic, if its principal connection is.
 \end{defi}

% For a foliated vector bundle $E$ over $M$, by the definition, it is  known that for any two connections $\nabla^1$ and $\nabla^2$ adapted to this foliated vector bundle, we have that \[\nabla^1_Vs=\nabla^2_Vs,\] for all $s\in\Gamma(M,E)$ and $V\in\Gamma(F)$.
 Using an adapted connection,  we define the basic sections by \[\Gamma_b(E)=\{s\in \Gamma(E)\big|~ \nabla_Xs\equiv0,~\mbox{for all }X\in \Gamma(F)\},\]
 where $\nabla$ is an adapted connection. It is  known that the space of  basic sections is independent of the choice of the adapted connection.
 \begin{defi}
   A transverse Clifford module $E$ is a complex vector bundle over $M$ equipped with a hermitian metric  satisfying the following properties:
 \begin{enumerate}
   \item  $E$ is a  bundle of $Cl(Q)$-modules, and the Clifford action $Cl(Q)$ on $E$ is skew-symmetry, i.e.  \[(s\cdot\psi_1,\psi_2)+(\psi_1,s\cdot\psi_2)=0,\]
 for any $s\in\Gamma(Q)$ and $\psi_1,\psi_2\in\Gamma(E)$;
   \item $E$ admits a basic metric-compatible connection, and this connection is compatible with the Clifford action.
 \end{enumerate}
 \end{defi}

 We say that $(M,F)$ admits a transverse \spinc structure, if $Q$ is \spinc and the associated spinor bundle  $S$ is a transverse Clifford module over $(M,F)$.

% {\bf Remark:} Actually, for any  complex Hermitian foliated bundle $E$,   such that $E$ is a transverse Clifford module, self-adjoint and equipped with a Hermitian Clifford connection $\nabla^E$, by \cite{KT} we can choose $\nabla^E$ to be a \emph{basic connection} which means that theconnection and curvature forms of $\nabla^E$ are  basic forms. In this paper, we always assume that the connection is basic.

%\hfill

% For a foliated \spinc principle bundle $P$ over $M$, the automorphism  group $G_b$ is a subgroup of $G=\{u:M\to U(1)\}$  satisfies that\[Vu=0,\] for all $V\in\Gamma(F)$ for $u\in G$.

% {\bf Example}:For $T^4$ with dense a $1$-dimensional foliation $V$.  We have that for any foliated \spinc structure, the automorphism group $G_b$ consists of all constant maps.To see this, we fix a point $x_0\in T^4$, choose a metric on $T^4$. For any $x\in T^4$, and $\epsilon>0$, there is a geodesic ball at $x_0$ with radius $\delta$, such that for any point $x$ in this ball, we have \[|f(x)-f(x_0)|<\epsilon,\]for any smooth function $f$. As the foliation is dense, there is a constant $T(\delta)$, such that $\exp(T(\delta)x_0)$ locates in this ball. Since for any $u\in G_b$, $u$ is constant along this orbit $\exp(tx_0)$, by the above argument, we get \[|u(x)-u(x_0)|<\epsilon,\] hence $u(x)=u(x_0)$.

\begin{defi}\label{defi-transverse-Dirac}Let $E$ be a  transverse Clifford module $E$ over $(M,F)$.
   Fixing a basic connection $\nabla^E$, we define the transverse  Dirac operator $\Dirac^T$ by $\Dirac^T=\sum^q_{i=1}e_i\cdot\nabla^E_{e_i}$ action on $\Gamma(E)$, where $\{e_i\}$ is a local orthonormal basis of $Q$.
\end{defi}
 Note that $\Dirac^T$ is not formally self adjoint in general, whose adjoint operator is $\Dirac^{T,*}=\Dirac^T-\tau_b$. We set $\Dirac_b=\Dirac^T-\frac12\tau_b$, which is called basic Dirac operator.
 %For the non-taut case, we define $\Dirac_b=\sum^q_{i=1}e_i\cdot\nabla^E_{e_i}-\frac12\tau_b$.
%By straightforward calculation,
This basic Dirac operator is a formally self-adjoint operator and maps the basic sections $\Gamma_b(E)=\{s\in\Gamma(M,E)\big|\nabla_Xs\equiv0,\mbox{ for any }X\in\Gamma(F)\}$ to itself.
%   \begin{prop}   The basic Dirac operator maps the basic sections $\Gamma_b(S)=\{s\in\Gamma(M,S)\big|\nabla_Xs\equiv0,\mbox{ for any }X\in\Gamma(F)\}$ to basic sections. \end{prop}
% \begin{pf}   By the straightforward calculation, one gets that  \begin{eqnarray*}     \nabla_X(\Dirac_b\psi)&=&\sum^m_{i=1}\nabla_X(e^i\nabla^E_{e_i}\psi)\\     &=&\sum^m_{i=1}(\nabla^T_X(e^i))\nabla^E_{e_i}\psi+     \sum^m_{i=1}e^i(\nabla^E_X\nabla^E_{e_i}\psi)\\     &=&\sum^m_{i=1}(\nabla^T_X(e^i))\nabla^E_{e_i}\psi+     \sum^m_{i=1}e^i(\nabla^E_{e_i}\nabla^E_X\psi)\\     &&+\sum^m_{i=1}e^i(\nabla_X(\nabla^E))     (e_i,\psi)+\sum^m_{i=1}e^i\nabla^E_{[X,e^i]}\psi.   \end{eqnarray*}   Letting $X\in\Gamma(F)$ and $\psi\in \Gamma_b(S)$, on the above formula we have that the second and third terms vanish,  the first and last terms cancel. \end{pf}

Let $E$ be a foliated vector bundle on $M$ equipped with a basic Hermitian structure
and a compatible basic connection $\nabla^E$. We define the basic $\|\|_{L^p_k}$-norm by
\[\|u\|_{L^p_k}=\sum^k_{j=0}(\int_M|(\nabla^E)^ju|^pdvol_M)^{\frac1p},\]
for any $u\in \Gamma_b(E)$. Let $L^p_k$ be the completion of $\Gamma_b(E)$ with respect to this norm.
One has the similar Sobolev embedding and Sobolev multiplication properties for basic sections, which are shown in \cite[Theorem 9, 10, 11]{KLW}.% To make this paper complete, here we just give the statements.

%\begin{thm}Suppose that $(M,F)$ is a closed oriented manifold with codimension-$m$ foliation $F$, then the following inclusions hold:\begin{itemize}  \item \[L^p_k\hookrightarrow L^q_l,\]  where integers $l,~k$ satisfying $0\leq l\leq k$ and $l-\frac mk\leq k-\frac mp$.  \item \[L^p_k\hookrightarrow C^l,\]  where $l<k-\frac mp$.\end{itemize}\end{thm}

%\begin{thm}%[ \cite{KA}]Let $0\leq l\leq k$, under the setting of above theorem, we have the following continuous maps: \begin{itemize}    \item \[L^p_k\times L^q_l\to L^q_l,\]    where $k - \frac mp > 0$ and $k -\frac mp > l- \frac mq$.    In particular, if $k=l,~p=q$ and $k-\frac mp>0$, then    \[L^p_k\times L^p_k\to L^p_k.\]    \item \[L^p_k\times L^q_l\to L^r_t,\]    where $k-\frac mp<0$, $l-\frac mq<0$ and $l$ satisfies $0\leq t\leq l$,$r$ satisfies $0<\frac tm+\frac1p-\frac km+\frac1q-\frac lm\leq \frac1r\leq1$.    \item \[(L^p_k\cap L^\infty)\times(L^q_l\cap L^\infty)\to (L^q_l\cap L^\infty),\]        where $k=\frac mp$ and $l-\frac mq\leq0$.    \item\[(L^p_k\cap L^\infty)\times L^q_l \to L^q_l,\]    where $l-\frac mp<0$. \end{itemize}\end{thm}

%Before the construction of the Bauer-Furuta invariant, we need to define what the transverse elliptic operator is.

%Let $E_1,~E_2$ be two foliated vector bundle on $M$ with compatible basic connections.

\begin{defi}\label{defi-transverse-elliptic}
Let $E_1$ and $E_2$ be two foliated vector bundles over $M$ with compatible basic connections. A differential operator $L:\Gamma(E_1)\to \Gamma(E_2)$, is called basic, if in any foliated chart $(x,y)\in U\times V$ with distinguished local trivialization of $E_1$ and $E_2$, then one locally writes
\[L|_{U\times V}=\sum a_\alpha(y)\frac{\partial^{|\alpha|}}{\partial^{\alpha_1} y_1\cdots\partial^{\alpha_q}y_q}.\]

   A basic differential operator $L$ defined as the above, is said to be  \emph{transverse elliptic},  if its transverse symbol is an
isomorphism away from the $0$-section, i.e. $\sigma(p,y)$ is an isomorphism
for any $p\in M$ and non-zero $y\in Q^*_p$.
\end{defi}

For any transverse elliptic operator, we have the regularity estimate \cite[Theorem 12]{KLW}.

\section{Analysis of basic forms}
Throughout this section, let $(M,F)$ be a closed oriented manifold with a taut Riemannian foliation.   We give some basic tools for analysis on Sobolev space of basic sections. The goal of this section is to prepare for the analysis on the moduli space of the later sections. Firstly,
%\subsection{General analysis of basic sections}
we recall the  unique continuation property  on Hilbert space(see Kroheimer and Mrowka \cite[Chpater 7, 14]{KM}).

 \begin{lemma}[c.f. {\cite[Lemma 7.1.3]{KM}}]\label{lemma-7.1.3}
     Let $z:[t_1,t_2]\to H$ be a solution to the equation
     \[\frac{d}{dt}z(t)+L(t)z(t)=f(t),\]
     where $L(t)$ is a first-order transverse  elliptic formal self-adjoint operator, $H$ is a Hilbert space and $f(t)$ is an element of $ H$ with $\|f(t)\|\leq C\|z(t)\|$ for some constant $C$. If $z(t)$ is zero at one-point, then it vanishes identically.
 \end{lemma}

    On the  finite cylinder $Z=[a,b]\times M$, we have the following trace theorem.
\begin{thm}[c.f. {\cite[Theorem B10]{Wehrheim}}]\label{trace-thm}
  Set $Z=[a,b]\times M$ and $1\leq p<\infty$ and $n=codim(F)$. In the case $p<n+1$  we assume $1\leq q
  \leq\frac{(n+1)p-p}{1+n-p}$, and in the case $p\geq1+n$ we assume $1\leq q<\infty$. Then,  the basic trace theorem holds for all basic functions of $Z$, i.e. $L^p_1(Z)\to L^q(\partial Z)$.
\end{thm}

%\begin{pf} For all basic function $h\in C^\infty(Z)$, we have  \begin{eqnarray*}    \int_{\partial Z}|h|^qdvol_{\partial Z}    &=&\int_Zd(|h|^q\iota_\nu dvol_Z)\\    &=&\int_Z|h|^q L_\nu dvol_Z+\int_Z\frac{\partial}{\partial\nu}|h|^qdvol_Z\\    &\leq&(\int_Z|h|^q+\int_Z|h|^{q-1}|\nabla h|)\\leq&C(\|h\|^q_{L^q(Z)}+(VolZ)^{\frac1s}\|h\|^{q-1}_{L^\rho(Z)}\|\nabla h\|_{L^p(Z)}),  \end{eqnarray*}  where $C$ is a finite number as $Z$ is compact. We choose $1\leq r,s\leq\infty$ such that the Sobolev $L^p_1\to L^\rho$ holds for $\rho=r(q-1)$ and  \[\frac1s+\frac1r+\frac1p=1. \]  \begin{itemize}    \item In case $q=1$, set $s$ to be the number $p^*$ satisfy $\frac1p+\frac1s=1$.    \item In case $p>n$ and $q>1$, we choose $r=\infty$ and  $s=p^*$.    \item In case $p=n$ and $q>1$, we choose $r=\max(\frac1{q-1},p^*)$.    \item In case $p<n$ and $q>1$, we choose $r=\frac{np+p}{(n+1-p)(q-1)}$ to make the Sobolev embedding hold.  \end{itemize} Hence, in either case we can find a constant $C$, such that \[\|h\|_{L^q(\partial Z)}\leq C\|h\|_{L^p_1(Z)}. \]\end{pf}

By the same idea of the \cite[Appendix B]{Hormander}, we have the following theorem.

\begin{thm}
  For $j>\frac12$, there is a continuous restriction map between the Sobolev basic sections,
  \[r:L^2_j(Z,E)\to L^2_{j-\frac12}(\partial Z,E_0),\]
  where $E$ is the pull-back foliated bundle of $E_0\to M$.
\end{thm}

Consider the equation
\[\begin{cases}
  \Delta_bu=f&Z\\
  \frac{\partial u}{\partial \nu}=g&\partial Z,
\end{cases}\]
where $\nu$ denotes the unit normal vector field along the boundary $\partial Z$,
with the condition $$\int_Zf+\int_{\partial Z}g=0. $$
Recall that   we have that $\Delta_bu=\Delta u$ for any basic function $u$(see \cite[Page 86]{Ton2}). %  We can consider a more general case.
Similar to the property of the Laplacian equation with Neumann boundary condition,
 the following theorems hold.

\begin{thm}\label{Neumann-estimate-thm}%[c.f. {\cite[Theorem 2.3]{Wehrheim}}]%Let $P$ be an elliptic transverse operator of order $2$.
For any basic function $u$ on $Z$, it holds that   \[\|u\|_{L^2_{k+2}}\leq C(\|\Delta_bu\|_{L^2_k}+\|u|_{\partial Z}\|_{L^2_{k+1/2}}+\|u\|_{L^2_{k+1}}).\]  Furthermore, if on the boundary $u$ satisfies $\frac{\partial u}{\partial \nu}|_{\partial Z}=0$, then one has that   \[\|u\|_{L^2_{k+2}}\leq C(\|\Delta_b u\|_{L^2_k}+\|u\|_{L^2_{k+1}}).\]
\end{thm}

%Before the proving, we give a condition on when  $u$ is $L^2$-orthogonal to the kernel of $P$. Suppose that $\frac{\partial u}{\partial\nu}\big|_{\partial Z}=0$, then one deduces that
%\begin{eqnarray*}  \int_Z|du|^2dvol_Z&=&\int_Zdu\wedge\bar*du\chi_F\\  &=&\int_{Z}(u\Delta_b u)dvol_Z+\int_{\partial Z}u\frac{\partial u}{\partial\nu}dvol_{\partial Z},\end{eqnarray*}where we used the fact that $dvol_{\partial Z}=\iota(\nu)dvol_Z$ for the second equality. This implies that the kernel of $P$ under the Neumann condition consist of all constant functions. Hence the condition that $u$ is $L^2$-orthogonal to the kernel of $P$ is equivalent to the condition\[\int_Zu=0.\]

\begin{pf} Since $\Delta_b u=\Delta u$ for the basic function $u$,  the proof follows by the standard theory, c.f. {\cite[Formula 7.37 ]{Tay}} for the first formula and {\cite[Formula 7.34]{Tay}} for the second one.
   \end{pf}

\begin{thm}\label{CN-thm}
  Let $f\in L^2_k(Z)$ and $g\in L^2_{k+1/2}(\partial Z)$ be   basic functions. If the formula $\int_Zf+\int_{\partial Z}g=0$ holds, then there is a   solution to the equation
  \[\begin{cases}
    \Delta_bu=f& Z\\
    \frac{\partial u}{\partial \nu}=g&\partial Z.
  \end{cases}\].
\end{thm}

\begin{pf} The idea is similar to the one of \cite[Theorem 3.1]{Wehrheim}.
  We choose an element $v\in
  L^2_{2+k}(Z)$ such that $\frac{\partial v}{\partial\nu}\big|_{\partial Z}=g\big|_{\partial Z}$, this can be realized  by letting $v=\phi(t)g$ for some smooth function $\phi(t)$ with support near  the boundary and near each the boundary we have $\phi(t)=t-a$ for $t\in[a,a+\epsilon)$ and $\phi(t)=t-b$ for $t\in(b-\epsilon,b]$.
  Now we have
  \[\int_Z(f-\Delta_bv)=\int_Zf+\int_{\partial Z}\frac{\partial v}{\partial \nu}=
  0.\]
  Thus, by the Theorem \ref{Neumann-estimate-thm} and the Rellich embedding,  there exists a solution $u_1\in L^2_{k+2}(Z)$ to the Neumann problem with $f$
  replaced by $f-\Delta_b v$. The solution is given by $u=u_1+v$.
\end{pf}

%\subsection{Analysis on foliation  with taut bundle-like metric}

% In this subsection, we assume that  $(M,F)$ is a oriented closed manifold with   oriented Riemannian foliation  and $(M,F)$ admits a taut bundle-like metric.

 \begin{thm}\label{Weitzenbock-thm}
  On $Z=[a,b]\times M$, %we assume that the foliation of $M$ is taut.
  if the basic one-form $a\in\Omega^1_b(Z)$ satisfies the condition $a(\nu)=0$ on $\partial Z$, then the following holds
  \begin{equation}
    \int_Z(\nabla_{tr}a,\nabla_{tr}a)+Ric^T(a,a)=
    \int_Z(d_ba,d_ba)+(\delta_ba,\delta_ba).\label{Weitzenbock-formula}
  \end{equation}
\end{thm}
The proof  is similar to the classical formula, which is stated by the following proposition.

  \begin{prop}[Jung {\cite{Jung}}]
   % Let $Z$ be $(M,F)$ or for  a taut bundle-like Riemannian foliation.
   %Suppose that $\alpha\in \Omega^1_b(Z)$. Then, 
   The formula
    \begin{equation}
      \Delta_b\alpha=(\nabla^T)^*\nabla^T\alpha+Ric^T(\alpha), \label{formula-ricci-laplacian}
    \end{equation}
    holds, for any basic one-form $\alpha$  on $(M,F)$ and $(I\times M,F)$.
  \end{prop}
 Recall that  $Z=I\times M$,  where $I\subset (-\infty,\infty)$ is a compact interval, we establish the following lemmas.

\begin{lemma}\label{lemma-poincare-dual}
%Under the above conditions,
We have that
   $H^r_b(Z)$ is dual to $H^{m+1-r}_{b,c}(Z)$, where $m=co\dim(F)$ in $M$ and $H^{l}_{b,c}(Z)$ denotes the basic deRham cohomology which is vanishing at the ends.
 \end{lemma}

 \begin{pf}
  Let $e$ be a $1$-form on $I$ with integral $1$, which is vanishing at the ends. Define the map
 \[e_*:\Omega^{*}_b(M)\to \Omega^{*+1}_{b,c}(Z),\]
 by\[\alpha\mapsto \alpha\wedge e.\]
 We set a map $\pi_*:\Omega^*_{b,c}(Z)\to \Omega^{*-1}_b(M)$ as the integration along the $I$-direction. Similar to the arguments of \cite[Proposition 4.6]{BT},  the induced cohomology map $e_*: H^{*-1}_b(M)\to H^*_{b,c}(Z)$ is an isomorphism, whose inverse is the induced cohomology map of $\pi_*$. Since $H^r_b(M)$ is dual to $H^{m-r}_b(M)$, we proved the lemma.
% Let $\alpha\in \Omega^r_{b}(Z)$ with $d_b\alpha=0$ and $\beta\in \Omega^{m-r}_{b}$ with compact support and $d_b\beta=0.$  Consider $\alpha'=\alpha+d_bf$ and $\beta'=\beta+d_b f'$. Then  \begin{eqnarray*}    \alpha'\wedge\beta'\wedge\chi_F&=&    \alpha\wedge\beta\wedge\chi_F+\alpha\wedge d_bf'\wedge\chi_F+    d_bf\wedge\beta'\wedge\chi_F.  \end{eqnarray*}  The identity $d_b\alpha=d\alpha=0$ implies that  \[d(\alpha\wedge f'\wedge\chi_F)=(-1)^r\alpha\wedge d_bf'\wedge\chi_F.\]  Similarly, $d_b \beta=0$ implies that  \[d(f\wedge\beta'\wedge\chi_F)=d_bf\wedge\beta'\wedge\chi_F+  (-1)^{r-1}d_bf\wedge\beta\wedge\chi_F.\]  By the condition that $f'$ and $\beta$ are compactly supported, we have that \[\int_Z \alpha\wedge\beta\wedge\chi_F=\int_Z \alpha'\wedge\beta'\wedge\chi_F. \]  So the paring is well-defined.
  \end{pf}

%Similarly, we have the following lemma.

  \begin{lemma}
  Under the above conditions,
    one has  the isomorphism
    \[H^1_b(Z)\cong H^1_b(M).\]
  \end{lemma}
  \begin{pf}
    %The idea is   the same as for the usual deRham cohomology.
    It is clear that $\pi^*_Z :H^1_b(M)\to H^1_b(Z)$ is an injective, where $\pi_Z$ denotes the canonical projection $Z\to M$.
      We set the map $i:M\to Z$, as $p\mapsto (p,t_1)$, where $t_1$ is the left endpoint. Assume that there is an element $[\omega]\in H^1_b(Z)$ such that $i^*([\omega])=0$, i.e. $i^*\omega=d_bf$ for some basic function on $M$. Setting $\omega'=\omega-d_b\pi^*_Zf$, we rewrite it as
    \[\omega'=\alpha+f'dt.\]
    The condition $d\omega'=0$ implies that
    \[d\alpha=0,~\dot\alpha-df'=0. \]
    Hence, it holds that $\alpha|_{\{t_1\}\times M}=0$, and $\alpha(p,t)=\int^t_{t_1}df'(s)ds=d\int^t_{t_1}f'(s)ds$, i.e. $[\omega']=[\omega]=0$. This implies the isomorphism $H^1_b(Z)\cong H^1_b(M)$.
  \end{pf}

  \noindent
We finish this section with  the Coulomb gauge fixing property:
  For any non-trivial homotopy map $u:Z\to S^1$, by Theorem \ref{CN-thm}, we can find an element $v$ of this homotopy class satisfying the equation
  \begin{equation}
    \begin{cases}
    \delta_b(v^{-1}d_bv)=0&\mbox{ in } Z,\\
    dv(\nu)=0& \mbox{ on } \partial Z.
  \end{cases}\label{eqn-gauge-Neumann-condition}
  \end{equation}
  Notice that for  any basic one  form $\alpha$ we have that $\delta_b\alpha=\delta\alpha$.
  Let $\Gamma^1_b$ be the lattice of $H^1_b(M)\cap H^1(M,\mathbb Z)$. $\Gamma^1_b$ also corresponds to a lattice of $H^1_b(Z)\cap H^1(Z,\mathbb Z)$. Choose a basis $\{a_i\}$ of this lattice, by   pairing we have  a basis $\{\beta_i\}$  of $H^3_{b,c}(Z)$, which is
  dual to $\{a_i\}$.

%if $b^1_b>1$, then for general metric and perturbation the basic Seiberg-Witten moduli space is a set of finite points. Moreover, by choosing an orientation of the basic deRham cohomology, the signed counting of this moduli space is independent of the generic choice of the general metric and the perturbation.

\section{Basic Seiberg-Witten invariant on codimension $3$ foliation}
In this section, we define the basic Seiberg-Witten invariant on manifolds with a codimension $3$ foliation under a certain condition. %Since many results are similar to the model of codimension $4$, here we omit their proofs.

\subsection{Basic Seiberg-Witten equations on codimension $3$ foliation}
%In this section, we consider closed oriented $4$-manifold with vanishing Euler number.By the definition of foliation, i.e. $[X,Y]\in\Gamma(F)$ for $X,~Y\in\Gamma(F)$, we know that the existence of one-dimensional foliation is equivalent to say that the manifold admits a non-zero vector field.Therefore, for manifold with trivial Euler number, we know that there is a one-dimensional foliation. We assume that $H^1_b(M)=0$ with respect to this foliation.Let $T$ be a non-vanishing vector field of $M$. We have a one-parameter group $\exp(tT)$ without fixing-point of $Diff(M)$, the closure of this group is a torus. For a metric $g$, by integrating $g$ over this torus, we get a invariant metric over the action this torus, especially over the action of this one-parameter group.Hence, we get a bundle-like metric, and by the argument in the formal section, we can make this metric taut.
In this subsection, we focus on the manifolds with a foliation satisfying the following assumption.

\begin{assump}\label{assum-main}
  Let $(M,F)$ be an oriented closed manifold with codimension 3 oriented Riemannian foliation $F$ and admits a transverse \spinc structure $\mathfrak s$. Suppose  that  $H^1_b(M)\cap H^1(M,\mathbb Z)\subset H^1(M)$ is a lattice of $H^1_b(M)$.
\end{assump}

%By the above lemma, we have that\[\int_M[u]\wedge( F_A-\bar*\eta)\wedge \chi_F\]is unique up to some positive-number multiplication.
Let $\mathcal A_b(\mathfrak s)$ be the space of basic \spinc connections.  We define the basic Seiberg-Witten equations for manifolds with a codimension-$3$ Riemannian foliation by
\begin{equation}
  \begin{cases}
  \Dirac_{b,A}\Psi=0,\\
 \frac12\bar*F_{A^t}-q(\Psi)=0,
\end{cases}\label{eqn-SW-equation-1}
\end{equation} for $(A,\Psi)\in \mathcal A_b(\mathfrak s)\times \Gamma_b(S)$.
Here,   we identify the    traceless endomorphism of the spinor bundle with the imaginary valued cotangent bundle(we use $q(\Psi)$ instead of $\rho^{-1}(\Psi\Psi^*)_0$ in the book \cite[Formula 4.4]{KM}), $A^t$ denotes the connection on the determinate bundle of $S$, see \cite[Notation 1.2.1]{KM} and $\Dirac_{b,A}$ denotes the basic Dirac operator twisted with the basic connection $A$. The basic gauge group \[\mathcal G_b=\{u:M\to U(1)\big|~L_Xu\equiv0,\mbox{ for all }X\in \Gamma(F)\}, \] acts on $\mathcal C_b(\mfs)= \mathcal A_b(\mathfrak s)\times \Gamma_b(S)$ as (see \cite[Formula 4.5]{KM}):
\[u:(A,\Psi)\mapsto (A-u^{-1}du, u\Psi).\]
%Here we  the notation $q(\Psi)$ instead of $\rho^{-1}(\Psi \Psi^*)_0$ in the book \cite[Formula 4.4]{KM}.
To construct the basic Seiberg-Witten invariant, we need to consider the moduli space, which is defined as below.

 \begin{defi}
   The moduli space $\mathcal M_g(M,F,\mfs)$ of the basic Seiberg-Witten equations on $(M,F,g,\mathfrak s)$ is the space of the solutions to the above Seiberg-Witten equations modulo the gauge transformation group $\mathcal G_b$. The moduli space $\mathcal M^*_g(M,F,\mfs)$ is the irreducible part of $\mathcal M_g(M,F,\mfs)$, i.e. the spinor field part is not identically zero.
 \end{defi}

Similar to the manifold case, we consider the following complex
 \[L^2_{2+k}(\Omega^0_{b}(M,i\mathbb R))\overset{G_{(A,\Psi)}}\longrightarrow
 L^2_{1+k}(\Omega^1_{b}(Y,i\mathbb R))\oplus L^2_{1+k}(\Gamma_{b}(S))\overset{L_{(A,\Psi)}}\longrightarrow L^2_{k}(\Omega^1_{b}(M,i\mathbb R))\oplus L^2_k(\Gamma_{b}(S)),\]
  where $G_{(A,\Psi)}f=(-d_bf,f\Psi)$ and $L_{(A,\Psi)}(a,\Phi)=(-\frac12\bar*da+q(\Psi,\Phi),\Dirac_{b,A}\Phi+a\Psi)$.
  %By the straightforward calculation,
  One has that \[L\comp G=0.\]
  To read off the virtual dimension of the moduli space,
 it is convenient to consider the  form operator(see \cite[Formula 2.7]{B.L. Wang}),
 \[Q_{(A,\Psi)}:L^2_{2+k}(\Omega^1_{b}(M,i\mathbb R))\oplus L^2_{2+k}(\Omega^0_{b}(M,i\mathbb R))\oplus L^2_{2+k}(\Gamma_{b}(S))\]\[\to
 L^2_{1+k}(\Omega^1_{b}(M,i\mathbb R))\oplus L^2_{1+k}(\Omega^0_{b}(M,i\mathbb R))\oplus L^2_{1+k}(\Gamma_{b}(S)),\]
 \begin{equation}
   Q_{(A,\Psi)}=\left(\begin{array}{cc}
   L_{(A,\Psi)}&G_{(A,\Psi)}\\
   G^*_{(A,\Psi)}&0\\
 \end{array}\right), \label{Hessian-operator}
 \end{equation}
 where $G^*_{(A,\Psi)}$ is the formal self-adjoint of $G_{(A,\Psi)}$.% defined by
 %\[G^*_{(A,\psi)}(a,\phi)=-2(\delta_ba+i Im\langle\psi,\phi\rangle).\]
 %We can this operator $Q_{(A,\psi)}$ Hessian operator.
 To show the smoothness of the moduli space, we need to perturb the above equations:
 Fixing a  basic perturbation $\eta\in i\Omega^1_b(M)$, we denote the moduli space of the perturbed basic Seiberg-Witten equations by $\mathcal M_{g,\eta}(M,F,\mfs)$, i.e.   the space of the solutions to equations
\[\begin{cases}
  \Dirac_{b,A}\Psi=0,\\
  \frac12\bar*F_{A^t}-q(\Psi)=\eta.
\end{cases}\]
modulo the gauge action.

\begin{defi}  We say that $[A,\Psi]\in \mathcal M_{g,\eta}(M,F,\mfs)$ is \emph{non-degenerate} if  \[\ker(L_{(A,\Psi)})/ImG_{(A,\Psi)}=0.\]\end{defi}

%We recall the result of Seiberg-Witten equation for  codimension $4$ foliation.

\begin{prop}
If $H^1_b(M)\cap  H^1(M,\mathbb Z)$ is a lattice of $H^1_b(M)$, then for a generic perturbation
  the irreducible moduli space of basic Seiberg-Witten equations is a compact manifold with formal dimension zero.
\end{prop}
\begin{pf}
To show that moduli space is a compact and smooth manifold, we can repeat the similar arguments of  \cite[Lemma 2.2.3, Lemma 2.2.6, Theorem 2.2.8]{B.L. Wang}.  % \ref{thm-sw-codim4}.
To prove that the formal dimension of the moduli space is zero, we   calculate the index the  operator \eqref{Hessian-operator}.
% which is zero, i.e. formal dimension of the moduli space is zero.
 %By the Weitzenb\"ock formula, we get\begin{eqnarray*}   \Dirac_{b,A}\Dirac_{b,A}\psi=\nabla^*_A\nabla_A\psi+(\frac14Scal^T)\cdot\psi+\frac14|\psi|^4, \end{eqnarray*} where $Scal^T$ denotes the scalar curvature of the transverse connection. By the Kato-inequality \begin{eqnarray*}   \Delta_b|\psi|^2&=&-\sum^q_{i=1}e^2_i(\psi,\psi)+\nabla^T_{e_i}e_i   (\psi,\psi)\footnote{In general case, i.e. without the taut condition, there   willbe an extrac term $\tau_b(\psi,\psi)$. By our assumption, it is vanishing.}\\   &=&\sum^q_{i=1}-2(\nabla^A_{e_i}\psi,\psi)   -(\nabla^A_{e_i}\nabla^A_{e_i}\psi,\psi)   -(\psi,\nabla^A_{e_i}\nabla^A_{e_i}\psi)+   (\nabla^A_{\nabla^T_{e_i}e_i}\psi,\psi)+   (\psi,\nabla^A_{\nabla^T_{e_i}e_i}\psi)\\   &=&-2|\nabla^A\psi|^2+2Re\left((\nabla^A)^*\nabla^A\psi,\psi\right), \end{eqnarray*} we get \[\frac12\Delta_b|\psi|^2 \leq(-\frac14Scal^T+\frac12|\eta|)|\psi|^2-\frac18|\psi|^4,\] therefore we get \[|\psi|^2\leq\max_M(0,-Scal^T). \] By the standard argument of the Seiberg-Witten theory \cite[Chapter 3]{Morgan}, we know that thethis moduli space is compact.
 Recall that the  operator \eqref{Hessian-operator} is a compact perturbation of the transverse elliptic operator
 \[\left(\begin{array}{ccc}
   \frac12\bar*d_b&-d_b&0\\
   -\delta_b&0&0\\
   0&0&\Dirac_{b,A}\\
 \end{array}\right).\]
 %which is a first-order transverse elliptic operator,
 However it is not formal self-adjoint in general.
 It is clear that the operator
  \[\left(\begin{array}{ccc}
   \frac12(\bar*d_b-\frac12\kappa)&-d_b&0\\
   -\delta_b&0&0\\
   0&0&\Dirac_{b,A}\\
 \end{array}\right)\]
 is a formal self-adjoint operator and the difference
 \[\left(\begin{array}{ccc}
   \frac12\bar*d_b&-d_b&0\\
   -\delta_b&0&0\\
   0&0&\Dirac_{b,A}\\
 \end{array}\right)-\left(\begin{array}{ccc}
   \frac12(\bar*d_b-\frac12\kappa)&-d_b&0\\
   -\delta_b&0&0\\
   0&0&\Dirac_{b,A}\\
 \end{array}\right)\]
 is a compact operator, hence they have the same index zero. To equip an orientation of the moduli space, we just need to equip an orientation of the determinant line bundle of \eqref{Hessian-operator}.
 %$$\left(\begin{array}{cc}  \bar*d_b&-d_b\\  -\delta_b&0\\ \end{array}\right):\Omega^1_b(M)\oplus \Omega^0_b(M)\to \Omega^1_b(M)\oplus \Omega^0_b(M)$$ one have an orientation of the moduli space.
\end{pf}

The above proposition implies that  the determine line bundle of the operator \eqref{Hessian-operator} is trivial over the moduli space $\mathcal M^*_{g,\eta}(M,F,\mfs)$. Hence, there is a natural orientation for the moduli space.   The basic  Seiberg-Witten invariant $SW_{g,\eta}(M,F,\mfs)$ on $(M,F)$ is defined by the signed counting of the moduli space $\mathcal M^*_{g,\eta}(M,F,\mfs)$.%, which is stated as below.

%\begin{defi}  Suppose that $(M,F)$ satisfies Assumption \ref{assum-main}. We fix a transverse \spinc structure $\mathfrak s$. Choosing a bundle-like metric $g$, a perturbation $\eta$ and an orientation of $H^1_b(M)$ and $\ker(d_\kappa+\delta_b)$, we define the basic Seiberg-Witten invariant as the signed counting the irreducible solutions. We denote it by $SW^*_{\eta,g}(\mathfrak s)(M)$\end{defi}

\subsection{Basic Seiberg-Witten invariant on codimension $3$ foliation}

In this subsection, we   show how the basic Seiberg-Witten invariant depends on  the bundle-like metrics and  basic perturbations.
We review the notion of the reducible solution.
Let $(A,\Psi)$ be a solution to the basic Seiberg-Witten equations \eqref{eqn-SW-equation-1}. When $\Psi=0$,  the basic Seiberg-Witten equations reduce to a single equation, $\frac12\bar *F_{A^t}=\eta$. If $\eta=0$, then the reducible class  is identified with the moduli space of the flat basic $U(1)$-connection of  $\det(\mathfrak s)$. We denote the first Chern class of $\det(\mathfrak s)$ by $c_1(\mathfrak s)$.

\begin{lemma}
  The equation $\frac12\bar *F_{A^t}=\eta$ has  a solution  if and only if  $d_b\bar*\eta=0$ and $\pi i[c_1(\mathfrak s)]=[\bar *\eta]$ in $iH^2_b(M)$. In particular, if $H^1_b(M)\cap H^1(M,\mathbb Z)$ is a lattice of $H^1_b(M)$, then the set of  reducible solutions modulo gauge action  is isomorphic to $H^1_b(M)/(H^1_b(M)\cap H^1(M,\mathbb Z))$.
\end{lemma}
\begin{pf}
  Obviously, $[\frac12F_{A^t}]=[\bar *\eta]$ is a necessary condition to solving the equation.
  Conversely, suppose that $\pi i[c_1(\mathfrak s)]=[\bar *\eta]$.   We fix a basic connection $A_0$ such that $[\frac12 F_{A^t_0}]=[\bar *\eta]$. It suffices to get a basic one-form $a$ such that $da=\bar *\eta-\frac12 F_{A^t_0}$. Since $\bar *\eta-\frac12 F_{A^t_0}$ is an exact basic form, such   one-form $a$ always exists. By choosing one solution $a_0$ of the above equation, one  represents all the others as
  \[a_0+\mbox{ closed one form}.\]
  Any two solutions $a_1$ and $a_2$ are equivalent, if and only if $a_1=a_2+u^{-1}du$ for some $u\in\mathcal G_b$.
\end{pf}

\begin{thm}
 Let $(M,F)$   satisfy the  Assumption \ref{assum-main}.
  When $b^2_b>1$,  then the basic Seiberg-Witten invariant is well-defined, i.e. it is independent of  the generic choice of the basic perturbations and  bundle-like metrics.
\end{thm}
\begin{pf}
 To show that the  basic Seiberg-Witten invariant is independent of the basic perturbations and  bundle-like metrics, we apply the similar proof in codimension $4$ case \cite[Chapter 5]{Morgan}. Here we give a sketch of  the proof.
  We denote by $\mathcal N=\{a\in i\Omega^1_b|~d_b\bar*a=0\}$, which is identified with the space of closed basic two-forms. Set
\[\mathcal W_{\mathfrak s}=\mathcal W_{\mathfrak s}(g)=\{\eta\in \mathcal N|~[\bar*\eta]=\pi ic_1(\mathfrak s)\}.\]
Notice that $\mathcal W_{\mathfrak s}$ is a codimension $b^2_b$ affine subspace of $\mathcal N$.
When $b^2_b>1$, $\mcWs$ is of codimension two or more. %For two perturbations $\eta_{i}\in \mathcal N\setminus \mcWs$ for $i=1,~2$, there are no reducible solutions. Because of $b^2_b>1$,
One can choose a generic path $\eta_s$ connecting these two perturbations $\eta_1$ and $\eta_2$, such that for each $s\in[1,2]$, $\eta_{s}\in \mathcal N\setminus \mcWs$. This completes the proof.
\end{pf}

%\subsection{Seiberg Witten invariant for foliation with $b^2_b=1$}
 At the end of this subsection, we   show the dependence on the basic perturbations and bundle-like metrics for the basic Seiberg-Witten invariant, when $b^2_b=1$.
Choose an orientation of the one dimensional  space $iH^2_b(M)$. There exists an unique unit $g$-harmonic basic two form $\omega$, i.e. $\|\omega\|_{L^2}=1$.
The wall $\mcWs$ is identified with  the solutions to the linear equation
\[(\bar*\eta-\pi ic_1(\mfs),\omega)=0.\]
We set $\mathcal N^\pm$ as $\mathcal N^\pm=\{\eta\in\mathcal N|~\pm(\bar*\eta-\pi i c_1(\mfs),\omega)>0\}$.
%We consider the codimension-$1$ subspace\[R=\{(g,\eta)\big|~ \eta=2\pi\bar *_g([c_1(\mfs)])+\bar*g d_b\nu\}.\]
%There is a projection map $\pi' :Met\times \mathcal N\to H^2_b(M)$, defined by\[(g,\eta)\mapsto(\bar*_g\eta)_h.\]
Consider a family of $(g_t,\eta_t)_{t\in[-1,1]}$ crossing the wall transversely once at $t=0$, such that $g_t$ is locally constant near $-1,~0$ and $1$. The one-parameter family moduli space
\[\mathcal M_{[-1,1]}=\bigcup_{t\in[-1,1]}\mathcal M_{g_t,\eta_t}(\mfs)\]
has   reducible solutions at $t=0$.
Under the Assumption \ref{assum-main}, the reducible space %$\mathcal M^0$,
$\mathcal M^0=\{(A,0)\big|~\frac12F_A=\bar*\eta^0+d_b\alpha\}/\mathcal G_b$
is identified with the  torus $T^{b^1_b}$, where $b^1_b\geq1$.
We decompose the connection $A$ as $A=A_0+\alpha+a$, where $A_0$ is a connection such that $\frac12F_{A^t_0}=\bar*\eta^0$ and $a$ is a harmonic one-form.

 \begin{lemma}
  Suppose that $\mathcal M^0$ is parameterized by $a\in H^1_b(M)/H^1_b(M)\cap H^1(M,\mathbb Z)$ with the decomposition $A=A_0+\alpha+a$. If $\ker(\Dirac_{b,A})=0$, then there is no irreducible solution connection $A$ in
  $\mathcal M_{[-1,1]}$.
\end{lemma}
\begin{pf}
The idea is similar to \cite[Lemma2.3.7]{B.L. Wang}.
  Let $(A_t,\Psi_t)$ be a family of solution to the perturbed basic Seiberg-Witten equations with  $(A_t,\Psi_t)|_{t=0}=(A,0)$. Near $(A,0)$ we write $(A_t,\Psi_t)=(A+a_t,\Psi_t)$. We let $(a_t,\Psi_t)$ satisfy
  \[\begin{cases}
    \delta_ba_t=0\\
    \bar*d_ba_t=q(\Psi(t))\\
    \Dirac_{b,A}(\Psi_t)+a_t\Psi_t=0\\
  \end{cases}\]
  near $t=0$. Locally, we write $\Psi_t=\sum_{i\geq1}t^i\Psi_i$ and $a_t=\sum_{i\geq1}t^ia_i$. Differentiating the third equation with respect to $t$, gives  the result.
\end{pf}

\begin{defi}
  A reducible solution is called regular, if $\ker\Dirac_{b,A}=0$.
\end{defi}
Therefore, by the above lemma
we can perturb the equation such that regular part $\mcM^0_{reg}$ of the reducible solutions is isolated from the irreducible solution in $\mathcal M_{[-1,1]}$.

In particular,   when the foliation is taut($H^1_b(M)\cong H^2_b(M)$), we can find a generic perturbation %as in\cite[Lemma 2.2.5]{B.L. Wang}
such that there are only finitely many points in $\mcM^{0}$  meeting the irreducible solutions in $\mcM_{[-1,1]}$ and the kernel of the associated twisted Dirac operator $\Dirac_{b,\omega+a_\theta}$ is of dimension $1$. Following the same arguments as in \cite[Proposition 2.3.8]{B.L. Wang}, yields the following proposition.

\begin{prop}
 Let $(M,F)$  satisfy the  Assumption \ref{assum-main}, and $F$ is taut.
  Suppose that $(g_{-1},\eta_{-1})$ and $(g_{1},\eta_{1})$ belong to the different parts which are separated by  $\ ofmWs$, choosing an orientation of $H^1_b(M)$. Then, we have the formula
   $$SW_{g_1,\eta_1} (M,F,\mathfrak s )-SW_{g_{-1},\eta_{-1}}(M,F,\mathfrak s)=
  SF(\Dirac_{A(\theta)}),$$
  where $A(\theta)$ defines  a connection  joining $A_{-1}, ~A_1$ and the $SF(\Dirac_{A(\theta)})$ denotes the spectral flow of the corresponding basic Dirac operator.
\end{prop}
%The proof is no different to \cite[Proposition 2.3.8]{B.L. Wang}, here we omit it.

\begin{pf}
 We choose a family of metrics and perturbations $(g_t,\eta_t),t\in[-1,1]$ such that it crosses the set $\mathfrak R=\{(g_Q,\eta)| \eta=\eta^0+\bar*_{g_Q}d\alpha,~[\eta^0]=\pi i[c_1(\mfs)],~\eta\mbox{ is harmonic }\}$ at $t=0$ with finite singular points in $\mathcal M^0$. Similar to \cite[Proposition 2.3.8]{B.L. Wang},   the difference $SW_{g_1,\eta_1} (M,F,\mathfrak s )-SW_{g_2,\eta_2}(M,F,\mathfrak s)$ is equal to the spectral flow of the twisted basic Dirac operator along $\mathcal M^0$, which proves the proposition.
\end{pf}

%When $b^1_b(M)=0$, we choose a taut metric, by the argument of \cite[Proposition 2.2.15]{B.L. Wang}, the wall\[\mathcal N=\{(g,\bar*da)\big| \ker(\Dirac^a_g)\neq0\}\]is a codimension one subset of $Met\times Z^1_b(M)$, where $Z^1_b(M)$ denotes the co-closed basic one-forms. We have the following proposition.

%\begin{prop}  We can choose a generic path $(g_t,\bar*da_t)\in\ Met \times  Z^1_b(M)$ for $t\in[-1,1]$ crossing the wall $\mathcal N$ transversely.  Moreover, we have  \[SW(M,\mathfrak s,g_1,\eta_1)-SW(M,\mathfrak s, g_{-1},\eta_{-1})  =SF(\Dirac^{a_t}_{g_t}),\]  where $SF(\Dirac^{a_t}_{g_t})$ is the spectral flow of the Dirac operator  $\Dirac^{a_t}_{g_t}$ for $t\in[-1,1]$.\end{prop}

% By the similar argument of Lim's work \cite{Lim},
Summarizing the above arguments,  yields  the following results.
 \begin{thm}
 Suppose that $(M,F)$ satisfies the Assumption \ref{assum-main}. Then  for  each  transverse \spinc structure $\mathfrak s$ and for a generic perturbation $\eta$ and a bundle-like metric $g$, we define the basic Seiberg-Witten invariant $SW_{\eta,g}(M,F,\mathfrak s)$ by the signed counting of the moduli space $\mathcal M^*_{g,\eta}(M,F,\mfs)$.
  Moreover, we have the properties:
   \begin{itemize}
  \item If $b^2_b>1$, then for a generic bundle-like metric and perturbation the basic Seiberg-Witten moduli space is  a smooth compact manifold, and  $SW_{g,\eta} (M,F,\mathfrak s )$  is generically independent of the choice of bundle-like metrics and  perturbations.
  \item If $b^2_b=1$, then $SW_{\eta,g} (M,F,\mathfrak s )$ depends only on the component of $H^2_b(M)\setminus \pi c_1(\mfs)$.
  %\item If $b^1_b=0$($b^2_b=0$), then for generic perturbation the only reducible solution is isolated in $C^\infty$-topology.
\end{itemize}

 \end{thm}

\section{Chern-Simons-Dirac functional and moduli space for foliation}
The purpose of this section is to give the preparation to show the compactness for the moduli space. %We will define the foliated Chern-Simons-Dirac functional in the first subsection and give some estimates for this functional in the second subsection.   The original idea was introduced by Kroheimer and Mrowka(see \cite[Chpater 5]{KM}).
\subsection{Chern-Simons-Dirac functional for foliation}
Throughout this subsection, the following assumption holds for  $(M,F)$.
\begin{assump}\label{assum-main-1}
  Let $(M,F)$ be an oriented closed manifold with codimension 3 oriented Riemannian foliation $F$ and admits a transverse \spinc structure $\mathfrak s$. Suppose  that  $H^1_b(M)\cap H^1(M,\mathbb Z)\subset H^1(M)$ is a lattice of $H^1_b(M)$ and $F$ is taut.
\end{assump}
Fixing a bundle-like metric and a transverse \spinc structure $\mathfrak s$,
we define the basic  Chern-Simons-Dirac functional over $M$ by
\[\L(A,\Psi)=-\frac18\int_M(A^t-A^t_0)\wedge(F_{A^t}+F_{A^t_0})\wedge\chi_F
+ \frac12\int_M(\Psi,\Dirac_A\Psi)dvol_M,\]
for any $(A,\Psi)\in \mathcal A_b(\mathfrak s)\times\Gamma_b(S)$.

The formal gradient of the Chern-Simons-Dirac functional is given by the following lemma.
\begin{lemma}
% We have the following identity
It holds that
  \[grad\L(A,\Psi)=(\frac12 \bar*(F_{A^t}+\frac12(A^t-A^t_0)\wedge\kappa_b)-
  q(\Psi)),\Dirac_{b,A}\Psi). \]
\end{lemma}
Note that the gradient is not gauge-invariant in general.

\begin{pf}
  Choosing a variation $(A+ta,\Psi+t\Phi)\in \mathcal A_b(\mathfrak s)\times\Gamma_b(S)$,  % by the straightforward calculation
  one deduces that
  \begin{eqnarray*}
    &&\partial_t(\L(A+ta,\Psi+t\Phi))\big|_{t=0}\\
    &=&-\frac18\int_M(2a\wedge(F_{A^t}+F_{A^t_0})+(A^t-A^t_0)\wedge 2da)\wedge
    \chi_F
    \\
    &&+ \frac12\int_M(\langle\Phi,\Dirac_A\Psi\rangle+
    \langle\Dirac_A\Psi,\Phi\rangle+\langle\Psi,a\cdot\Psi\rangle)dvol_Q\wedge\chi_F\\
    &=&\int_M a\wedge(-\frac12F_{A^t}-\frac14(A^t-A^t_0)\wedge\kappa_b+
    \bar*q(\Psi))\wedge\chi_F
    +\int_M Re\langle\Dirac_A\Psi,\Phi\rangle dvol_Q\wedge \chi_F\\
    &=&\int_M(a,\frac12\bar*(F_{A^t}+\frac12(A^t-A^t_0)\wedge\kappa_b)-q(\Psi))dvol_M+
    \int_M Re\langle\Dirac_A\Psi,\Phi\rangle dvol_M,
  \end{eqnarray*}
    where we used Rummler formula \eqref{formula-Rummler} to deduce the second identity.
\end{pf}

We say $(A,\Psi)$ is a critical point of $\L$, if its gradient vanishes at $(A,\Psi)$.
\noindent
For any gauge action $u\in\mathcal G_b(M,S^1)$, we have
\[\L(u(A,\Psi))-\L(A,\Psi)=\frac12\int_Mu^{-1}du\wedge (F_{A^t})\wedge \chi_F.\]
In general, the above term  does  depend on the choice of the representation of the cohomology class $[u]=[\frac{-i}{2\pi}u^{-1}du]$.
When $F$ is taut, the critical points of the basic Chern-Simons-Dirac-functional coincide with the solutions of the basic Seiberg-Witten equations \eqref{eqn-SW-equation-1}.
% When $(M,F)$ admits a  taut bundle-like metric,
We consider the gradient flow to the
Chern-Simons-Dirac functional
\[\frac{d}{dt}(A(t),\Psi(t))=-grad\L(A(t),\Psi(t))\]
for a path $(A(t),\Psi(t))$ of configuration space $\mathcal C_b(M,F,\mfs)=\mathcal A_b(M)\times\Gamma_b(S)$.
Let $\mathcal M([\alpha],[\beta])$ denote the moduli space of trajectories connecting the critical points up to gauge, i.e. the solutions to the basic Seiberg-Witten equations
\[\begin{cases}
  \frac12F^+_{A^t}=q(\Psi),\\
  \Dirac^+_A\Psi=0,
\end{cases}\]
on $\mathbb R\times M$ modulo the gauge action, such that $[A(t),\Psi(t)]\to[\alpha]$ as $t\to-\infty$ and $[A(t),\Psi(t)]\to[\beta]$ as $t\to\infty$,
 where $[\alpha],~[\beta]$ denote the gauge equivalence classes of the critical points. The components of this space have different dimensions corresponding to the different lifts of $[\alpha],~[\beta]$. This is a manifestation of the fact that the quotient space $\mathcal B_b(M,F,\mfs)=\mathcal C_b(M,F,\mfs)/\mathcal G_b$ may have non-trivial fundamental group.
We have a decomposition
\[\mathcal M([\alpha],[\beta])=\bigcup_{z\in\pi_1([\alpha],[\beta])}
\mathcal M_{z}([\alpha],[\beta])\]
as the union over the moduli spaces in a given relative homotopy class, where $\pi_1([\alpha],[\beta])$ denotes all the homotopy classes of paths joining $[\alpha]$ and $[\beta]$.
 Given two critical points  $\alpha,~\beta\in\mathcal C(M,F,\mfs)$, we define  the quantity
 \[gr(\alpha,\beta)\in\mathbb Z\]
 by the spectral flow of the Hessian operator \eqref{Hessian-operator} of a path  connecting them. This is well defined because the spectral flow is invariant under homotopy and the configuration space $\mathcal C(M,F,\mfs)$ is simply connected. Moreover, such a number computes the formal dimension of the moduli space of trajectories connecting $\alpha$ and $\beta$ by the path $z$, i.e. $\dim\mathcal M_{z}([\alpha],[\beta])$.

\subsection{Compactification  of the moduli space}

In this subsection, we give a   compactification  of the moduli space on the cylinder. The original idea is using the energy functional, which was introduced by Kroheimer and Mrowka \cite[Chapter 5]{KM}. Here we gave a foliated version of their work.
 Recall that $Z=I\times M$,   $I\subset (-\infty,\infty)$, and $\alpha,~\beta$ are the critical points of the $\L $. Denoting by $\mathcal M([\alpha],[\beta])$ the moduli spaces of trajectories in $\mathcal B_b(M,F,\mfs)$  and $\check{\mathcal M} ([\alpha],[\beta])$ the unparameterized moduli space. We  show that by adding broken trajectories it can be compactified, which is denoted by $\check{\mathcal M}^+ ([\alpha],[\beta])$.

   \begin{lemma}
     Let $(Z=I\times M,F)$ be a compact taut Riemannian foliation with a taut bundle-like product metric $g$. For a basic one-form $\alpha $ on $Z$, satisfying the boundary condition $(\alpha,\nu)=0$, where $\nu$ denotes the outward unit vector field. If there exists a constant $C_0$ such that $\int_Z\beta_i\wedge \alpha\wedge\chi_F\in[ -C_0, C_0]$ for each $\{\beta_i\}$, where $\{\beta_i\}$ is a basis of $H^3_{b,c}(Z)$. Then, there are constants $C_1$ and $C_2$ such  that
      $$\|\alpha\|^2_{L^2_1(Z)}\leq C_1\int_Z(|\delta_b\alpha|^2+|d_b\alpha|^2)+C_2.$$
   \end{lemma}
   \begin{pf}
Recall that by Theorem \ref{Weitzenbock-thm}, we have
\begin{equation}
  \int_Z|\nabla^T\alpha|+Ric^T(\alpha,\alpha)=\int_Z|\delta_b\alpha|^2+|
  d_b\alpha|^2. \label{formula-basic-ric-laplacian-1}
\end{equation}
Because of the product metric on this finite cylinder,  $Ric^T$ has a uniform bound. Moreover, by the $L^2_1$-bound of $\alpha$, the estimate follows from the proof of \cite[Lemma 5.1.2]{KM}.
   \end{pf}

 Recall that $(M,F)$ is a closed oriented taut Riemannian foliation with codimension 3, and it admits a transverse   \spinc structure. For the spinor bundle $S_Z=S^+\oplus S^-$ on $Z=I\times M$, we take $S_Z=S\oplus S$. For the Clifford multiplication $\rho_Z:TZ\to Hom(S_Z,S_Z)$, we take
 \[\rho_Z(\partial_t)=\left(\begin{array}{cc}
   0&-1\\
   1&0\\
 \end{array}\right),~ \rho_Z(v)=\left(\begin{array}{cc}
   0&-\rho(v)^*\\
   \rho(v)&0\\
 \end{array}\right),\]
 for $v\in Q$. A time-dependent \spinc connection $B$ on $S$ gives a \spinc connection $A$ on $S_Z$,
whose $t$ component is an ordinary differentiation, i.e.
\[\nabla_A=\frac{d}{dt}+\nabla_B.\]
We call the connection  $A$ in temporal gauge.
  With a basic connection $A$ as above, we have the basic Dirac operator
\[\Dirac^+_{b,A}:\Gamma_b(S^+)\to \Gamma_b(S^-),~\Dirac^+_{b,A}=\frac{d}{dt}+\Dirac_B.\]
 For a general basic connection $A$, we write
\[A=B+(cdt),\]
where $c$ is a basic function. The corresponding basic Seiberg-Witten equations for $(A,\Psi)$ are written as
\begin{equation}
\begin{cases}
  \frac{d}{dt}B-dc=-(\frac12\bar*F_{B^t}-q(\Psi)),\\
  \frac{d}{dt}\Psi+c\Psi=-\Dirac_B\Psi.
\end{cases}
\end{equation}
We define the analytic energy by
\begin{equation}
    \mathcal E^a(A,\Psi)=2(\L(t_1)-\L(t_2))+\int_Z|\dot\gamma(t)+
    grad(\L)(A,\Psi)|^2,\label{formula-analytic energy}
  \end{equation}
  where $\gamma=(A,\Psi)$.
  By the transverse Weitzenb\"ock formula \cite{GK} and similar arguments of \cite[Section 4.5]{KM}, one has  that
\begin{eqnarray*}
  \mathcal E^a(A,\Psi)%&=&  \int_Z|\nabla_t\Phi|^2+|\Dirac_B\Phi|^2+|\dot{B}-dc|^2+\int_Z|\bar*F_{B^t}-q(\Phi)|^2\\
  &=&\int_Z(|\nabla_A\Psi|^2+ \frac{Scal^T}{4}|\Psi^2|+\frac14|\Psi|^4)+\frac14\int_Z|F_{A^t}|^2.
  %+\int_Z(|\dot{B}-dc|^2+|F_B|^2.
\end{eqnarray*}

The topological energy $\mathcal E^{top}$ is defined by the twice of  drop of the basic Chern-Simons-Dirac  functional on the cylinder, i.e.
\[\mathcal E^{top}(A,\Psi)=2(\L(t_1)-\L(t_2)).\]

%\begin{prop}  We have that  \begin{equation}    \mathcal E^a(A,\Phi)=2(\L(t_1)-\L(t_2))+\int_Z|\dot\gamma(t)+    grad(\L)(A,\Phi)|^2.\label{formula-analytic energy}  \end{equation}\end{prop}
%\begin{pf}By the identity\[[\nabla_t,\Dirac_B]=\dot{B}-dc,\]we have that  \begin{eqnarray*}    &&\mathcal E^a(A,\Phi)-\int_Z|\dot\gamma(t)+    grad(\L)(A,\Phi)|^2\\    &=&2\int_Z\langle\nabla_t\Phi,\Dirac_B\Phi\rangle+    \int_Z(\dot{B}-dc,\bar*F_{B^t}+q(\Phi))\\    &=&\int_Z2\langle\nabla_t\Phi,\Dirac_B\Phi\rangle+\langle\Phi,\dot{B}-    dc\Phi\rangle+2\int_Z(\dot{B}-dc)\wedge F_{B^t}\wedge\chi_F\\    &=&\int_Z(2\langle\nabla_t\Phi,\Dirac_B\Phi\rangle+\langle\Phi,    [\nabla_t,\Dirac_B]\Phi\rangle)+2\int_Z\dot{B}\wedge F_B\wedge\chi_F\\    &=&\int_Z(\langle\nabla_t\Phi,\Dirac_B\Phi\rangle+\langle \Phi,\nabla_t    \Dirac_B\Phi\rangle)    +\int^{t_2}_{t_1}\frac{d}{dt}\int_M(B-B_0)\wedge(F_B+F_{B_0})\wedge\chi_F\\    &=&2(\L(t_1)-\L(t_2)).  \end{eqnarray*}\end{pf}

For convenience, we denote the Seiberg-Witten map by $\mathfrak F$.% By the idea of {\cite[Theorem 5.1.1, Corollary 5.1.8]{KM}}, we have the following theorem.

\begin{thm}\label{thm-coro5.1.8}
  Let $\gamma_n\in\mcC(Z)$ be a sequence of solutions to  the  basic Seiberg-Witten equations on manifold with the codimension $4$ foliation, $\mathfrak F(\gamma)=0$, on $Z=[t_1,t_2]\times M$. Suppose that $$\L(\gamma_n(t_1))-\L(\gamma_n(t_2))\leq C.$$
  Then, there is a sequence of gauge transformations $u_n\in \mcG(Z)$, such that, after passing to a subsequence, the sequence of transformed solutions $\{u_n\gamma_n\}$ converges in the $C^\infty$ topology on $[t'_1,t'_2]\times M$ for a smaller interval $[t'_1,t'_2]$ in the interior of $[t_1,t_2]$.
\end{thm}

\begin{pf}
  The formula \eqref{formula-analytic energy} implies that there is  a uniform bound for the analytic energy. We take a gauge-fixing action by Theorem \ref{Neumann-estimate-thm}.  The arguments of   \cite[Theorem 5.1.1]{KM} shows that there are uniform bounds for the following terms
  \[\|\Phi_n\|_{L^4}\leq C,~ \|F_{A_n}
  \|_{L^2}\leq C.\]
  For each $a_n=A_n-A_0$, there is  a gauge action $u_n$, to make $a^1_n=u_n(a_n)$ satisfy the Coulomb-Neumann condition(see \cite[Page102]{KM}).  %For  each homotopy class we  choose an element satisfying the equations \eqref{eqn-gauge-Neumann-condition},
  Hence, for each $a_n$ we have  a gauge action $v_n$ satisfying  the equations \eqref{eqn-gauge-Neumann-condition} such that there is a uniform constant $C$(independent of $n$) to make the following estimate holds:
  \begin{equation}
    \int_Z\beta_i\wedge a^1_n\wedge\chi_F\in [-C,C],
  \end{equation}
  where $\{\beta_i\}_{1\leq i\leq b^1_b}$ is dual basic to the lattice $H^1_b(Z)\cap H^1(Z;\mathbb Z)$. %A well-know theorem says that  :   If $\{x_n\}$ is a weakly convergent sequence in a Banach space, then the norm of the limit $x$ can only be smaller than or equal to the limit inferior of the norms of the elements of the sequence.
 This implies that there is a uniform $L^2_1$-bound on $(A_n-A_0,\Psi_n)$ up to gauge. Let $(A,\Psi)$ be the weak limit of $(A_n,\Psi_n)$. One has that
\[\sup\mathcal E^a(A_n,\Psi_n)\geq \mathcal E^a(A,\Psi).\]
 The drop of the Chern-Simons-Dirac functional is bounded above, which  implies that
 up to subsequence the following holds
 \[\mathcal E^a(A_n,\Psi_n)\to \mathcal E^a(A,\Psi).\]
  Therefore, the sequence  $(A_n-A_0,\Psi_n)$ converges in $L^2_1$ up to a subsequence.  The remainder  argument is parallel to \cite[Page 107-108]{KM}.
   \end{pf}

We also need to define the perturbation, which is similar to \cite[Section 10, 11]{KM}.%, we also need to define the perturbation.
Let $\mathcal V_k(Z)$ be the $L^2_k$-completion of $\Omega^+(Z)\oplus \Gamma(Z,S^-)$.

\begin{defi}[c.f. {\cite[Definition 10.5.1]{KM}}]\label{defin-perturbation}
  A perturbation $\mathfrak q$ is called $k$-\emph{tame}, if  it is a formal gradient of a continuous $\mcG_b(M)$-invariant function on $\mcC(M)$ and satisfies the following properties:
  \begin{enumerate}
    \item the corresponding codimension $4$ perturbation $\hat{\mathfrak q}$ defines a smooth section $\hat{\mathfrak q}: \mcC_k(Z)\to \mathcal V_k(Z)$(see \cite[Formula 10.2]{KM}), where $Z=I\times M$;
    \item for each $i\in[1,k]$, the codimension $4$ perturbation $\hat{\mathfrak q}$ defines a continuous section
        $\hat{\mathfrak q}: \mcC_j(Z)\to \mathcal V_j(Z)$;
    \item the derivative
    $$D\hat{\mathfrak q}: \mcC_k(Z)\to
    Hom(T\mcC_k(Z),\mathcal V_k(Z))$$
    extends to a map
    $$D\hat{\mathfrak q}: \mcC_k(Z)\to
    Hom(T\mcC_j(Z),\mathcal V_j(Z)),$$
    for $j\in[-k,k]$;
    \item we have the estimate
    \[\|\mathfrak q(A,\Psi)\|_{L^2}\leq C(\|\Psi\|_{L^2}+1),\]
    for some constant $C$ and each $(A,\Psi)\in\mcC_k(M)$;
    \item for any reference connection $A_0$, we have
    \[\|\hat{\mathfrak q}(A,\Psi)\|_{L^2_{1,A}}\leq f_1(\|A-A_0,\Psi\|_{L^2_{1,A_0}}),\]
    where $f_1$ is a real function and $(A,\Psi)\in\mcC_k(Z)$;
    \item $\mathfrak q$ defines a $C^1$-section
    $\mcC_1(M)\to\mathcal T_0$.
  \end{enumerate}
  We say $\mathfrak q$ is \emph{tame}, if it is $k$-\emph{tame} for all $k\geq2$.
\end{defi}

Using the Theorem \ref{thm-coro5.1.8} and the same idea of \cite[Proposition 16.2.1]{KM}, one establishes the following proposition.

 \begin{prop}[c.f. {\cite[Proposition-16.2.1]{KM}}]\label{cmpt-prop}
   For any $C>0$, there are only finitely many $[\alpha],~[\beta]$ and pathes $z$ with
   $\mathcal E^{top}_q(z)\leq C$, such that the space $\check{\mathcal M}^+_z([\alpha],[\beta])$ is non-empty. Furthermore, each $\check{\mathcal M}^+_z([\alpha],[\beta])$ is compact.
 \end{prop}

\section{Compactness of moduli space }
%In this section, we will give a gluing theorem for.
The purpose of this section is for the preparation to construct the monopole Floer homologies in the next section. %In the first subsection, we will recall some facts about the gluing model for general Hilbert space; in the second subsection we  will give the gluing theorem for the blow-down model; in the third subsection, we

\subsection{Gluing trajectories}
%\subsubsection{Preparation}
In this subsection we show the gluing theorem in gauge version.
 %Let $X=M$ or $X=I\times M$ be the compact oriented manifold with oriented taut  Riemannian foliation $F$.
 We set the blow-up configuration space of $(M,F)$  as the space of triples, see \cite[Chapter 9]{KM}
 \[\mathcal C^\sigma_k(M,F,\mfs)=\{(A,s,\phi)\big| (A,\phi)\in \mathcal C_k(M),~ s\in\mathbb R^{\geq0},~\|\phi\|_{L^2}=1\}.\]
 We define   the quotient space $\mathcal B^\sigma_k(M,F,\mfs)=\mathcal C^\sigma_k(M,F,\mfs,F,\mfs)/\mathcal G_{k+\frac12}$, which is a Hilbert manifold with boundary.
  The basic Seiberg-Witten equations naturally extend to the equations
  \[\begin{cases}
    \Dirac_{b,A}\phi=0\\
    \frac12\bar*F_{A^t}=s^2q(\phi).
  \end{cases}\]
  We define \[\tilde{\mathcal C}^\tau_k(Z,F,\mfs)\subset \mathcal A_k(Z,\mfs)\times L^2_k(I,\mathbb R)\times L^2_k(S^+) \]
  to be the subset consisting of triples $(A,s,\phi)$ with $\|\phi(t)\|_{L^2(M)}=1$ for each $t\in I$. We have an involution map $\mathfrak i: \tilde{\mathcal C}^\tau_k(Z,F,\mfs)\to \tilde{\mathcal C}^\tau_k(Z,F,\mfs)$ defined by $(A,s,\phi)\mapsto (A,-s,\phi)$.
  Similarly, we define the $\tau$-module for the configuration space for  $(Z,F)$, where $Z=I\times M$
  \[{\mathcal C}^\tau_k(Z,F,\mfs)\subset \mathcal A_k(Z,\mfs)\times L^2_k(I,\mathbb R)\times L^2_k(S^+) \]
  to be the subset consisting of triples $(A,s,\phi)$ with
  \[s(t)\geq0,~\|\phi(t)\|_{L^2(M)}=1,\]
  for each $t\in I$. %A configuration $(A,s,\phi)$ on the blow-up space $\mathcal C^\sigma_k(I\times M,F,\mfs)$ gives rise an element of $\tau$-model or a path of configurations on $M$ given by $$(A(t),s\|\phi(t)\|_{L^2(M)},\phi(t)/\|\phi(t)\|_{L^2(M)}),$$  provided that $\|\phi(t)\|_{L^2(M)}\neq0$ for each $t\in I$.
There is a well-defined map
$$\pi: \mathcal C^\sigma_k(M,F,\mfs)\to \mathcal C_k(M,F,\mfs),~(A,s,\phi)\mapsto (A,s\phi).$$
By this map,  any vector field on $\mathcal C_k(M,F,\mfs)$ lifts to a vector field on $\mathcal C^\sigma_k(M,F,\mfs)$.
In order to get the
transversality condition, we need to add a perturbation $\mathfrak p$ as defined in \eqref{defin-perturbation} on $grad(\L)$. The sum $grad(\L)+\mathfrak p$ is a gauge invariant and gives rise to a vector field $v^\sigma_q$,
\[v^\sigma_q:\mathcal B^\sigma_{k+1/2}(M,F,\mfs)\to \mathcal T_{k-1/2}(M),\]
where $\mathcal T_{k-1/2}(M)$ denotes the $L^2_{k-/2}$-completion of the tangent bundle of $\mathcal B^\sigma_{k-1/2}(M,F,\mfs)$.   $grad^\sigma(\L)$ is defined as follows
 \[grad^\sigma(\L)(A,r,\psi)=(\frac12\bar*F_{A^t}-r^2q(\psi),\Lambda(A,r,\psi)r,\Dirac_A\psi-
 \Lambda(A,r,\psi)\psi),\]
 where $\Lambda(A,r,\psi)=\langle\psi,\Dirac_A\psi\rangle_{L^2}$. %If $(A,0)$ is a critical point of $\L$, then the equation for $\phi$ is  \[\Dirac_A\phi=\Lambda(A,0,\phi)\phi=0,\]  i.e. $\phi$ is a an eigenvector of $\Dirac_{A}$.
  A  trajectory $\gamma(t)=(A(t),r(t),\psi(t))$ is a solution to the equations
 \begin{equation}
   \begin{cases}
     \frac12\frac{d}{dt}A=-(\frac12\bar*F_A-r^2q(\psi)),\\
     \frac{d}{dt}r=-\Lambda(A,r,\psi)r,\\
     \frac{d}{dt}\psi=-(\Dirac_A\psi-
 \Lambda(A,r,\psi)\psi)).
   \end{cases}
 \end{equation}
 We call the perturbation $\mathfrak q$ admissible, if all critical points  of  $v^\sigma_q$ are nondegenerate and moduli spaces of the flow lines connecting them  are regular(see \cite[Defition 22.1.1]{KM}).
We categorize the set $C$ of critical points in $\mathcal B^\sigma_k(M,F,\mfs)$ into the disjoint union of three subsets:
\begin{itemize}
  \item $C^o$, the set of irreducible points;
  \item $C^s$, the set of reducible boundary stable(where the spinor part locates on the positive eigenspace part)   critical points;
  \item $C^u$, the set of reducible  boundary unstable(where the spinor part locates on the negative eigenspace part) critical points.
\end{itemize}

We set
\[M([\mfa],[\mfb])=\bigcup_{z\in\pi_1([\mfa],[\mfb])}
 M_z([\mfa],[\mfb])\]
as the union over the moduli spaces in a given relative homotopy class, where $\pi_1([\mfa],[\mfb])$ denotes the homotopy class of path connecting the two critical points in the quotient space. We recall two notions which are given in  \cite[Page 261]{KM}.

\begin{defi}
  We say that a moduli space $M([\mfa],[\mfb])$ is \emph{boundary-obstructed}, if $[\mfa],~[\mfb]$ are reducible, $[\mfa]\in C^s$ and $[\mfb]\in C^u$.
\end{defi}
\begin{defi}
   %We choose an element $[\gamma]\in M_z([\mfa],[\mfb])$.
   When the moduli space $M_z([\mfa],[\mfb])$ is not boundary-obstructed, we say that $\gamma$ is \emph{regular}, if the linearized Seiberg-Witten map $Q_\gamma$ along $\gamma$ is surjective for $[\gamma]\in M_z([\mfa],[\mfb])$. In the boundary-obstructed case, we say $\gamma$ is \emph{regular}, if the linearized Seiberg-Witten map restriction along the boundary $Q^\partial_\gamma$ is surjective for $[\gamma]\in M_z([\mfa],[\mfb])$. We say that $ M_z([\mfa],[\mfb])$  is regular if its elements are all regular.
\end{defi}
We topologize the space of unparameterized broken trajectories as follows\cite[Page 276]{KM}.
Choose an element $[\breve\gamma]=([\breve\gamma_1],\cdots,[\breve\gamma_n])\in \breve M^+_z([\mfa],[\mfb])$, where  $[\breve\gamma_i]\in \breve M_{z_i}([\mfa_{i-1}],[\mfa_i])$ is represented by a trajectory
\[[\gamma_i]\in  M_{z_i}([\mfa_{i-1}],[\mfa_i]).\]
Let $U_i\in \mcB^\tau_{k,loc}(Z,F,\mfs)=\mathcal C^\tau_{k,loc}(Z,F,\mfs)/\mcG_{k+1,loc}$ be an open neighborhood of $[\gamma_i]$ and $T\in \mathbb R^+$ be a positive number. We define $\Omega=\Omega(U_1,\cdots,U_n,T)$ to be the subset of $\breve M_z([\mfa][\mfb])$ consisting of unparameterized broken trajectories
$\breve\delta=([\breve\delta_1],\cdots,[\breve\delta_m])$ satisfying the following condition: there exists a map
$$(\iota,s):\{1,\cdots, n\}\to \{1,\cdots,m\}\times \mathbb R$$
such that
\begin{itemize}
  \item $[\tau_{s(i)}\delta_{\iota(i)}]\in U_i$,
  \item if $1\leq i_1<i_2\leq n$, then either $\iota(i_1)<\iota(i_2)$ or $\iota(i_1)=\iota(i_2)$ and $s(i_1)+T\leq s(i_2)$.
\end{itemize}
%\begin{defi}[c.f. {\cite[Page 276]{KM}}]\label{defin-separating-property}   Let $\mathcal U$ be a collection of gauge-invariant open neighborhoods $U_a\subset \mathcal C_{k,b}(I\times M,F,\mfs)$ of the critical points. We say that $\mathcal U$ has the \emph{separating property} for $I\times M$, if the following statements hold: there should exist neighborhoods $V_{[a]}\subset   \mathcal B_{k-1,b}(M,F,\mfs)$ of the critical points $[a]\in \mathcal B_{k-1,b}(M,F,\mfs)$ such that:   \begin{itemize}     \item the sets $V_{[a]},~V_{[b]}$ are disjoint whenever $[a]\neq[b]$,     \item each $V_{[a]}$ is path-connected and simply connected,     \item if $\gamma$ belongs to $U_a$, then the gauge-equivalent class of $\gamma(t)$ is in $V_{[a]}$, for all $t\in I$.   \end{itemize} \end{defi} Note that the same definition  can be repeated verbatim in the blown-up context.

%At the last of this subsection, we review the perturbation \cite[Section 10]{KM}. A perturbation $\mathfrak q$ is regarded as a section\[\mathfrak q:\mcC(M)\to \mathcal T_0,\] where $\mathcal T_0$ denotes the tangent bundle of $\mcC(M)$.

%\subsubsection{Compactness of the blow-up model}
To prove the compactness theorem for the blow-up model, we need to get a bound of $\Lambda$  of the moduli space.
For any $C>0$ and any $[\mfa],~[\mfb]$ with energy $\mathcal E^a_q\leq C$ for which $M^+([\mfa],[\mfb])$ is non-empty,   the space of broken trajectories $[\gamma]\in M^+([\mfa],[\mfb])$ with energy $\mathcal E_q\leq C$ is compact.
For a trajectory $\gamma^\tau\in M([\mfa],[\mfb])$, we define $K(\gamma^\tau)$ to be the total variation of $\Lambda_{\mathfrak q}$ by(see \cite[Section 16.3]{KM})
\[K(\gamma^\tau)=\int_{\mathbb R}|\frac{d\Lambda_{\mathfrak q}(\gamma^\tau)}{dt}|dt.\]
Set
\[K_+(\gamma^\tau)=\int_{\mathbb R}\left(\frac{d\Lambda_{\mathfrak q}(\gamma^\tau)}{dt}\right)^+dt.\]
Proposition \ref{cmpt-prop}    gives a bound on the number of components for which the blow-down is non-constant.
To get the energy bound, we  need the proposition below.
 \begin{prop}[c.f. {\cite[Proposition 16.1.4]{KM}}]\label{prop-16.1.4}
    The space $\breve M^+([\mfa],[\mfb])$ of broken trajectories %$[\breve\gamma]\in\breve M^+([\mfa],[\mfb])$
    with topology energy $\mathcal E_{\mathfrak q}(\breve\gamma)\leq C$ is compact.
  \end{prop}

Set two spaces $Z^T$ and $Z^\infty$ as $Z^T=[-T,T]\times M$ and $Z^\infty=(\mathbb R^{\leq}\times M)\coprod(\mathbb R^{\geq}\times M)$ respectively. Let $\mfa$ be a critical point on $\mathcal C_{k}(M,F,\mfs)$, we write $\gamma_\mfa$ as a translation-invariant solution on $Z^T$ or $Z^\infty$ in temporal gauge.
We assume that $\mfa=(A_0,\Phi_0)$ is non-degenerate by choosing a generic perturbation. We define the quotient space
\[\mathcal B_k(Z^\infty,[\mfa])=\mathcal C_k(Z^\infty,[\mfa])/\mathcal G_{k+1}(Z^\infty),\]
where $C_k(Z^\infty,[\mfa])=\{\gamma\in\mathcal C_{k,loc}(Z^\infty,F,\mfs)\big|
\gamma-\gamma_\mfa\in L^2_{k,A_0}\}$ and $\mathcal G_{k+1}=\{u\in \Gamma_b(Z^\infty,S^1)\big| u\in L^2_{k+1,loc}~,1-u\in L^2_{k+1}\}$.
Let \[\mathcal K_{s,\mfa}(M)\]
be the $L^2_{s}$-completion of the complement $\mathcal K_{\mfa}$ to the gauge orbit, where $\mathcal K_{\mfa}$ is the orthogonal complement to the gauge-orbit. Similarly, we denote $\mathcal K^\sigma$  % and $\mathcal S^\sigma$
by the blow-up model of $\mathcal K$. Similar to \cite[Proposition 9.3.4]{KM},  the proposition below holds. % and $\mathcal S$ respectively.

\begin{prop}\label{prop-slice}%[c.f. {\cite[Proposition 9.3.4]{KM}}]
  Let $\mathcal J_{k,\gamma}$ be the image of $d_\gamma: L^2_{k+1,b}(M,i\mathbb R)\to \mathcal T_{k,\gamma}$, via
  $\xi\mapsto (-d\xi,\xi \Phi_0)$, where $\mathcal T_{k,\gamma}$ denotes the tangent space at $\gamma$. As $\gamma$ varies over $\mathcal C^*_k(M)$, we define $\mathcal K_{k,\gamma}$ to be the subspace of $\mathcal T_{k,\gamma}$, which is orthogonal to $\mathcal J_{k,\gamma}$ with respect to the $L^2$-inner product. Then, we have the decomposition
  \[\mathcal T_{k,\gamma}=\mathcal J_{k,\gamma}\oplus\mathcal K_{k,\gamma}.\]
\end{prop}

\begin{pf}
  We denote by $\gamma=(A_0,\Phi_0)$. By doing integral by part,  we define $\mathcal K_{k,\gamma}=\{(a,\phi)\big|
  -\delta_ba+iRe\langle i\Phi_0,\phi\rangle=0\}$.  It is clear that $\mathcal K_{k,\gamma}$ is orthogonal to $\mathcal J_{k,\gamma}$. We need to show that
  \[\mathcal T_{k,\gamma}=\mathcal J_{k,\gamma}\oplus\mathcal K_{k,\gamma}.\]
  It is sufficient to show that for any $(a,\phi)$, there is a unique solution to the equation
  \[
    \delta_b((a,\phi)+d_\gamma(\xi))=0,
  \]
  which is equivalent to
  \[\Delta_b\xi+|\Phi_0|^2\xi=c,\]
  where $c=G^*_\gamma(a,\phi)$.
  Since $\Phi_0$ is non-zero, this equation has a unique solution by Theorem \ref{CN-thm}.
\end{pf}

Following the arguments of \cite[Proposition 9.3.5, 9.4.1]{KM},   we can show the similar decompositions for $\sigma$-model and $\tau$-model.

By doing the integral part with the taut condition, one has that  the slice\[S_{k,\mfa}(Z^T)\subset \mathcal C_{k,b}(Z^T)\] can be represented by \[S_{k,\mfa}(Z^T)=\{(A_0+a,\Phi)\big| -\delta_ba+iRe\langle i\Phi_0,\Phi\rangle=0,~(a,\overrightarrow{n})\big|_{\partial Z^T}=0.\},\]where $\gamma_\mfa=(A_0,\Phi_0)$. %% we need a proof for such a decomposition.
Similarly, we define the slices $\mathcal S^\sigma_{k,\mfa}$ and $\mathcal S^\tau_{k,\mfa}$(see \cite[Page 144, Page 147]{KM}).
%Note that even without the taut condition, the slice decomposition also holds.
We can run the arguments in \cite[Section 18.4]{KM} in the foliation case.
Here we give a sketch.
The boundary of $Z^T$ is $\bar{M}\amalg M$, we have the restriction map\[r:\mathcal{\tilde C^\tau}_{k,b}(Z^T)\to \mathcal{\tilde C^\tau}_{k-\frac12,b}(\bar{M}\amalg M)\times L^2_{k-\frac12,b}(\bar{M}\amalg M,i\mathbb R),\]where the second component is defined by the normal component of the basic connection $A$ at the boundary and
$\bar M$ is a copy of $M$ with the reversing orientation by reversing the orientation of $Q$.
For the non-degeneracy of the Hessian operator \eqref{Hessian-operator} at $\mfa$, we have the decomposition $\mathcal K^\sigma_{k-\frac12}\big|_{\mfa}= \mathcal K^+\oplus \mathcal K^-$.
Let $H^-_M$ and $H^-_{\bar M}$ be two subspaces defined by\[ H^-_M=\{0\}\oplus \mathcal K^-\oplus L^2_{k-\frac12}(M,i\mathbb R),~H^-_{\bar M}=\{0\}\oplus \mathcal K^+\oplus L^2_{k-\frac12}(M,i\mathbb R).\]
%where we used the decomposition\[\mathcal T_{k-\frac12,b}\big|_{\mfa}\oplus L^2_{k-\frac12,b}(M,i\mathbb R)=\mathcal J_{k-\frac12,b}\big|_{\mfa}\oplus\mathcal K_{k-\frac12,b}\big|_{\mfa}\oplus L^2_{k-\frac12,b}(M,i\mathbb R).\]
We define $H=H^-_M\oplus H^-_{\bar M}$ and define $\Pi_M= \Pi^-_M\oplus \Pi^-_{\bar M}$ by the projection to the space $\mathcal K^-\oplus \mathcal K^+$, i.e.
\[\Pi_M:\mathcal T_{k-\frac12}\big|_{\mfa}(M\amalg \bar M)
\oplus L^2_{k-\frac12}(M\amalg\bar M,i\mathbb R)\to \mathcal K^-\oplus \mathcal K^+.\]
For $\gamma $ in a small neighborhood of $\gamma_\mfa\in \mathcal C^\tau_{k,b}(Z,F,\mfs)$ on $Z=Z^T$ or $Z^\infty$, we consider the equations
\[\begin{cases}
  \mathfrak F_q(\gamma)=0,\\
  \gamma\in S^\tau_{k,\mfa}(Z),\\
  (\Pi_M\comp r)(\gamma)=h,
\end{cases}\]
where $h\in H$. We  write
the equations as
\[\begin{cases}
  (Q_{\gamma_\mfa}+\alpha)\gamma=0,\\
  (\Pi_M\comp r)\gamma=h,
\end{cases}\]
where $Q_{\gamma_\mfa}$ is defined by $D_{\gamma_\mfa}\mathfrak F_q$ %\oplus G^*_{\gamma_\mfa}$
 and $\alpha$ denotes the remainder terms.
We  write $Q_{\gamma_\mfa}=\frac{d}{dt}+L_b$, let $H^\pm_L$ be the spectral subspaces of $L_b$ in $L^2_{\frac12}(M)$.
 By Proposition \ref{prop-17.2.7},  the linear map
\begin{equation}
  (Q_{\gamma_\mfa},\Pi^-_L\comp r) \label{Fredholm-formula}
\end{equation}
is an isomorphism, where $\Pi^-_L$ is the spectral projection with kernel $H^+_L$.
Let $K$ denote the kernel of $Q_{\gamma_\mfa}$, the domain can be decomposed as $C\oplus K$. We rewrite the above operator as
\[\left(\begin{array}{cc}
  Q_{\gamma_\mfa}\big|_C&0\\
  *&(\Pi^-_L\comp r)\big|_K\\
\end{array}\right).\]
The isomorphism of  two components on diagonal implies that the matrix is an invertible operator. Thus, we   verified the  abstract hypothesis \cite[Hypothesis 18.3.1]{KM}. %\ref{hypo-ivertible}.and  Hypothesis \ref{hypo-continuous}.
 By the definition of the tame perturbation, the abstract hypothesis  \cite[Hypothesis 18.3.3]{KM} follows.  Setting $\mathcal K=\mathcal K^+\oplus \mathcal K^-$,
 there is  an $ \eta_1 > 0$ and  two maps from the $B_{\eta_1}(\mathcal K)$ to the slices parameterizing the  subsets of the set of solutions, i.e.
\[u(T,\cdot): B_{\eta_1}(\mathcal K)\to S^\tau_{k,\gamma_\mfa}(Z^T)\cap \mathfrak F^{-1}_q(0),~
u(\infty,\cdot): B_{\eta_1}(\mathcal K)\to S^\tau_{k,\gamma_\mfa}(Z^\infty)\cap \mathfrak F^{-1}_q(0),\]
for some positive number $\eta_1>0$.
By the parallel arguments of the proof of \cite[Theorem 18.2.1]{KM},  one has the proposition below.
%By Theorem \ref{thm-gluing-1}, one has the proposition below.
%\begin{thm}[c.f. {\cite[Theorem 18.2.1]{KM}}]\label{thm-18.2.1}  There exists $T_0$, such that for all $T\geq T_0$, wee can find smooth maps\[u(T,\cdot): B (\mathcal K)\to M(Z^T),~u(\infty,\cdot): B (\mathcal K)\to M(Z^\infty),\]which are diffeomorphisms from a  ball $B (\mathcal K)$ of $\mathcal K$ onto neighborhoods of theconstant solution $[\gamma_{\mfa}]$. These can be chosen so that the map\[\mu_T:B (\mathcal K)\to\mcB_{k-1/2}(\partial Z^T) \]defining by composing $u(T,-)$ with the restriction maps to the boundary is a smooth embedding of $B (\mathcal K)$ for $T\geq T_0$. Moreover, we have that, as function on $[T_0,\infty)\times B (\mathcal K)$, the map $(T,h)\mapsto\mu_T(h)$ is smooth for finite $T$, and $\mu_T\to\mu_\infty$ in $C^\infty_{loc}$. Finally, there is an $\eta$(irrelative to $T$), such that the images of the maps $u(T,-)$ can be taken to contain all $[\gamma]\in M(Z^T)$ with $\|\gamma-\gamma_\mfa\|_{L^2_{k,\mfa}(Z^T)}\leq \eta$, denoted by $M_\eta(Z^T,[\mfa])$\end{thm}

%We need the next proposition to completing the gluing theorem.

\begin{prop}%[c.f. {\cite[Proposition 18.4.1]{KM}}]
There is an $\eta_0$, such that all $\eta<\eta_0$, there is a number $\eta'$, independent of $T$, such that:  \begin{enumerate}
 \item the map    \[\mu:\{\gamma\in  S^\tau_{k,\gamma_\mfa}(Z^T)\big |\|\gamma-\gamma_\mfa\|_{L^2_k}\leq\eta\}\to \tilde{\mathcal B^\tau_k}(Z^T)   \]is a diffeomorphism onto its image,   where $\tilde{\mathcal B^\tau_k}(Z^T) $ denotes the quotient of $\tilde{\mathcal C^\tau_{k}}(Z^T)$ by the gauge action;
     \item the image of the above map contains all gauge-equivalent classes $[\gamma]\in \tilde{\mathcal B^\tau_k}(Z^T) $ represented by the elements $\gamma$ satisfying        \[\|\gamma-\gamma_\mfa\|_{L^2_k}\leq \eta'.\]  \end{enumerate}\end{prop}
          By the similar arguments, we have the following   propositions for the foliated case.

%  \begin{prop}  Let $x_1$ and $x_2$ be two solutions of the  basic Seiberg-Witten equations on $Z=[t_1,t_2]\times M$. Suppose there exists some $t_0$ of the interval such that $x_1(t_0)$ and $x_2(t_0)$ are gauge-equivalent on $\{t_0\}\times M$. Then $x_1$ and $x_2$  are gauge-equivalent on $Z$.  \end{prop}

  \begin{prop}[c.f. {\cite[Proposition 17.2.7]{KM}}]\label{prop-17.2.7}    Let $Z=(-\infty,0]\times M$ and $D_0:C^\infty(Z,E)\to L^2(Z,E)$ be a  transverse elliptic operator of the form    \[D_0=\frac{d}{dt}+L_0,\]    where $L_0:C^\infty(M,E_0)\to C^\infty(M,E_0)$ is a transverse self-adjoint elliptic operator on $M$. Suppose  that the spectrum of $L_0$ does not contain zero. Then, the operator    \[D_0\oplus(\Pi_0\comp r):L^2_j(Z,E)\to L^2_{j-1}(Z,E_0)\oplus (H^-_0\cap L^2_{j-\frac12}(M,E_0))\]    is an isomorphism for all $j\geq1$, where $\Pi_0$ denotes the projection to the negative eigen-vector part of $L_0$. Moreover, it holds that  $H^-_0\cap L^2_{j-\frac12}(M,E_0)=Im(r\big|_{\ker(D_0)})$.  \end{prop}

%  \begin{thm}[c.f. {\cite[Theorem 17.1.3]{KM}}]\label{thm-fredholm-regularity}    Let $Z=I\times M$ be a finite cylinder and $D :L^2_j(Z,E)\to L^2_{j-1}(Z,E)$ be an basic elliptic operator of the form    \[D=\frac{d}{dt}+L_0,\]    where $L_0:C^\infty(M,E_0)\to C^\infty(M,E_0)$ is a basic self-adjoint elliptic operator on $M$. Let $\Pi_0$ denote the projection to the negative eigen-vector part of $L_0$. Suppose  that the spectrum of $L_0$ does not contain zero. Then, the following statements hold:    \begin{enumerate}      \item For $1\leq j\leq k$,      \[D\oplus (\Pi_0\comp r):L^2_j(Z,E)\to L^2_{j-1}(Z,F)\oplus (H^-_0\cap      L^2_{j-\frac12}(M,E_0))\]      is Fredholm.      \item If $\{u_i\}$ is a bounded sequence in $L^2_j(Z,E)$ and $\{Du_i\}$ is a Cauchy sequence in $L^2_{j-1}$, then $\{(1-\Pi_0)r(u_i)\}$ has a convergent subsequence in $H^+_0\cap L^2_{j-\frac12}$.      \item  If $u\in L^2_j(Z,E)$ for $j\leq k$ and the image of $u$ under  $D\oplus (\Pi_0\comp r)$ lies in $L^2_{k-1}(Z,F)\oplus (H^-_0\cap          L^2_{k-\frac12}(M,E_0))$, then $u\in l^2_k(Z,E))$.    \end{enumerate}  \end{thm}

  The proof only needs the parametrix patching and regularity. %, thus we omit it here.
   Let $Z=I\times M$ be a  closed  finite cylinder, suppose that $I=I_1\cup I_2$ with $I_1\cap I_2=\{0\}$. We denote by $Z=Z_1\cup Z_2$, where $Z_1=I_1\times M$ and $Z_2=I_2\times M$. Let $D:C^\infty(Z,E)\to L^2(Z,E)$ be a transverse  elliptic operator of the form\[D=\frac{d}{dt}+L_0+h(t),\]where $L_0$ is a self-adjoint operator on $M$, $h: L^2_j(Z,E)\to L^2_{j}(Z,E)$ is a bounded operator. We set $D_1,~D_2$ as  the restriction of these operators to the two subcylinders respectively, and set
  \[H^i_{j-\frac12}\subset L^2_{j-\frac12}(\{0\}\times M,E_0)\]as  the image of the $\ker(D_i)$ under the restriction map\[r_i: L^2_{j}(Z_i,E)\to L^2_{j-\frac12}(\{0\}\times M,E_0). \]Denoting by $D_0=-\frac{d}{dt}+L_0$.
  We have the following lemma.
%\begin{lemma}[c.f. {\cite[Lemma 17.2.9]{KM}}]  Let $\Pi_0$ be the spectral projection of the self-adjoint elliptic operator $D_0$, whose image is the negative subspace $H^-$. If $b\in H^1_{j-\frac32}$ for some $2\leq $ and finite cylinder $Z=I\times M$,  then one has that   \[(1-\Pi_0)b\in L^2_{j-\frac12}.\] \end{lemma}

% \begin{pf}   Let $v$ be the solution to $D_0v=0$ on $Z$, such that $v\big|_{\{0\}\times M}   =\Pi_0b$. This implies that $v\in L^2_{j-1}(Z_1,E)$ by \cite[Theorem 17.1.4]{KM}. From the definition of $H^1_{j-\frac32}$, there exists $v'\in L^2_{j-1}(Z_1,E)$ with $D_1v'=0$ and $v'\big|_{\{0\}\times M}=b$. We have   \[D_0(v-v')=h(v')\in L^2_{j-1}(Z_1,E_0),\]   and $\Pi_0(v-v')\big|_{\{0\}\times M}=0$. These two conditions imply that $v-v'\in L^2_j$ on a neighborhood of $\{0\}\times M$ in $Z_1$ on the half-infinite cylinder.  Applying $(1-\Pi_0)$ to the boundary value of $v-v'$, we get that $(1-\Pi_0)b\in L^2_{j-\frac12}$. \end{pf}

 \begin{prop}[c.f. {\cite[Proposition 17.2.8]{KM}} ]\label{prop-surjective-on-subinterval}
   Suppose that $D:L^2_j(Z)\to L^2_{j-1}(Z)$ is surjective for $2\leq j$. Then, we get the decomposition
   \[L^2_{j-\frac12}(\{0\}\times M,E_0) = H^1_{j-\frac12}+H^2_{j-\frac12}. \]
   Conversely, if the above formula holds and $D_1,~D_2$ are surjective, then $D$ is surjective.
 \end{prop}
% \begin{pf}   When $h=0$, since $H^1_{j-\frac12}$ and $H^2_{j-\frac12}$ contain the negative and positive spectral subspaces of $L_0$,  the result trivially holds.   In general, let $u_i\in L^2_1(Z_i)$ for $i=1,2$ and satisfy $D_0u_i=0$ with $R_1u_1-R_2u_2=a$. Let $u\in L^2(Z,E)$ be the element that equals to $u_i$ on $Z_i$ for $i=1,~2$. Since $D$ is surjective, there exists $w\in L^2_1(Z)$ satisfying $Dw=hw$. We define $w_i$ be the restriction of $w$ to $Z_i$, and $u'_i=u_i-w_i$. Then, by $R_1w_1=R_2w_2$, one gets that   \[R_1u'_1-R_2u'_2=a,\]   and   \[D_i(u'_i)=D u_i-Dw_i=hu_i-Dw_i=0,\]   which exhibits $a$ as $R_1u'_1-R_2u'_2\in H^1_{\frac12}+H^2_{\frac12}$.   Now, suppose that $j>1$, and $a\in L^2_{j-\frac12}(\{0\}\times M,E_0)$, by the identity   \[a=a_1+a_2,\]   with $a_i\in H^i_{\frac12}$. By an induction hypothesis that $a_i\in H^i_{j-\frac32}$, we get $(1-\Pi_0)a=(1-\Pi_0)a_1+(1-\Pi_0)a_2$. The previous lemma tells that $(1-\Pi_0)a_1\in L^2_{j-\frac12}$, so $(1-\Pi_0)a_2\in L^2_{j-\frac12}$, for $a\in L^2_{j-\frac12}$. By applying the lemma with opposite signs, we get $\Pi_0a_2\in L^2_{j-\frac12}$, which completes the proof. \end{pf}

  Let $Z$ be a finite cylinder.
 We set $\tilde M(Z)=\{[\gamma]\in\tilde{\mathcal B^\tau_k}(Z)|~\mathfrak F^\tau_{\mathfrak q}(\gamma)=0 \}$.
  We have the following theorem.
 \begin{thm}[c.f. {\cite[Theorem 17.3.1]{KM}}]\label{thm-17.3.1}
    The subspace $\tilde M(Z)\subset \tilde{\mathcal B^\tau_k}(Z)$ is a closed Hilbert submanifold.
The subset $M(Z)$ is a Hilbert submanifold with boundary, i.e. it is
identified with the quotient of $\tilde M(Z)$ by the involution $\mathfrak i$.
 \end{thm}

Let $[\gamma]\in\tilde M(Z)$ and let $\bar\mfa$ and $\mfa$ be the restrictions of $\gamma$ to the two boundary components. We have the restriction maps
\[R_+:\tilde M(Z)\to\mcB^\sigma_{k-1/2}(M),~R_{-}:\tilde M(Z)\to \mcB^\sigma_{k-1/2}(\bar M).\]
\begin{thm}[c.f. {\cite[Theorem 17.3.2]{KM}}]\label{thm-17.3.2}
  Let $\gamma,~\mfa$ and $\bar\mfa$ be as above, and let $\Pi:\mathcal K^\sigma_{\bar\mfa}(\bar M)\oplus \mathcal K^\sigma_{\mfa}( M)\to \mathcal K^-_{\bar \mfa}(\bar M)\oplus \mathcal K^-_{\mfa}( M)$ be the projection with kernel $\mathcal K^+_{\bar \mfa,\mfa}(\bar M\amalg M)$. Then, the two composition maps
  $\Pi\comp (DR_-,R_+)$ and $(1-\Pi)\comp (DR_-,R_+)$ are Fredholm and compact respectively, where $DR_-$ and $DR_+$ denote the  derivatives of $R_-$ and $R_+$ respectively.
\end{thm}

We can prove the foliated version of \cite[Lemma 16.5.3, Proposition 16.5.2,  Proposition 16.5.5]{KM}. By Lemma \ref{lemma-7.1.3}, Proposition \ref{prop-surjective-on-subinterval}, Theorem \ref{thm-17.3.1} and Theorem \ref{thm-17.3.2},
  there is no difficulty  to apply  similar  arguments % of \cite[Section 19.2,Section 19.3]{KM}. Repeat the similar arguments
   of \cite[Section 19.1, Section 19.2, Section 19.3 and Section 19.4]{KM}, one establishes the following theorems.

\begin{thm}[c.f. {\cite[Theorem 19.5.4]{KM}}]\label{thm-19.5.4}  Suppose that the moduli space $M_z([\mfa],[\mfb])$ is $d$-dimensional and contains irreducible trajectories, such that the moduli space  $\breve M^+([\mfa],[\mfb])$ is a $(d-1)$-dimensional space stratified by manifolds(see \cite[Definition 16.5.1]{KM}). Let $M'\subset \breve M^+([\mfa],[\mfb])$ be any component of the codimension-1 stratum. Then along $M'$, the moduli space is either a $C^0$-manifold with boundary, or has a codimension-1 $\delta$-structure in the sense of  \cite[Definition 19.5.3]{KM}. The latter
occurs only when $M'$ consists of 3-component broken trajectories, with the
middle component boundary-obstructed .
\end{thm}

\subsection{Finite result on moduli space}
In this subsection, we will give some properties for the compactified moduli space, which are necessary to construct the basic monopole  Floer  homologies without using the Novikov ring.

 Recall that a reducible critical point $\mfa$ corresponds to a pair $(\alpha,\lambda)$, where $\alpha=(B,0)=\pi(\mfa)$ is a critical point in
$\mathcal C_k(M,F,\mfs)$, and $\lambda$ is an element of the spectrum of
$\Dirac_{B,\mathfrak q}$. For such $\mfa$, we define $\iota(\mfa)$(see \cite[Fromula 16.2]{KM}) by
\[\iota(\mfa)=\begin{cases}
  |Spec(\Dirac_{B,\mathfrak q})\cap [0,\lambda)|,&\lambda>0,\\
  1/2-|Spec(\Dirac_{B,\mathfrak q})\cap[0,\lambda]|,&\lambda<0.
\end{cases}\]
We denote by $M^{red}([\mfa],[\mfb])\subset M([\mfa],[\mfb])$ the subspace consisting of all the
reducible trajectories. For the moduli space of reducible trajectories, we have simple structure, i.e. it is always a manifold without boundary.
For its dimension, we have the formula(see \cite[Formula 16.9]{KM})
\[\dim(M^{red}_z([\mfa],[\mfb]))=\bar{gr}_z([\mfa],[\mfb])=gr_z([\mfa],[\mfb])-
o[\mfa]+o[\mfb],\]
where $o[\mfa]=0$ when $[\mfa]\in C^s$ and $o[\mfa]=1$ when $[\mfa]\in C^u$.
For an irreducible critical point $\mfa$, we set $\iota(\mfa)=0$. If $[\mfa]$ and $[\mfa']$ are two critical points whose images under $\pi$ equal to the same critical point $[\alpha]\in\mathcal B_k(M)$, then we have the following identity
\[gr_{z_0}([\mfa],[\mfa'])=2(\iota(\mfa)-\iota(\mfa'))\]
for a trivial homotopy class $z_0$(see \cite[Formula 16.3]{KM}).
%is the virtual dimension of the moduli space $\breve M^{red}([\mfa],[\mfa'])$.
%The lemma below gives a bound for the energy for $b^1_b=0$.

\begin{lemma}\label{lemma-16.4.4}
  Suppose that  all  moduli spaces are regular and there is positive number $C_0>0$ such that \begin{equation}
    \mathcal E^{top}_{\mathfrak q} (z_u)+C_0gr_{z_u}([\mfa],[\mfa])=
    0,\label{qunatity-zero}
  \end{equation}
   where $z_u$ is the closed loop joins $\mfa$ to $u\mfa$ for any $u\in\mathcal G_b(M)$.
   Then, there exists a constant $C$ such that for every $[\mfa],[\mfb]$ an each broken trajectory $[\breve\gamma]\in \breve M^+([\mfa],[\mfb])$, we have the energy bound
  \[\mathcal E^{top}_{\mathfrak q}(\gamma)\leq C+C(\iota([\mfa])-\iota([\mfb])).\]
\end{lemma}
\begin{pf}
The idea is similar to \cite[Lemma 16.4.4]{KM}.
  Let $[\breve \gamma]=([\breve\gamma_1],\cdots,[\breve \gamma_l])$ be a broken trajectory in $\breve M^+_z([\mfa],[\mfb])$ with $[\breve\gamma_i]\in \breve M_{z_i}([\mfa_{i-1}],[\mfa_i])$. The space $\breve M_{z_i}([\mfa_{i-1}],[\mfa_i])$ is non-empty, and it is manifold of dimension $1$, possibly  with boundary.
  We have that  $\dim (\breve M_{z_i}([\mfa_{i-1}],[\mfa_i]))$ is either $gr_{z_i}([\mfa_{i-1}],[\mfa_i])-1$ or $gr_{z_i}([\mfa_{i-1}],[\mfa_i])$. In either case, $gr_{z_i}([\mfa_{i-1}],[\mfa_i])\geq0$.  By adding the grading, it holds that $$gr_z([\mfa_0],[\mfa_l])\geq0.$$
  The energy $\mathcal E^a_{\mathfrak q}(\gamma)$ is  equal to $\mathcal E^{top}_{\mathfrak q}(z)$, defined by the twice of the change in $\L$ along any path $\tilde\zeta$ in $\mcC^\sigma(M,F,\mfs)$ whose image $\zeta$ in $\mcB^\sigma(M,F,\mfs)$ belongs to the class $z\in \pi_1(\mcB^\sigma(M,F,\mfs),[\mfa],[\mfb])$.
  By the condition, we have that the quantity
  $\mathcal E^{top}_{\mathfrak q} (w)+C_0gr_{w}([\mfa],[\mfa])$ depends only on $[\mfa]$ and $[\mfb]$ not on the homotopy class.
  This implies that the term \begin{equation}
    \mathcal E^{top}_{\mathfrak q} (w)+C_0(gr_{w}([\mfa],[\mfa])-2\iota(\mfa)+2\iota(\mfb))
    \label{quantity-bound}
  \end{equation}depends only on the critical points $[\alpha]=[\pi\mfa]$ and $[\beta]=[\pi\mfb]$. Since there are only finitely many critical points in $\mcB(M,F,\mfs)$, there is a constant  $C$ such that this quantity is at most $C$, which proves the lemma.
\end{pf}

\noindent{\bf Remark}:
In particular when $b^1_b=0$, the above lemma automatically holds.

In the $3$-manifold case, i.e. $F=0$, Kroheimer and Mrowka consider the quantity
\[\mathcal E_{\mathfrak q}(z)+4\pi^2 gr_z([\mfa],[\mfb]),\]
 where $z$ is a homotopy class   connecting $[\mfa],[\mfb]$. When $b^1 >0$,  the difference of the above quantity for two different homotopy classes is the class of a closed loop whose lift to the configuration space joins $\mfa$ to $u\mfa$, for some gauge action $u$. By the Atiyah-Singer theory on closed oriented manifolds, the difference is zero. For the basic Dirac operator $\Dirac_A$, Br\"uning, F. W. Kamber, K. Richardson gave an expression for its index \cite{BKR}. They showed that
  \[Ind(\Dirac)=\int_{\bar M_0/\bar F}A_{0,b}|\tilde{dx}|+\sum^r_{j=1}\beta(M_j),\]
  \[\beta(M_j)=\frac12\sum_\tau\frac1{n_\tau rank(W^\tau)}(-\eta(D^{S^+,\tau}_j)+h(D^{S^+,\tau}_j))\int_{\bar M_j/\bar F}A^{\tau}_{j,b}(x)|\tilde{dx}|,\]
  where the integrands $A_{0,b},~A^{\tau}_{j,b}(x)$ are similar to the Atiyah-Singer integrands and notations are explained in their paper.   Here we pose a question.

\noindent{\bf Question}: For $b^1_b>1$, under what topological condition, there is  a   constant $C_0$, such that for any non-degenerate critical point $\mfs$, we have $$C_0gr(\mfa,u\mfa)-\int_M[F_A]\wedge[u]\wedge\chi_F=0,$$ where $A$ is the connection component of $\mfa$ and $u\in \mathcal G_b(M)$.

%To get the energy bound, we still need the proposition below. \begin{prop}[c.f. {\cite[Proposition 16.1.4]{KM}}]\label{prop-16.1.4}    The space of broken trajectories $[\breve\gamma]\in\breve M^+([\mfa],[\mfb])$ with topology energy $\mathcal E_{\mathfrak q}(\breve\gamma)\leq C$ is compact.  \end{prop}
%Notice that the above proposition does not require that $b^1_b=0$, i.e. it holds for the general case.
%When $b^1_b=0$, we establish following proposition for the moduli space.
 % \begin{pf}  Consider a sequence $\{[\gamma^n]\}\subset M([\mfa],[\mfb])$ with bounded energy say $C$.We  choose a separating collection of neighborhoods $\mathcal U=\{U_\mfc\}$ for the interval $I=[-1,1]$ in the blown-up configuration space. Let $V_\mfc\subset \mcB^\sigma_{k}(M)$ be the corresponding neighborhoods.  Lemma  \ref{lemma-16.3.2} provides $\epsilon_1$ and $\epsilon_2$ such that for any trajectory $\gamma$ with $\mathcal E^I_{\mathfrak q}(\gamma)\leq \epsilon_1$ and $K^I_{\mathfrak q}\leq\epsilon_2$, we have that $\gamma_I\in U_\mfa$. Let $K_0$ be a uniform upper bound for $K(\gamma^n)$. We have that for each $n$, there are at most $2C/\epsilon+2K_0/\epsilon_2$ integers $p$ such that \[\tau_p\gamma^n\not\in \bigcup_{\mfc}U_{\mfc}.\] The rest of the proof is exactly same as Proposition \ref{cmpt-prop}.  \end{pf}
% By the arguments of the proof of \cite[Proposition 16.4.3]{KM}, we have that
 By Lemma \ref{lemma-16.4.4},
 we have the following proposition, which is analog to \cite[Proposition 16.4.1]{KM}.
 \begin{prop}\label{prop-16.4.1}
  Suppose  all the moduli spaces $M_z([\mfa],[\mfb])$ are regular and \eqref{qunatity-zero} holds. Then, there are finitely many homotopy classes $z$ for which  space $\check M^+([\mfa],[\mfb])$ is non-empty.
 \end{prop}

 \begin{prop}\label{prop-16.4.3}%[~{\cite[Proposition 16.4.3]{KM}}]
  Suppose  \eqref{qunatity-zero} holds and the  moduli space $M_z([\mfa],[\mfb])$ is regular. If $c_1(\mfs)=0\in H^2_b(M)$, then for a given $[\mfa]$ and $d\geq0$, there are finitely many pairs $([\mfb,]z$) such that the moduli space  $\check M^+([\mfa],[\mfb])$ is non-empty and of dimension $d$. If $c_1(\mfs)\neq0\in H^2_b(M)$, then for a generic perturbation there are only finitely many
triples $([\mfa],[\mfb],z)$ for which the moduli space  $\check M^+([\mfa],[\mfb])$ is non-empty.
\end{prop}
\begin{pf} The idea is no different to \cite[Proposition 16.4.3]{KM}, here we just give a sketch of the proof to the
case $c_1(\mfs)=0$. The functional $\L$ descends to a well-defined function on
  $\mcB^\sigma(M,F,\mfs)$, which is pulled back from $\mcB(M,F,\mfs)$.
   Since the image of critical points in $\mcB(M,F,\mfs)$ is finite, $\L$ takes finitely many values, the energy $\mathcal E^{top}_q$ of a trajectory is the twice of the drop of $\L$, so there is
a uniform bound on the energy of all solutions. The expression \eqref{quantity-bound} depends on $[\pi\mfa],~[\pi\mfb]$, which takes only finitely many values. By the condition that the dimension is bounded and $[\mfa]$ is fixed,   we have that $\iota([\mfb])$ is uniform bounded, which leaves finitely many choices.
\end{pf}

%The proof is to combine Proposition \ref{prop-16.4.3} and Lemma \ref{lemma-16.4.4}, the method is similar to the first part of the  proof of  \cite[Proposition 16.4.3]{KM}, here we omit it.

\section{Basic monopole Floer homologies on manifold with codimension $3$ foliation }

In this section, we  show the main result of this paper, i.e. to construct the basic monopole Floer homologies.

 \subsection{Basic Seiberg-Witten Floer homology for $b^1_b>1$}
  In this subsection, we assume that  $(M,F)$ is an oriented closed taut Riemannian foliation admitting a transverse \spinc structure and whose basic first deRham cohomology is nontrivial. Note that $\mathcal B(M,F,\mfs)$ is not simply connected in general, which implies that the index of critical points
  $$gr([\mfa],[\mfb])\in\mathbb Z$$
  might  not be
well defined. However we can still define relative gradings.
On the other hand,   the components of the moduli space $\mathcal M([\mfa],[\mfb])$, trajectories connecting the
critical points mod the gauge action, might have different dimensions corresponding
to the different lifts of $[\mfa]$ and $[\mfb]$, where $[\mfa]$ and $[\mfb]$ are the gauge equivalence classes
of the critical points.  Recall that
we   decompose the space of trajectories
$\mcM([\mfa],[\mfb])=\bigcup_{z\in\pi_1([\mfa],[\mfb])}
 \mcM_z([\mfa],[\mfb])$
as the union over the moduli spaces of different  relative homotopy classes, where $\pi_1([\mfa],[\mfb])$ denotes the homotopy class of path connecting the two critical points in the quotient space.

For one critical point $[\mfa]\in\mathcal B(M,F,\mfs)$, we might have different lifts in $\mathcal C(M,F,\mfs)$, say $\mfa$ and $u\mfa$, we can measure their spectral by the following index,
\[gr(\mfa,u\mfa)=Ind(\Dirac^+_{u}),\]
where $Ind(\Dirac^+_{u})$ denotes the index of the basic Dirac operator on the product space $(M\times S^1,F)$.

\begin{prop}
  The index $Ind(\Dirac^+_{u})$ defined above lifts to a homomorphism
  \[Ind:\pi_0(\mathcal G)\to \mathbb Z.\]
\end{prop}
\begin{pf}
  We need to show that for different critical points the index is unchanged, since it is clear to see that  for the same homotopy class, the index is well-defined. For another critical point $\mfb$, the connection difference is a one-form, which is a compact operator. Hence the index lifts to  a homomorphism.
\end{pf}

We define
\[d(\mathfrak s)=gcd(Ind:\pi_0(\mathcal G)\to \mathbb Z).\]
%Fix a critical point $\theta$(when $b^1_b=0$, we choice the unique reducible critical point as the reference point), for critical point $\mfa$, we set\[i(\mfa)=gr(\mfa,\theta)\in \mathbb Z_{d(\mathfrak s)}. \]
For two distinct irreducible critical points $\mfa$ and $\mfb$, we denote by $ \mcM^i([\mfa],[\mfb])$ the dimension $i$ component of ${\mathcal M}([\mfa],[\mfb])$. Let  $\check{\mathcal M}([\mfa],[\mfb])$ be the unparameterized space of ${\mathcal M}([\mfa],[\mfb])$, i.e. $\check{\mathcal M}([\mfa],[\mfb])={\mathcal M}([\mfa],[\mfb])/\mathbb R$. At the irreducible critical  points, the slice decomposition of Proposition \ref{prop-slice} holds.
By Theorem \ref{thm-19.5.4}, we have the following proposition.
\begin{prop}[c.f. {\cite[Corollary 3.1.24]{Wang}}]
  Suppose that $[\mfa_0],~[\mfa_2]$ are two irreducible critical points with the relative index $gr([\mfa_2],[\mfa_0])=2\bmod d(\mfs)$. Then the boundary of $\breve\mcM^2([\mfa_0],[\mfa_2])$ consists of union
  \[\bigcup_{[\mfa_1]\in Crit}\breve\mcM^1([\mfa_0],[\mfa_1])\times \breve\mcM^1([\mfa_1],[\mfa_2]),\]
  where $\mfa_1$ runs over critical points with $gr([\mfa_1],[\mfa_0])=1\bmod d(\mfs)$ and $Crit$ denotes the set of irreducible critical points in $\mcB(M,F,\mfs)$.
\end{prop}
We define the relative Floer complex $C(M)$ by the complex generated by the irreducible critical points of Chern-Simon-Dirac functional with grading given by the relative indices $\mathbb Z_{d(\mfs)}$ or $\mathbb Z$, i.e.
\[C(M)=\bigoplus_{\mfa\in Crit}\mathbb Z_2\mfa.\]
The boundary operator  of the complex is defined by
\[\partial :C(M)\to C(M),~\partial([\mfa])=\sum_{[\mfb]}\sharp\breve\mcM^1([\mfa],[\mfb]),\]
%\[\partial([\mfa])=\sum_{[\mfb]}\sharp\breve\mcM^1([\mfa],[\mfb]),\]
where $\sharp\breve\mcM^1([\mfa],[\mfb])\in \mathbb Z_2$ denotes the signed number of points in $\ \breve\mcM^1([\mfa],[\mfb])$ mod $2$.
\begin{lemma} \label{lemma-partial^2}
  $\partial^2=0$.
\end{lemma}

\begin{pf}
 By definition we have that
 $\partial^2 ([\mfa])=\sum_{[\mfb],} \sharp\breve\mcM^1([\mfa],[\mfb])\sharp\breve\mcM^1([\mfb],[\mfc])([\mfc])$, where $[\mfb]$ runs over the irreducible critical points with relative index $1$. By the above proposition, it is known that each term $\sharp\breve\mcM^1([\mfa],[\mfb])\sharp\breve\mcM^1([\mfb],[\mfc])([\mfc])$  is the sum of the number of oriented boundary points of a compact 1-dimensional manifold, which is zero.
\end{pf}

We define the basic Seiberg-Witten Floer homology  as  $HF(M,F,\mfs,\eta,g)=\ker(\partial)/Im(\partial)$, which is $\mathbb Z_{d(\mfs)}$-relative grading(or $\mathbb Z$-grading). % By the similar arguments, we have the following proposition.

\begin{prop}\label{prop-independent}
  Suppose that $(M,F)$ satisfies the Assumption \ref{assum-main-1}. Then, for $b^1_b>1$, we have that the relative Floer homology  is independent of the taut bundle-like metrics  and perturbations. We denote the basic Seiberg-Witten Floer homology group by $HF(M,F,\mfs)$.
\end{prop}

\begin{pf}
For a   bundle-like metric $g$, it corresponds to a triple
\[g~\leftrightarrow~(g_F,g_Q,s),\]
where $g_F$ is the leafwise restriction, $s$ corresponds to the decomposition
 \[s: Q\to TM,~\pi_Q\comp s=Id_Q\]
 and $g_Q$ is the transverse restriction.
 It is clear to see that the domain $\mathcal A_b\times \Gamma_b(M,S)$ and the Seiberg-Witten equations \eqref{eqn-SW-equation-1} are independent of the leafwise metric $g_F$ and the decomposition $s$.  For two distinct   leafwise metrics $g_F$ and $g'_F$ with the same $(p,l)$, the two Sobolev spaces $L^p_l$ and $L'^p_l$ are mutually  equivalent to each other. Therefore, we have that the basic Seiberg-Witten Floer homology groups are invariant under the leafwise metric. For two distinct decompositions $s$ and $s'$, we can apply the same argument, as the character form $\chi_F$ only depends on the leafwise metric $g_F$ and the decomposition $s$.

 The remaining part is to verify that the Floer homology group is independent of the generic basic perturbation and metric $g_Q$. The idea is exactly the same as Floer's original  proof \cite{Floer}.
  \end{pf}

\subsection{Basic monopole Floer homologies  for $b^1_b=0$}
 The purpose  of this subsection is to construct the basic monopole Floer homologies and show that they are independent of the perturbations and taut bundle-like metrics in a special case  $b^1_b=0$ .
%\subsection{Floer homologies }
We  define
the basic monopole Floer homologies $\overline{HM}_*(M,F,\mfs;\mathbb F)$, $\widecheck{HM}_*(M,F,\mfs;\mathbb F)$ and $\widehat{HM}_*(M,F,\mfs;\mathbb F)$ by the homologies of the chain complexes freely generated by $\bar C=C^s\cup C^u,~\check{C}=C^o\cup C^s,~\hat{C}=C^o\cup C^u$ respectively(see \cite[Section 22]{KM}), where $\mathbb F=\mathbb Z_2$.  The differentials on them are given in components as
\[\bar\partial =\left(\begin{array}{cc}
  \bar\partial^s_s&\bar\partial^u_s\\
  \bar\partial^s_u& \bar\partial^u_u\\
\end{array}\right),~
\breve\partial =\left(\begin{array}{cc}
  \partial^o_o&\partial^u_o\bar\partial^s_u\\
  \partial^o_s& \bar\partial^s_s+\partial^u_s\bar\partial^s_u\\
\end{array}\right),~
\hat\partial =\left(\begin{array}{cc}
  \partial^o_o&\partial^u_o \\
  \bar\partial^s_u\partial^o_s& \bar\partial^u_u+\bar\partial^s_u\partial^u_s\\
\end{array}\right).\]
The linear maps
\[\partial^o_o:C^o\to C^o,~\partial^o_s:C^o\to C^s,\]
\[\partial^u_o:C^u\to C^o,~\partial^u_s:C^u\to C^s\]
 are defined by the formula
 \[\partial^o_o[\mfa]=\sum_{[\mfb]\in C^o}\sharp \breve\mcM([\mfa],[\mfb])[\mfb],~[\mfa]\in C^o,\]
  where $\sharp \breve\mcM([\mfa],[\mfb])\in \mathbb F$ is the signed counting number, the other three are defined similarly. By considering the number $\sharp \breve\mcM^{red}([\mfa],[\mfb])$, we  similarly define the linear  maps
 \[\bar\partial^s_s:C^s\to C^s,~\bar\partial^s_u:C^s\to C^u,\]
 \[\bar\partial^u_s:C^u\to C^s,~\bar\partial^u_u:C^u\to C^u.\]

 When $b^1_b=0$, it is clear that
for a given $[\mfa]$, there are finitely many pairs $([\mfb],z)$ such that  the moduli space $ M_z([\mfa],[\mfb])$ is non-empty and of dimension $1$.
% \end{prop}
 %Combining  with Proposition \ref{prop-16.4.3}, one has the following proposition.
 \begin{prop}[c.f. {\cite[Proposition 22.1.4]{KM}}]
   \[\bar\partial^2=0,~\breve\partial^2=0,~\hat\partial^2=0.\]
 \end{prop}
 \begin{pf}
   The proof is by showing that $\bar\partial^2=0$, which is the same as the blow-down case(Lemma \ref{lemma-partial^2}), and the following identities
   \begin{enumerate}
     \item $\partial^o_o\partial^o_o+\partial^u_o\bar\partial^s_u \partial^o_s=0$;
     \item $\partial^o_s\partial^o_o+\bar\partial^s_s\partial^0_s+
     \partial^u_s\bar\partial^s_u\partial^o_u=0$;
     \item $\partial^o_o\partial^u_o+\partial^u_o\bar\partial^u_u+
     \partial^u_o\bar\partial^s_u\partial^u_s=0$;
     \item $\bar\partial^u_s+\partial^o_s\partial^u_o+\bar\partial^s_s\partial^u_s
         +\partial^u_s\bar\partial^u_u+\partial^u_s\bar\partial^s_u\partial^u_s=0$.
   \end{enumerate}
   Each of the four formulas is proved by considering a moduli space $\breve M_z([\mfa],[\mfb])$ of dimension $1$. By Theorem \ref{thm-19.5.4}, we can run the similar arguments of the proof in  \cite[Proposition 22.1.4]{KM}.
 \end{pf}

% We have the following proposition for the reducible moduli spaces.\begin{prop}[c.f. {\cite[Proposition 16.6.1]{KM}}]  Suppose that $ M^{red}([\mfa],[\mfb])$ is nonempty and of dimension $d$. Then the space of unparametrized broken reducible trajectories $\breve M^{red,+}([\mfa],[\mfb])$ is a compact $(d-1)$-dimensional space stratified by manifolds(see \cite[Definition 16.5.1]{KM}). The top stratum only consists  of $\breve M^{red}  ([\mfa],[\mfb])$. The $(d-l)$-dimensionalstratum consists of the spaces of unparametrized broken trajectories with$l$ factors:\[\breve M^{red}  ([\mfa_0],[\mfa_1])\times\cdots\times  \breve M^{red} ([\mfa_{l-1}],[\mfa_l]).\]\end{prop}

%\subsection{Indepedence}For the remaining part of this subsection, we show that the basic monopole Floer homology groups  are independent of the generic taut bundle-like metric and basic perturbation.
We give a grading for these homologies.
Let $\mathcal P$ be the space of the perturbations. We define $\mathbb J$ by the quotient of $\mathcal B^\sigma(M,F,\mfs)\times\mathcal P\times \mathbb Z/\sim$, where the equivalent relation $\sim$ is defined as follows(see \cite[Section 22.3]{KM}):
 for any two elements $([\mfa],\mathfrak q_1,m),~([\mfb],\mathfrak q_2,n)\in \mathcal B^\sigma(M,F,\mfs)\times\mathcal P\times \mathbb Z$, let $\zeta$ be a path joining $[\mfa]$ and $[\mfb]$ and $\mathfrak p$ be a path of perturbation joining $\mathfrak q_1$ and $\mathfrak q_2$. We have a Fredholm operator $P_{\zeta,\mathfrak p}$ as defined on \eqref{Fredholm-formula}, we say that $([\mfa],\mathfrak q_1,m)\sim ([\mfb],\mathfrak q_2,n)$, if there is a path $\zeta$ such that
 \[Ind(P_{\zeta,\mathfrak p})=n-m.\]
 The map $([\mfa],\mathfrak q,m)\mapsto ([\mfa],\mathfrak q,m+1)$ descends to $\mathbb J$, and raises to an action of $\mathbb Z$.

  Note that  the  above construction of the index set $\mathbb J$ is  also available when $b^1_b>0$.  Let $\mathfrak q$ be a fixed admissible perturbation, for a critical point $[\mfa]$, we define its grading by
 \[gr([\mfa])=([\mfa],\mathfrak q,0)/\sim\in \mathbb J.\]
 For reducible critical points, we define the modified grading by
 \[\bar{gr}([\mfa])=\begin{cases}
   gr([\mfa])&[\mfa]\in C^s\\
   gr([\mfb])-1&[\mfa]\in C^u.
 \end{cases}\]

% \begin{lemma}%[c.f. {\cite[Lemma 22.3.2, Lemma 22.3.3]{KM}}]    When $b^1_b=0$, the $\mathbb Z$ action on $\mathbb J$ is transitive and free. For each $j\in\mathbb J$, the free abelian groups $\breve C_j,~\hat C_j$ and $\bar{C}_j$ are finitely generated. \end{lemma}

 We  show that   the basic monopole Floer homologies is independent of  the generic perturbations and    bundle-like metrics.
Let $W=[0,1]\times M$ and $W^*=(-\infty,0]\times M\cup W\cup [1,\infty)\times M$. To tell the distinguish, we denote the left boundary
$\{0\}\times M$ with metric by $Y_-$ and perturbation and the right boundary $\{1\}\times M$ with another metric and perturbation by $Y_+$. We consider the moduli space $M([\mfa],W^*,[\mfb])$, as defined in \cite[Section 25]{KM}. Using broken trajectories, we denote its compactification by $M^+([\mfa],W^*,[\mfb])$.

Fix a positive  integer $d_0$, we consider a pair  $([\mfa],[\mfb])$ for which the moduli space $M([\mfa],W^*,[\mfb])$ or $M^{red}([\mfa],W^*,[\mfb])$ has dimension $d_0$ at most. To prove the independence of the metrics, we need to define the homomorphism  maps by the trivial cobordism, which are given by
counting the number of solutions in the zero-dimensional moduli spaces. % We can compute the cohomology of $\mathcal B^\sigma_{k,loc}$ using \v{C}ech cochains carried by open covers $\mathcal U$ that are transverse to all these moduli spaces(See \cite[Lemma 21.2.1]{KM}). Let  \[u\in C^d(\mathcal U;\mathbb F)\]  be a \v{C}ech cochain with $d\leq d_0$. If $M([\mfa],W^*,[\mfb])$ has dimension $d$, then there is a well-defined evaluation(see \cite[21.3]{KM})  \[\langle u,M([\mfa],W^*,[\mfb])\rangle\]  by counting solutions in $d$ dimensional moduli spaces.
We define linear operators,
  \[m^o_o:%C^0(\mathcal U;\mathbb F)\otimes
   C^o_*(Y_-)\to C^o_*(Y_+),~m^o_s:%C^0(\mathcal U;\mathbb F)\otimes
  C^o_*(Y_-)\to C^s_*(Y_+)\]
 % \[m^o_s:%C^0(\mathcal U;\mathbb F)\otimes  C^o_*(Y_-)\to C^s_*(Y_+)\]
 \[m^u_o:%C^0(\mathcal U;\mathbb F)\otimes
 C^u_*(Y_-)\to C^o_*(Y_+),~m^u_s:%C^0(\mathcal U;\mathbb F)\otimes
  C^u_*(Y_-)\to C^s_*(Y_+)\]
 %\[m^u_s:%C^0(\mathcal U;\mathbb F)\otimes  C^u_*(Y_-)\to C^s_*(Y_+)\]
 by
 \[m^o_o(-)=\sum_{[\mfa]\in C^o(Y_-)}\sum_{[\mfb]\in C^o(Y_+)}
 \sharp M([\mfa],W^*,[\mfb]),\]
 for the first one and by the similar formulas  for the other three. Similarly, we define operators on the reducible part of the Floer complexes: we have an operator
 \[\bar m:%C^0(\mathcal U;\mathbb F)\times
 \bar C_*(Y_-)\to \bar C_*(Y_+)\]
 \[\bar m=\left(\begin{array}{cc}
   \bar m^s_s&\bar m^u_s\\
   \bar m^s_u&\bar m^u_u\\
 \end{array}\right)\]
 where $\bar m^s_s(-)=\sum_{[\mfa]\in C^s(Y_-)}\sum_{[\mfb]\in C^s(Y_+)}
 \sharp M([\mfa],W^*,[\mfb]) $, and the other three entries  are defined similarly. On $\check{C}_*$, we define
 \[\check m:%C^0(\mathcal U;\mathbb F)\otimes
 \check{C}_*(Y_-)
 \to \check{C}_*(Y_+)\]
 by the formula
 \[\check m=\left(\begin{array}{cc}
   m^o_o&m^u_o\bar\partial^s_u(Y_-)+\partial^u_o(Y_+)\bar m^s_u\\
   m^o_s&\bar m^s_s+m^u_s\bar\partial^s_u(Y_-)+\partial^u_s(Y_+)\bar m^s_u\\
 \end{array}\right),\]
 where  $\partial^u_o(Y_\pm)$, for example, denotes the operator $\partial^u_o$ on $Y_\pm$. On $\hat C_*$, we define
 \[\hat m: %C^0(\mathcal U;\mathbb F)\otimes
  \hat C_*(Y_-)\to \hat C_*(Y_+)\]
 by the formula
 \[\hat m=\left(\begin{array}{cc}
   m^o_o&m^u_o\\
   \bar m^s_u\partial^o_s(Y_-)+\bar\partial^s_u(Y_+)m^o_s&
   \bar m^u_s+\bar m^s_u\partial^u_s(Y_-)
 \end{array}\right).\]

%By {\cite[Lemma 25.3.6]{KM}}, we have the following lemma.
By considering the zero-dimension moduli space, we have the following proposition.

 \begin{prop}\label{prop-cobordism}
   The operators $\breve m,~\hat m$ and $\bar m$ satisfy the identities:
   \[\begin{cases}
     \breve\partial(Y_+)\breve m_{-+}=\breve m_{-+}( \partial(Y_-)),\\
     \hat\partial(Y_+)\hat m_{-+}=\hat m_{-+}( \bar\partial(Y_-)),\\
     \bar\partial(Y_+)\bar m_{-+}=\bar m_{-+}(  \bar \partial(Y_-)).
   \end{cases}\]
   %for $u\in C^d(\mathcal U)$ and $\delta$ denotes the \v{C}ech coboundary map $\delta: C^d(\mathcal U)\to C^{d-1}(\mathcal U)$.
   In particular, we give rise to the operators
   \[\begin{cases}
     \breve m_{-+} :  \widecheck{HM}_*(Y_-)\to \widecheck{HM}_*(Y_+)\\
     \hat m_{-+} : \widehat{HM}_*(Y_-)\to \widehat{HM}_*(Y_+)\\
     \bar m_{-+}:\overline{HM}_*(Y_-)\to \overline{HM}_*(Y_+).
   \end{cases}\]
   Moreover, the above operators only depend on the data of $Y_-$ and $Y_+$.
   %for all open cover $\mathcal U$ of $\mcB^\tau_{k,loc}(W^*)$ transverse to all the moduli spaces of dimension less than $d_0$.
 \end{prop}
 Note that since we focus on the zero dimension part, the proof is much easier than {\cite[Proposition 25.3.8]{KM}}.
 %By the previous argument, we have the following theorem.
% \begin{thm}[ {\cite[Proposition 25.3.8]{KM}}]   Let $(g_0,\mathfrak q_0)$ and $(g_1,\mathfrak q_1)$ be two pairs of metric  and perturbation on $W$, which are isometric in a collar of the boundary to the same cylindrical metric. Suppose that the corresponding moduli spaces on $W^*$ are regular in both cases. Then, there is an operator   \[\breve K:C^d(\mathcal U)\otimes \breve C_*(Y_-)\to \breve C_*(Y_+)\]   satisfying the identity   \[(-1)^d\breve\partial \breve K(u\otimes-)=-\breve K(\delta u\otimes-)+\breve K(u\otimes\breve\partial-)+(-1)^d\breve m(0)(m\otimes-)-   (-1)^d\breve m(1)(u\otimes-).\]   Similarly, we have the other two chain maps $\hat K$ and $\bar K$. \end{thm}

The last step is to prove the composition law. Let $Y_-$, $Y_0$ and $Y_+$ be  the same $(M,F)$ with three metrics and basic perturbations, let $W_1$ be the cobordism from  $Y_-$ to $Y_0$ such that near each collar, the metric of $W_1$ is the product metric and $W_2$ be the cobordism from $Y_0$ to $Y_+$ with the same condition on the metric. %We set \[W(S)=W_1\cup([0,S]\times Y_0)\cup W_2.\] Let $\mathfrak q^-$, $\mathfrak q^0$ and $\mathfrak q^+$ be the admissible tame perturbations on $(M,F)$. By cut-off functions, we extend these to $t$-dependent perturbations $\mathfrak q$ supported in the collars of the four components of the boundaries of $W_-$ and $W_+$. When $W_-$ and $W_+$ are joined together to form the composite cobordism $W$, the perturbations near the two copies $Y_0$ match. % so we have a well-defined perturbation of the  on $W$, compatible with the restriction maps to the two components. Let $[\mfa_-]$ and $[\mfa_+]$ be the critical points on $Y_-$ and $Y_+$ respectively, we consider a moduli space \[M([\mfa_-],W(S)^*,[\mfa_+])\] on the cylindrical-end  manifold $W(S)^*$. As $S$ varies, the moduli space parametrized by $S\in[0,\infty)$, \[M([\mfa_-],[\mfa_+])=\bigcup_{S\in[0,\infty)}\{S\}\times M([\mfa_-],W(S)^*,[\mfa_+]).\]

Repeat the same argument  in \cite[Section 26.1]{KM}. We have the composition law below for the cobordisms.

\begin{prop}[c.f. { \cite[Proposition 26.1.2]{KM}}]
Let $(M,F)$  satisfy the Assumption \ref{assum-main-1}. Fix a transverse \spinc structure.
  Let $Y_-,~Y_0,~Y_+$ be three data of   bundle-like metrics and basic perturbations, and let $W_{-0}$ be the cobordism from $Y_-$ to $Y_0$, $W_{0+}$ be the cobordism from $Y_0$ to $Y_+$ and $W_{-+}$ be the composition of the $W_{-0}$ and $W_{0+}$. Suppose that $m_{-0}$, $m_{0+}$ and $m_{-+}$
 are the operators in Proposition \ref{prop-cobordism}. Then  we have that
 \[m_{0+}\comp m_{-0}=m_{-+}.\]
 \end{prop}
% Note that the proof is much easier than { \cite[Proposition 26.1.2]{KM}}, since we only consider the zero-dimension moduli space.
 The above proposition implies the corollary below.

\begin{cor}
  The monopole Floer homologies are independent of the generic choice of the perturbation and  the  bundle-like metric, which are denoted by $\widecheck{HM}_*(M,F,\mfs;\mathbb F)$, $\widehat{HM}_*(M,F,\mfs;\mathbb F)$ and $\overline{HM}_*(M,F,\mfs;\mathbb F)$.
\end{cor}
 \begin{prop}[c.f. {\cite[Proposition 22.2.1]{KM}}]
    Let $(M,F)$ satisfy the Assumption \ref{assum-main-1}. Then,  there is an exact sequence
   \[...\overline{HM}_*(M,F,\mfs;\mathbb F)\overset{i_*}\to
   \widecheck{HM}_*(M,F,\mfs;\mathbb F)\overset{j_*}\to \widehat{HM}_*(M,F,\mfs;\mathbb F)\overset{p_*}\to \overline{HM}_*(M,F,\mfs;\mathbb F)
   \overset{i_*}\to...\]
   in which the maps $i_*,~j_*$ and $p_*$ arise from the chain-maps
   \[i:\bar C\to \breve{C},~j:\breve{C}\to \hat{C},~
   p:\hat{C}\to \breve{C},\]
   which are defined by
   \[i=\left(\begin{array}{cc}
     0&-\partial^u_o\\
     1&-\partial^u_s\\
   \end{array}\right),~
   j=\left(\begin{array}{cc}
     1&0\\
     0&-\bar\partial^s_u\\
   \end{array}\right),~
   p=\left(\begin{array}{cc}
     \partial^o_s&\partial^u_s\\
     0&1\\
   \end{array}\right).\]
   Here $i$ and $j$ are genuine chain maps, however $p$ is a anti chain map, i.e.  $p\hat\partial+\hat\partial p=0$.
 \end{prop}
 We review the completion of graded groups, c.f. \cite[Defintion 3.1.3]{KM}.
 Let $G_*$ be an abelian group graded by the set $\mathbb J$ equipped with a  $\mathbb Z$-action. Let $O_a(a\in A)$ be the set  of free $\mathbb Z$-orbits in $\mathbb J$ and fix an element $j_a\in O_a$ for each $a$. Consider the subgroups
 \[G_*[n]=\bigoplus_a\bigoplus_{m\geq n} G_{j_a-m},\]
 which form a decreasing filtration of $G_*$. We define the negative completion of $G_*$ as the topological group $G_\bullet\supset G_*$ obtained by completing with respect to this filtration. We define the negative completions
 \[\widecheck{HM}_\bullet(M,F,\mfs;\mathbb F),~\widehat{HM}_\bullet(M,F,\mfs;\mathbb F),~\overline{HM}_\bullet(M,F,\mfs;\mathbb F),\]
 of the basic monopole Floer homologies defined as the above. If we want to consider all transverse \spinc
structures at the same time, we need to consider the completed basic monopole Floer homology
\[\widecheck{HM}_\bullet(M,F;\mathbb F)=\bigoplus_\mfs\widecheck{HM}_\bullet(M,F,\mfs;\mathbb F).\]
We make similar definitions for $\widehat{HM}_\bullet(M,F;\mathbb F)$ and $ \overline{HM}_\bullet(M,F;\mathbb F)$.

 By the previous argument, we construct the basic monopole Floer homology groups and prove that these   homology groups are independent of the basic perturbations and bundle-like metrics.
 
 {\bf Remark:
}
Even though what we did is to construct the basic monopole Floer homologies with $\mathbb Z_2$ coefficient, as argued in \cite[Section 20.4-20.5]{KM}, one can define the basic monopole Floer homology with integer coefficient or any commutative ring.

\subsection{Basic monopole Floer homologies   for $b^1_b>0$}

In this subsection, we construct the basic monopole Floer homologies %$\widecheck{HM}$, $\widehat{HM}$ and $\overline{HM}$
for $b^1_b>0$ with the Novikov ring.
We recall the notion of local system $\Gamma$ over a topological space $X$.

\begin{defi}
  A local system on a  topological space  $X$, is a system to distribute abelian groups $\{\Gamma_a\}$ for each point $a\in X$, such that for each relative homotopy class of paths $z$ from $a$ to $b$, there is  an isomorphism
  \[\Gamma(z):\Gamma_a\to \Gamma_b\] satisfying the composition law for the composition of two paths.
\end{defi}
%Suppose that $X$ is a manifold and $f$ is Morse function on $X$, we define a graded abelian group $C_*(f,\Gamma)=\bigoplus_{a\in Crit(f)}\Gamma_a$.
We review the classical results of \cite[Section 22, 29, 30]{KM}. %let $Y$ be a closed oriented $3$-manifold,
We   choose $\Gamma$ to be a local system of abelian groups on $\mcB^\sigma_b(M,F,\mfs)$, such that to each point $[\mfa]\in \mcB^\sigma_b(M,F,\mfs)$ there is an associated group $\Gamma[\mfa]$ and to each homotopy class $z$ of the paths from $[\mfa]$ to $[\mfb]$, there is an associated isomorphism $\Gamma(z):\Gamma[\mfa]\to\Gamma[\mfb]$. %We define \[C^o(\Gamma)=\bigoplus_{[\mfa]\in C^o}\mathbb Z\Lambda[\mfa]\otimes \Gamma[\mfa],\]\[C^s(\Gamma)=\bigoplus_{[\mfa]\in C^s}\mathbb Z\Lambda[\mfa]\otimes \Gamma[\mfa],\]\[C^u(\Gamma)=\bigoplus_{[\mfa]\in C^u}\mathbb Z\Lambda[\mfa]\otimes \Gamma[\mfa],\] where $\Lambda[\mfa]$ is a two-element set which is defined in \cite[Definition 20.3.1]{KM}.

%\begin{defi}  Let $G_*$ be a graded abelian group graded by a set $\mathbb J$,  $O_\alpha$($\alpha\in A$) be  the free $\mathbb Z$-orbits in $\mathbb J$, and an element $j_\alpha\in O_\alpha$ be given for each $\alpha$. Let $$G_*[n]=\bigoplus_\alpha\bigoplus_{m\geq n}G_{j_\alpha-m}.$$  This is a decreasing filtration of $G_*$, i.e. $G_*[n]\supset G_*[n-1]$. We define the negative completion of $G_*$ to the topological group $G_\cdot\supset G_*$ obtained by completing with respect to this fibration.\end{defi}
By using the above local system $\Gamma$, the boundary maps
% $\partial^o_o$ a
are well-defined. For instance, we consider
\[C^o(Y,\mfs,c,\Gamma)=\bigoplus_{[\mfa]}%\mathbb Z\Lambda[\mfa]\otimes
\Gamma[\mfa], \] where $[\mfa]$ denotes the irreducible critical point.
We define the partial by
\[\partial^o_o =\sum_{[\mfa]}\sum_{[\mfb]}\sum_z\sum_{[\gamma]}\in \check M_z([\mfa],[\mfb])\otimes \Gamma(z),\]
where the sum is over all the moduli space $M_z([\mfa],[\mfb])$ with dimension $1$ and $[\mfb]$ denotes the irreducible critical point. The contribution for a given pair of critical points takes the form
\begin{equation}
  \sum_zn_z\Gamma(z),\label{eqn-summand}
\end{equation}
where $z$ runs through all relative homotopy classes satisfying the conditions
$gr_z ([\mfa],[\mfb]) = 1 $.
Before proceeding, we review some definitions and notions which are given in \cite[Section 30]{KM}.
\begin{defi}[c.f. {\cite[Definition 30.2.1]{KM}}]
  Let $\mathcal E^{top}_{\mathfrak q}$ be a corresponding perturbation of the topological energy. A subset $S\subset \pi_1([\mfa],[\mfb])$ is called $c$-finite, where $\pi_1([\mfa],[\mfb])$ denotes the homotopy classes of paths joining $[\mfa]$ and $[\mfb]$ in $\mcB^\sigma(M,F,\mfs)$, if the following conditions are satisfied:
  \begin{itemize}
    \item for all $C$, $S \cap \{z|\mathcal E^{top}_{\mathfrak q}(z)\leq C\}$ is finite;
    \item there exists $d\geq0$ such that $|gr_z([\mfa],[\mfb])|\leq d$ for all $z\in S$.
  \end{itemize}
\end{defi}

We consider a local system of complete topological abelian groups $\Gamma$ on $\mcB^\sigma(M,F,\mfs)$, i.e. each $\Gamma[\mfa]$ is a complete topological group and each homomorphism $\Gamma(z):\Gamma[\mfa]\to\Gamma[\mfb]$ is continuous.  Assume that $0\in\Gamma[\mfa]$ has a neighborhood basis consisting of subgroups, such that $\Gamma[\mfa]$ is a complete filtered group, which is filtered by the open subgroups. Let $Hom(\Gamma[\mfa],\Gamma[\mfb])$ be the group of continuous homomorphisms, equipped with the compact-open topology. A neighborhood basis for $0$ in $Hom(\Gamma[\mfa],\Gamma[\mfb])$ consists of subgroups
\[\Omega(N,V)=\{k:\Gamma[\mfa]\to\Gamma[\mfb]|k(N)\subset V\}, \]
where $N$ runs over all precompact subsets of $\Gamma[\mfa]$ and $V$ runs all open subgroups of $\Gamma[\mfb]$. Note that a subset $N\subset \Gamma[\mfa]$ is precompact if and only if $(N+U)/U$ is finite for all open subgroups $U$ of $\Gamma[\mfa]$.

\begin{defi}
  A countable series $\sum_{k\in K}k$ of $Hom(\Gamma[\mfa],\Gamma[\mfb])$ is said to be equicontinuous, if for each open subgroup $U$ $\subset\Gamma[\mfb]$, there exists an open subgroup $V$ such that $k(V)\subset U$ for each $k\in K$.
\end{defi}

\begin{defi}[c.f. {\cite[Definition 30.2.2]{KM}}]
   A local system of complete filtered abelian groups $\Gamma$ is called $c$-complete, if it satisfies the following properties for each $[\mfa],~[\mfb]$:
  \begin{itemize}
    \item for any $c$-finite set $S\subset \pi_1(\mcB^\sigma,[\mfa],[\mfb])$, the set $\{\Gamma(z)|z\in S\}\subset Hom(\Gamma[\mfa],\Gamma[\mfb])$ is equicontinuous;
    \item  for any $c$-finite set $S\subset \pi_1(\mcB^\sigma,[\mfa],[\mfb])$,  $\Gamma(z)$ converges to zero as $z$ runs through  $S$ in the compact-open topology.
  \end{itemize}
\end{defi}
Notice that there might be infinitely many nonzero terms in the form \eqref{eqn-summand},
we set the support as
\[supp(n)=\{z| ~n_z\neq0\}.\]
By the definition of c-complete, we have that
\[supp(n)\cap \{z|\mathcal E^{top}_{\omega,\mathfrak q}(z)\leq C\}\] is finite. Using the completeness of the local system, we have that the form
the form \eqref{eqn-summand} is convergent. Similarly, one can verify that the maps $\check\partial$, $\hat\partial$ and $\bar\partial$ are well-defined.
Combining with the equicontinuous property of the local system $\Gamma$, the proofs of $\check\partial^2$, $\hat\partial^2$ and $\bar\partial^2$ go through as the non-exact perturbation of 3 manifold case(see \cite[Section 30.2]{KM}).

\noindent We give an example of such a local system, e.g. a Novikov ring \cite{Novikov}. We have a homomorphism
\[\mathcal E^{top}:\pi_1(\mcB (M,F,\mfs))\to \mathbb R,~z\mapsto \mathcal E^{top}(z),\] where $\mathcal E^{top}(z)$ denotes the difference of the Chern-Simons-Dirac functional  between the different representatives of the quotient point. Since $\pi_1(\mcB^\sigma(M,F,\mfs))\cong \mathbb Z^{b^1_b}$, we can choose  a basis $\{z_i\}_{1\leq i \leq b^1_b}$, such that each element $z$ can be written as
$z=k_1z_1+\cdots k_{b^1_b}z_{b^1_b},$
where $k_i\in\mathbb Z$ for $i=1,\cdots, b^1_b$ and $\mathcal E^{top}(z_i)\geq 0$. Moreover, we may assume that for  $1\leq i\leq l$,  $\mathcal E^{top}(z_i)>0$.   It is not hard to see that such a basis $\{z_i\}_{1\leq i \leq l}$ is independent of the metric $g$.
Choosing a commutative ring $R$(e.g. $\mathbb Z_2$), we define $R[t,,t^{-1}]$ by
 \[R[t,t^{-1}]=\{\sum _{-k\leq i\leq K} r_it^i|\mbox{ only for finitely many }i,~r_i\ne0\}.\]For $k\in\mathbb Z$, let $U_{-k}$ be the $R$-module spanned by the generators $t^i$,  satisfying $i\leq -k$. Using these as open neighborhoods of $0$, we form the completion $\bar R[t,t^{-1}]$, i.e.  each element is of the form $$\sum^C_{i= -\infty} r_it^i.$$  We define a local system by taking at $[\mfa_0]$ to be $\bar {R}[t,t^{-1}]$, and specifying that for each closed loop $z$ based at $[\mfa_0]$, the automorphism $\Gamma(z)$ be the multiplication by   $t^{-(k_1+\cdots k_{l})}$ for $z=k_1z_1+\cdots  k_lz_l+\cdots k_{b^1_b}z_{b^1_b}$. %This is a $c$-complete local system(see \cite[Definition 30.2.2]{KM}).

 Repeat the parallel arguments of previous subsection, together with Proposition \ref{prop-16.4.1} and Proposition \ref{prop-16.4.3}, we have that:
 \begin{thm}
  Let $(M,F)$ satisfy the Assumption \ref{assum-main-1}. For a complete local system $\Gamma$, e.g. a Novikov ring, we can construct the basic monopole Floer homologies. Moreover,
 The monopole Floer homologies are independent of the generic choices of the perturbation and  the  bundle-like metric, which are denoted by $\widecheck{HM}_*(M,F,\mfs;\Gamma)$, $\widehat{HM}_*(M,F,\mfs;\Gamma)$ and $\overline{HM}_*(M,F,\mfs;\Gamma)$. Moreover, if \eqref{qunatity-zero} holds, then for any local system, we have the well-defined basic monopole Floer homologies.
 \end{thm}

 In general, we   consider the (non-exact)perturbed basic Chern-Simons-Dirac functional defined as below:
 given a class $c\in H^2_b(M)$, we write
 \[\L_\omega(A,\Psi)=\L(A,\Psi)-\frac12\int_M(A^t-A^t_0)\wedge\omega\wedge\chi_F,\]
 where $\omega\in \frac{2\pi}ic$.
% We apply the non-exact perturbation method to construct the basic monopole Floer homology groups  for $(Y,\mfs)$.  Let $\omega$ be a closed imaginary-valued two-form belonging to the class $[\frac{2\pi}{i}c]$, we define $f_\omega(B,\Psi)=\int_Y(B-B_0)\wedge \omega$ and  \[\L_\omega(B,\Psi)=\L(B,\Psi)+f_\omega(B,\Psi).\] By the straightforward calculation, we have that \[\L_\omega(B,\Psi)-L_\omega(u(B,\Psi))=-4\pi^2\langle c_1(\mfs)-c ,[u]\rangle. \] Hence, we have the relationship \[\L_\omega(B,\Psi)-L_\omega(u(B,\Psi)+4\pi^2 gr(\mfa,u\mfa)=0,\] where $(B,\Psi)$ denotes the pull-back image of a critical point $\mfa$ on $\mcB(Y,\mfs)$.
 It is known that  a critical point $(A,s,\psi)$ in the blow-up model $\mcC^\sigma(M,F,\mfs)$ is defined by
\begin{equation}
  \begin{cases}
    \frac12\bar*(F_{A^t}-\omega)=s^2q(\psi),\\
    \Dirac_A\psi=0.
  \end{cases}\label{eqn-non-exact-SW-29.2}
\end{equation}
The corresponding perturbed  equations for $(A,s,\phi)\in\mcC^\tau(\mathbb R\times M)$ are defined by
\begin{equation}
  \begin{cases}
    \frac12(F^+_{A^t}-\omega^+)=s^2q(\phi),\\
    \frac{d}{dt}s+\Lambda(A,s,\phi)s=0,\\
    \Dirac^+_A\phi-\Lambda(A,s,\phi)\phi=0.
  \end{cases}\label{eqn-non-exact-SW-29.3}
\end{equation}
%We set   $\mathcal E^{top}_\omega=2(\L_\omega(t_1)-\L_\omega(t_2))$, for $Z=[t_1,t_2]\times M$, and   $\omega$ be the pull-back on $Z$, it is known that $\mathcal E^{an}_\omega(A,\Phi)=\mathcal E^{top}_\omega(A,\Phi)+\|\mathfrak F_\omega(A,\Phi)\|^2_{L^2}$, where
%\begin{eqnarray*}  \mathcal E^{an}_\omega(A,\Phi)=\frac14\int_Z|F_{A^t}-4\omega|^2+\int_Z|\nabla_A\Phi|+\frac14  \int_{Z}(|\Phi|^2+Scal/2)^2-\int_{Z}Scal^2/16+2\int_Z\langle \Phi,\omega\Phi\rangle.\end{eqnarray*}
%For a general non-exact perturbation $grad(f_\omega)+\mathfrak q$ on $Y$, we can write $\L=\L_\omega+g$, where $g$ is gauge-invariant and has gradient $\mathfrak q$, we define the topological energy\[\mathcal E^{top}_{\omega,\mathfrak q}(A,\Phi)=\mathcal E^{top}_\omega+2g(A,\Phi).\]
Applying the standard argument of the manifold case, we have the following lemma.

\begin{lemma}  If $c\neq c_1(\mfs)$, then there are no reducible critical points of $\L_\omega$.\end{lemma}
With the non-exact perturbation,  the space of broken trajectories space $\check M^+_z([\mfa],[\mfb])$ can be defined as the manner of exact perturbation. %Following the argument of the proof of Proposition \cite[Proposition 16.4.3]{KM}, one can establish the following proposition.

%\begin{prop}   Suppose that all the moduli spaces $M_z([\mfa],[\mfb])$ for the   perturbation $\mathfrak q$ are regular.  Then, for a given $[\mfa]$ and $d_0\geq0$, there are only finitely many pairs $([\mfb],z)$ for which the moduli space $\check{M}^+_z([\mfa],[\mfb])$ is non-empty and has dimension no more than $d_0$. In additional, if $c\neq c_1({\mfs})$ and the perturbation has  been chosen so that there are no reducible solutions, then for a given $[\mfa]$, there are only finitely many pairs $([\mfb],z)$ for which the moduli space $\check{M}^+_z([\mfa],[\mfb])$ is non-empty. \end{prop}

Following the same strategy of the construction in Section 6(or see \cite[Section 20-Section 26]{KM}),% we have the basic monopole Floer homologies \[\widecheck{HM}_*(M,F,\mfs,\omega;\Gamma), ~ \widehat{HM}_*(M,F,\mfs ,\omega;\Gamma),~ \overline{HM}_*(M,F,\mfs,\omega;\Gamma).\]By the argument of the previous subsection or \cite[Section 26]{KM},
 we have the following theorem. % to say that the above homology groups are independent of the choice of the representatives $\omega$, but dependent of the class $c$.

\begin{thm}
  Let $\Gamma$ be a complete local system, e.g. a Novikov ring, and $\L_\omega$ be a non-exact perturbation for the Chern-Simons-Dirac functional defined as above. Then we have  the basic monopole Floer homologies
  \[\widecheck{HM}_*(M,F,\mfs,c;\Gamma), ~ \widehat{HM}_*(M,F,\mfs ,c;\Gamma),~ \overline{HM}_*(M,F,\mfs,c;\Gamma),\]
  where $c\in H^2_b(M)$. These homologies depend only on the isomorphism class of the \spinc structure $\mfs$, $c$ and $(M,F)$, however they are independent of the bundle-like metrics and  the basic perturbations.
\end{thm}

  At the end of this subsection,
  we give a necessary condition to avoid the complete local system.
% Before we proceed, we give a review of classical non-exact perturbation(see \cite[Section 29-Section 30]{KM}).
%Let $Y$ be a closed oriented three manifold and $\mfs$ be a \spinc structure. Choosing a class $c\in H^2(Y,\mathbb R)$, we consider the (non-exact)perturbed Chern-Simons-Dirac functional, \[\L_\omega(B,\Psi)=\L(B,\Psi)-\frac12\int_Y(B-B_0)\wedge\omega,\] where $\omega\in \frac{2\pi}ic$ and $(B,\Psi)\in \mcC(Y,\mfs)$. By the argument of \cite[Section 30]{KM}, we can define the basic monopole Floer homologies with Novikov ring as coefficient.

% \begin{thm}[c.f. {\cite[Section 29.2]{KM}}]\label{thm-perturbation-no-need-Novikov}   Let $(Y,\mfs,c)$ be given as above. Suppose that there is a constant $t\in\mathbb R$ such that   \[2\pi^2c_1(\mfs)-c=t2\pi^2c_1(\mfs).\]   Then, the Novikov ring is not necessary, i.e. for any local coefficient, one can define the basic monopole Floer homologies. \end{thm}
% The proof is given in \cite[Section 29.2]{KM}.

% Similar to Theorem \ref{thm-perturbation-no-need-Novikov}, we have the following theorem.

% \begin{thm}   Let $Y$ be a closed oriented $3$-orbifold defined as above, and $(\mfs,c)$ as above. Suppose that there is a some constant $t$ such that $$2\pi^2c_1(\mfs)-c=2\pi^2t(c_1(\mfs)-c(\mfs,Y)).$$ Then, the Novikov ring coefficient is not necessary. \end{thm}
 \begin{thm}
   Let $(M,F,\mfs,c)$ be as above. Let $g$ be a   bundle like metric and $\chi_F$ be the character form of the foliation. Suppose that there is a constant $t$ such that the identity holds
   \[-\int_M(c_1(\mfs)-c)\wedge[u]\wedge\chi_F+t\cdot gr(\mfa,u\mfa)=0,\]
   for a non-degenerate critical point. Then, with any local coefficient $\Gamma$ we have the basic monopole Floer homologies
   \[\overline{HM}_*(M,F,\mfs,c;\Gamma),~\widehat{HM}_*(M,F,\mfs,c;\Gamma),~
   \widecheck{HM}_*(M,F,\mfs,c;\Gamma).\]
 \end{thm}
 \begin{pf}
   Here we give a sketch of the proof. The idea is to show that $\sum_zn_z\Gamma(z)$ is of finitely many sum, for each $z\in M_z([\mfa],[\mfa])$, where $M_z([\mfa],[\mfb])$ is a moduli space of dimension $1$ and $[\mfa],~[\mfb]$ are two regular critical points. It is sufficient to show a foliated version of %Kronheimer Mrowka's proposition
   \cite[Proposition 29.2.1]{KM}, which is stated as below.
 \end{pf}

 \begin{prop}\label{prop-29.2.1}
   Let $(M,F,\mfs,c)$ be as above, let $g$ be a  bundle like metric and $\chi_F$ be the character form of the foliation. Suppose that there is a constant $t$ such that the identity holds
   \[-\int_M(c_1(\mfs)-c)\wedge[u]\wedge\chi_F+t\cdot gr(\mfa,u\mfa)=0,\]
   for a non-degenerate critical point. Then, we have the following:
   \begin{enumerate}
     \item When $t\leq0$, then for a given $[\mfa]$ and a non-negative integer $d_0$, there are only finitely many  pairs $([\mfb],z)$ such that the moduli space $\check{M}^+_z([\mfa],[\mfb])$ is non-empty and of dimension at most $d_0$.
     \item When $t>0$, then for a given $[\mfa]$, there are only finitely many pairs $([\mfb],z)$ such that the moduli space $\check{M}^+_z([\mfa],[\mfb])$ is non-empty.
   \end{enumerate}
 \end{prop}
% The proof of the above proposition is to combine Proposition  \ref{prop-16.4.3} and follow the arguments of \cite[Proposition 29.2.1]{KM}.

 \begin{pf}
%By the hypothesis,  there is a  number $t$ such that  \[\mathcal E^{top}_{c}(z)+tgr_z([\mfa],[\mfb])\] is independent of homotopy class of the path $z$. We divide the proof into two cases.
\begin{itemize}
  \item When $t>0$, we repeat the same argument of Proposition \ref{prop-16.4.3} to get the conclusion.
  \item
     When $t=0$, and the moduli space $M_{z_i}([\mfa],[\mfb])$ are non-empty. It is known that the image of the critical-points set under the blow-down map $\pi:\mcB^\sigma(M,F,\mfs)\to \mcB(M,F,\mfs)$  is a set of  finite points. We may assume that $\pi[b_i]=[\beta]$ for all $i$. $\L_\omega$ descends to  a single-valued function on $\mcB(M,F,\mfs)$, hence the energy  of the trajectories in  these moduli spaces has  an up-bound. For the blow-down case, Proposition \ref{cmpt-prop} implies that there are only finitely many choices for the homotopy class of the path $\pi(z_i)$ in $\mcB(M,F,\mfs)$. In addition, the  dimension $d_0$ gives a lower-bound and up-bound for $\iota([b_i])$,  there are only finitely many $[b_i]$.
  \item When $t<0$,  there is a negative number $t$ such that
  $\mathcal E^{top}_{\mathfrak q}(z)+tgr_z([\mfa],[\mfb])$ is independent of $z$. To give a bound for the dimension, it suffices to give a bound for $gr_z([\mfa],[\mfb])$. Since $t<0$, we have an above bound for $\mathcal E^{top}_{\omega,\mathfrak q}(z)$. Since the dimension is bounded by $d_0$ and $[\mfa]$ is fixed,  there are finitely many pairs $([\mfb],z)$ such that $\iota([\mfb])$ is bounded above and below, and $gr_z([\mfa],[\mfb])\geq0$. The energy bound implies that only finitely many
   moduli spaces which are non-empty, and there are only finitely many critical points in the absence of reducibles, so the conclusion also holds.
\end{itemize}
 %   \end{itemize}
\end{pf}

It is known that $gr(\mfa,u\mfa)$ equals to the index of basic Dirac operator on $M\times S^1$,  by \cite{APS}. We rewrite the above formula in Proposition \ref{prop-29.2.1} as
\[\int_M(c_1(\mfs)-c)\wedge[u]\wedge\chi_F+t(\cdot\int_{\bar M_0\times S^1/\bar F}A_{0,b}|\tilde{dx}|+\sum^r_{j=1}\beta(M_j\times S^1))=0,\]
where
  \[\beta(M_j\times S^1)=\frac12\sum_\tau\frac1{n_\tau rank(W^\tau)}(-\eta(D^{S^+,\tau}_j)+h(D^{S^+,\tau}_j))\int_{\bar M_j\times S^1/\bar F}A^{\tau}_{j,b}(x)|\tilde{dx}|,\]
    the integrands $A_{0,b},~A^{\tau}_{j,b}(x)$ are similar to the Atiyah-Singer integrands,  $\bar M_0\times S^1$ is the principal domain of $M\times S^1$ and $\bar M_j\times S^1$'s are the finite desingularities of $M\times S^1$, more details are explained in the paper \cite{BKR}.

% Since, we have two types of the   monopole Floer homology groups, a natural question will be posed: What is the relationship between these two types?

%We have the theorem below, which is similar to \cite[Section 31]{KM}. \begin{thm}\label{thm-orbifold-Floer-homology-isomorphism}   Let $Y$ be an oriented closed $3$-orbifold, $\mfs$ be an orbifold \spinc structure and $c$ be  defined as above. Suppose that $\Gamma$ is a complete local system, e.g. Novikov ring and $c=c_1(\mfs)$. Then, we have the following isomorphisms:   \[\widecheck{HM}_*(M,F,\mfs,c;\Gamma)=\widecheck{HM}_*(M,F,\mfs;\Gamma).\] \end{thm}

\section{Examples}
In this section, we will give a family of   manifold with foliation satisfying the Assumption \ref{assum-main-1}.

\subsection{Fibration and orbifold}

%At the beginning of this subsection, we   focus on manifolds which admit a fibration over a closed oriented three manifold.

The easiest model is to consider $M=Y\times F$, where $Y$ is a closed oriented $3$  manifold and $F$ is a closed oriented manifold. Given a metric $g_Y$ and  a \spinc  structure $\mathfrak s$ of $Y$, by pulling back, one has  a data $(M,F,\pi^*g_Y\oplus g_F, \pi^*\mathfrak s)$, where $\pi:M\to Y$. Such a manifold with foliation $(M,F)$ satisfies the  Assumption \ref{assum-main-1}. We can generalize the global product model to the local product model, i.e. the fibration over $Y$.

\noindent Let $Y$  be a closed oriented $3$  manifold, and $M\to Y$ be a fibration over $Y$, such that $M$ is closed and oriented. Fix  a metric $g_Y$ and  a \spinc structure $\mathfrak s$ of $Y$, via pulling back, we have a bundle like metric and a transverse \spinc structure, still denoted by $\mathfrak s$. Since the volume form of $Y$ is closed, by pulling back, one has that $H^3_b(M)\neq0$. We have that $(M,F)$ satisfies the Assumption \ref{assum-main-1}, by Proposition \ref{prop-taut}. By the identification between the basic forms(sections) of $M$  and the forms(sections) of $Y$, one establishes the proposition below.

 \begin{prop}\label{prop-fibration}
   Let $(M,F,\mathfrak s)$ be defined as above.  Then, 
    %\begin{itemize}
     % \item when $b_1(Y)>1$, the basic Seiberg-Witten Floer homology $HF(M,F,\mathfrak s)$ is isomorphic to the Seiberg-Witten Floer homology of $(Y,\mathfrak s)$;
      %\item in general,
       the  basic monopole   Floer homology groups $\overline{HM}_*(M,F,\mathfrak s)$, $\widehat{HM}_*(M,F, \mathfrak s)$, $\widecheck{HM}_*(M,F, \mathfrak s)$ are isomorphic to the basic monopole  Floer homology groups $\overline{HM}_*(Y,\mathfrak s)$, $\widehat{HM}_*(Y, \mathfrak s)$, $\widecheck{HM}_*(Y,\mathfrak s)$ respectively with any coefficient.
  %  \end{itemize}
 \end{prop}

One can generalize the model of  fibration over manifold to the model of fibration over orbifold. First, we recall the notion of orbifold, which was first introduced by Satake \cite{Satake}.
\begin{defi}[c.f. \cite{Bald}]
 An $n$-dimensional orbifold $Y$ is a Hausdorff space $|Y|$ together with an atlas $(\{U_i\},\{\phi_i\},\{\tilde U_i\},\{\Gamma_i\}$, with transition maps $\{\phi_{ij}\}$,  which satisfies
 \begin{itemize}
     \item $\{U_i\}$ is locally finite;
     \item $\{U_i\}$ is closed under finite intersections;
     \item For each $U_i$, the finite group $\Gamma_i$ actions smoothly and effectively on a connected open subset $\tilde U_i\subset \mathbb R^n$, and there is a homeomorphism $\phi_i:\tilde U_i/\Gamma_i\to U_i$;
     \item If $U_i\subset U_j$, then there exists a monomorphism $f_{ij}:\Gamma_i\to \Gamma_j$ and a smooth embedding $\phi_{ij}: \tilde U_i\to \tilde  U_j$ such that for any $g\in\Gamma_i,~x\in\tilde U_i$, we have that $\phi_{ij}(g\cdot x)=f_{ij}(g)\cdot \phi_{ij}(x)$ and the following diagram commutes: \[\xymatrix{    \tilde U_i \ar[d]\ar[r]^{\phi_{ij}} &\tilde U_j\ar[d] \\    \tilde U_i/\Gamma_i\ar[r]^{f_{ij}}\ar[d]_{\phi_i}  & \tilde U_j/\Gamma_j\ar[d]_{\phi_j}\\
     U_i\ar[r]& U_j}\]
   where $f_{ij}$ is induced by the monomorphism and the canonical projection.
     %\[\xymatrix{\tilde U_i \ar[d]\ar[r]^{\phi_{ij}}\tilde U_j\ar[d]\\}\]
 \end{itemize}
\end{defi}

An $n$-dimensional orbifold bundle over $Y$ is defined in the similar manner.

\begin{defi}[c.f. \cite{Bald}]
   An orbifold $E$ is called an orbifold    bundle over $Y$, if there exists  a smooth orbifold map $p:E\to Y$, such that
    \begin{itemize}
      \item there is an atlas $(\{V_i\},\{\tilde V_i\},G_i)$ of $E$, satisfying $V_i=p^{-1}(U_i)$ and $\tilde V_i=\tilde U_i\times E_0$, where $(\{U_i\},\{\phi_i\},\{\tilde U_i\},\{\Gamma_i\} $ is an atlas of $Y$ and $E_0$ is a standard fiber;
      \item the following diagram commutes
      \[\xymatrix{    \tilde U_i\times E_0 \ar[d]\ar[r]^{\tilde p} &\tilde U_i\ar[d] \\    \tilde V_i/G_i\ar[d] & \tilde U_i/\Gamma_i\ar[d] \\
     V_i\ar[r]^p& U_i}\]
    \end{itemize}
    where $\tilde p$ is a $(G_i,\Gamma_i)$-equivariant map.
 \end{defi}
When $G_i$ acts freely, $E$ becomes a manifold, e.g. the frame bundle of an oriented orbifold(see \cite[Theorem 1.3]{ALR}).
Let $Y$ be an oriented closed $3$-orbifold. Suppose the singular set $\Sigma Y=\{x\in Y| ~G_x\neq{1}\}$ is a set of disjoint union of finite circles, where $G_x$ denotes the isotropy group at $x$. We rewrite $$\Sigma Y=\cup_{1\leq i\leq n}l_i$$ and each circle $l_i$ is assigned a positive integer $\alpha_i$ given by its isotropy group $\mathbb Z_{\alpha_i}$. Let $D$ be the unit disk and $\mathbb Z_{\alpha_i}$ act on it by rotation.  Near each $l_i$, we have an atlas,
\[\phi_i:(S^1\times D,S^1\times\{0\})\to (U_i,l_i),\] where $\phi_i$ induces a homeomorphism from $((S^1\times D)/\mathbb Z_{\alpha_i}, S^1\times\{0\})$ to $(U_i,l_i)$. It is known that $TY$ always lifts to an orbifold $spin^c$-bundle for such a $3$-orbifold. The definition of the Seiberg-Witten invariant can be generalized to 3-orbifold, see  Baldridge \cite{Bald} and Chen \cite{Chen}. For Seiberg-Witten invariant, we have the following proposition, which is similar to the manifold case.

\begin{prop}
  Let $Y$ be a closed oriented $3$-orbifold and $M\to Y$ be a fibration over $Y$. Suppose that $\mathfrak s$ is a transverse \spinc structure which comes from the pull-back \spinc structure of  $Y$ and $M$ is a closed oriented manifold. Then, we have that   basic Seiberg-Witten invariant of $M$ is equal to the Seiberg-Witten invariant of $Y$, for $b^1(Y)>1$.
\end{prop}

%The proof is no essential difference to the codimension $4$-case(see \cite[Theorem 28]{KLW}, \cite[Section 6.7]{Morgan}), here we omit it.

\noindent Under tensor product, the topological isomorphism classes of orbifold line bundles form a group. We give a local description for each class of  such a group. We have an orbifold line bundle over $Y$, which is a  trivial line bundle over $Y\setminus \Sigma Y$, and over each $U_i$, it is given by $(S^1\times D\times \mathbb C)/\mathbb Z_{\alpha_i}$, where  $  \mathbb Z_{\alpha_i}$ action is defined by,
\[a\cdot(t,w,z)\mapsto (t,e^{\frac{2\pi ia}{\alpha}}w,e^{\frac{2\pi ia}{\alpha}}z),\] for each element $a\in \mathbb Z_{\alpha_i}$. This bundle is glued together by a transition function $\varphi(t,w)=w$ on the overlap $\partial(S^1\times D)$. Each $l_i$ generates a line bundle $E_i$. Let $L$ be a line bundle over $Y$. There is a collection of integers $\{\beta_1,\cdots,\beta_n\}$ such that
\begin{itemize}
  \item $0\leq\beta_i<\alpha_i$, for each $i=1,\cdots,n$;
  \item the bundle $L\otimes E^{-\beta_1}_1\cdots\otimes E^{-\beta_n}_n$ is trivial over each neighborhood of $l_i$.
\end{itemize}
 By forgetting the orbifold structure, it can be naturally identified with a smooth line bundle (denoted by $|L|$) over the smooth manifold $|Y |$. We will list some necessary results of such orbifolds.% to define the basic monopole Floer homologies.

\begin{thm}[Baldridge \cite{Bald}]
  The tangent bundle $TY$ lifts to an orbifold \spinc bundle.
\end{thm}
%\begin{pf}  We need to show the split of the tangent bundle. Since each $l_i\times D/\mathbb Z_{\alpha_i}$ comes with a $\mathbb Z_{\alpha_i}$-invariant oriented vector field tangent to $l_i$ at each point of $D/\mathbb Z_{alpha_i}$. This vector field induces a non-zero section $\partial Y'\to TY|_{\partial Y'}$. We remove an extra $S^1\times D$ from the interior of $Y'$ and put a similar nonzero section on the boundary. It is known that  the obstruction to extending the section into the interior of  \[Y''=Y\setminus (S^1\times D\amalg l_i\times D/\mathbb Z_{\alpha_i})\]  corresponds to an element of $H^3(Y'',\partial Y'';\pi_2(S^2))\cong \mathbb Z$. Using the relation  \[[\partial(S^1\times D)]=-\sum_i[\partial(l_i\times D/\mathbb Z_{\alpha_i})],\]  the obstruction can be removed by changing the framing on the boundary of $S^1\times D$. \end{pf}

\begin{lemma}[Chen \cite{Chen}]
  Let $Y$ be defined as above. Then we have that
  \[\pi_0(C^\infty(Y,S^1))\cong H^1(|Y|,\mathbb Z).\]
\end{lemma}
%\begin{pf}  Given any $f\in Map(Y,S^1)$, it induces a function $f'$ on $|Y|$. We set $[f]$ as the element in $H^1(|Y|,\mathbb Z)$. We need to show that the map  \[\pi_0(Map(Y,S^1))\to H^1(|Y|,\mathbb Z)\] is onto and one-to-one.  \noindent  Let $[u]\in H^1(|Y|,\mathbb Z)$, there is a smooth function $u'$ of $|Y|$ which represents $[u]$. We put back $u$ on $Y$, which is denoted by $u$, since locally the projection $Y\to |Y|$ is given by $(t,z)\mapsto (t,z^\alpha)$ for some $\alpha>\geq1$, where $t\in\mathbb R$ and $z\in\mathbb C$. One can verify that $u$ represents $[u]$, which shows that  the map  \[\pi_0(Map(Y,S^1))\to H^1(|Y|,\mathbb Z)\] is onto.  \noindent Let $u_1,~u_2\in Map(Y,S^1)$ such that $[u_1]=[u_2]$. The induced continuous functions $u'_1,~u'_2$ on $|Y|$ are homotopic-equivalent. We perturb $u'_1,~u'_2$  into the smooth functions $u'_{1,+\epsilon},~u'_{2-\epsilon}$ on $|Y|$ through a family functions $\tilde u_t$, with $1\leq t\leq 1+\epsilon$ and $2-\epsilon\leq t\leq 2$ respectively. The pull back of $\tilde u_t$ to $Y$ are continuous, which can be perturbed to  a family of smooth functions, so that are equivariant with respect to the local group actions. This implies that $u_1$ and $u_2$ are homotopic-equivalent to each other, i.e. $u_1$ and $u_2$  define the same class in $\pi_0(Map(Y,S^1))$.\end{pf}

\begin{prop}\label{prop-foliated-deRham}
  Let $Y$ be the orbifold  as before. Then, we have the following isomorphism
  \[H^*(|Y|,\mathbb R)\cong H^*_{dR}(Y,\mathbb R). \]
\end{prop}
\begin{pf}
  We have the fine resolution below for orbifold $Y$,
  \[0\to\mathbb R\to \mathcal A^0\overset{d}\to\mathcal A^1\cdots\]
  of the constant sheaf $\mathbb R$. By the double complex argument, we have the isomorphism $$\check H^*(Y,\mathbb R)\cong H^*_{dR}(Y,\mathbb R), $$
  where the first cohomology group is the \v{C}eck-cohomology group. Since we can find a finite covering $\{\mathcal U_i\}$, such that all non-empty intersections of finitely-many sets are contractible,  \v{C}ech cohomology is   isomorphic to the singular cohomology of the of the underlying  space $|Y|$. Thus, we have that
  \[H^*(|Y|,\mathbb R)\cong \check H^*(Y,\mathbb R)\cong H^*_{dR}(Y,\mathbb R). \]
\end{pf}

For an oriented closed $3$-orbifold $Y$ with a metric $g$ and \spinc structure $\mathfrak s$ whose determinant line bundle has the Seifert data $(b,\beta_1,\cdots, \beta_n)$, one  defines the Chern-Simons-Dirac functional
\[\L(A,\Psi)=-\frac18\int_Y(A^t-A^t_0)\wedge(F_{A^t}+F_{A^t_0})
+\frac12\int_Y(\Psi,\Dirac_A\Psi)dvol_Y,\]
for any $(A,\Psi)\in \mathcal C(Y,\mfs)$. Let $u\in\mathcal G(Y)$, we have that
\[\L(A,\Psi)-\L(u(A,\Psi))=-\frac12\int_Yu^{-1}du\wedge F_{A^t_0}=-2\pi^2\langle c_1(\mathfrak s),[u]\rangle,\] where $c_1(\mathfrak s)=[\frac{i}{2\pi}F_{A^t_0}]$ and $[u]=[\frac{-i}{2\pi}u^{-1}du]$. Similar to the manifold case, we define the critical points of the Chern-Simons-Dirac functional and the blow-up configuration space. % Suppose that $\mfa\in Crit^\sigma(\L)$ is non-degenerate, it is known that the grading between $\mfa$ and $u\mfa$ is equal to\[gr(\mfa,u\mfa)=spec(\Dirac_{A,u(t)A})=Ind(\Dirac_{u(t)}|_{Y\times S^1})=\langle c_1(\mathfrak s),[u]\rangle-\sum_i\frac{\beta_i}{\alpha_i}\int_{l_i}[u],\]where $A$ is the connection component of $\mfa$, $spec(\Dirac_{A,u(t)A})$ denotes the spectral flow between $\Dirac^A$ and $\Dirac^{uA}$ and we applied the Index formula on orbifold(see Kawasaki\cite{Kawasaki}) for the third equality.
By the Proposition \ref{prop-foliated-deRham} and the Poincar\'e duality,  it is known that there is a unique second cohomology class $c(\mfs,Y)\in H^2(Y,\mathbb R)$ such that
%\[\int_Yc(\mfs|_{\Sigma(Y)})\wedge [u]=\sum_i\frac{\beta_i}{\alpha_i}\int_{l_i}[u],\]for any $u\in Map(Y,S^1)$. Such a class $c$ is uniquely determined by $c_1(\mfs)$ and the singular point set $\Sigma Y$. We rewrite the formula
 \[gr(\mfa,u\mfa)=\langle c_1(\mfs)-c(\mfs,Y),[u]\rangle, \]
 where $gr(\mfa,u\mfa)$ denotes the grading between $\mfa$ and $u\mfa$ for a  non-degenerate critical point $\mfa\in Crit^\sigma(\L)$. Using a complete local system $\Gamma$, we can construct the monopole Floer homologies for $(Y,\mfs)$.

 When $c(\mfs,Y)$ is propositional to $c_1(\mfs)$, i.e. there is a real constant $k$ such that
 \[c(\mfs,Y)=kc_1(\mfs).\] Suppose that $k\neq1$, then we can find a real constant $t$, such that
 \[\L_\omega(A,\Psi)-L_\omega(u(A,\Psi)+t gr(\mfa,u\mfa)=0,\]
 which is equivalent to the formula
\begin{equation}
    -c_1(\mfs)+t(c_1(\mfs)-c(\mfs,Y))=0.\label{formula-orbifold}
  \end{equation}
% \begin{defi}   We say $c(\mfs,Y)$ is monotone, if $c(\mfs,Y)$ is proposition to $c_{1}(\mfs)$.  It is called positively or negatively monotone, if there is a positive or negative real number $t$, such that  \[\L_\omega(B,\Psi)-L_\omega(u(B,\Psi)+t gr(\mfa,u\mfa)=0.\] Moreover, if such a number $t$ is zero, we call $c(\mfs,Y)$ is balanced. \end{defi}

 \begin{prop}\label{prop-29.2.1-energy-bound}
 Let $(Y,\mfs)$ be a closed oriented $3$-orbifold as above.
  Suppose that  and  all the moduli spaces $M_z([\mfa],[\mfb])$ for the   perturbation $\mathfrak q$ are regular and  the formula \eqref{formula-orbifold} holds for each non-degenerate critical $\mfa$ and $(A,\Psi)=\pi(\mfa)$.  Then, the following holds:
  \begin{enumerate}
     \item When $t\leq0$, then for a given $[\mfa]$ and a non-negative integer $d_0$, there are only finitely many  pairs $([\mfb],z)$ for which the moduli space $\check{M}^+_z([\mfa],[\mfb])$ is non-empty and of dimension at most $d_0$.
     \item When $t>0$, then for a given $[\mfa]$, there are only finitely many pairs $([\mfb],z)$ for which the moduli space $\check{M}^+_z([\mfa],[\mfb])$ is non-empty.$([\mfb],z)$ for which the moduli space $\check{M}^+_z([\mfa],[\mfb])$ is non-empty and has dimension no more than $d_0$.
  \end{enumerate} % \item Suppose that $c(\mfs,Y)$ is balance. Then, for a given $[\mfa]$ and $d_0\geq0$, there are only finitely many pairs $([\mfb],z)$ for which the moduli space $\check{M}^+_z([\mfa],[\mfb])$ is non-empty and has dimension no more than $d_0$.    \item Suppose that $c(\mfs,Y)$ is negatively monotone, and the perturbation has been chosen so that  there are no reducible solutions. Then, the above conclusion holds.    \item Suppose that $c(\mfs,Y)$ is positively monotone, and the perturbation has been chosen so that  there are no reducible solutions. Then,  for a given $[\mfa]$, there are only finitely many pairs $([\mfb],z)$ for which the moduli space $\check{M}^+_z([\mfa],[\mfb])$ is non-empty. \end{itemize}
\end{prop}
The proof is similar to Proposition \ref{prop-29.2.1}, here we omit it.

 The space of broken trajectories $\check M^+_z([\mfa],[\mfb])$ can be identified with the manifold model. This space is still compact for fixed $[\mfa]$, $[\mfb]$ and $z$ as in \cite[Theorem 16.1.3]{KM}. We apply the same arguments of \cite[Section 20-Section 25]{KM} or of the previous section to establish the following theorem.

 \begin{thm}
  Let $\Gamma$ be any local system of abelian groups on $\mcB^\sigma(Y,\mfs)$ and let $(Y,\mfs)$ be a closed oriented $3$-orbifold as above. Suppose that the the formula \eqref{formula-orbifold} holds for each non-degenerate critical point. Then we construct  the basic monopole Floer homologies
  \[\widecheck{HM}_*(Y,\mfs; \Gamma), ~ \widehat{HM}_*(Y,\mfs;\Gamma),~ \overline{HM}_*(Y,\mfs;\Gamma).\]
\end{thm}

%By using the above local system $\Gamma$, the maps $\partial^o_o$ el at  are well-defined. For instance, we consider\[C^o(Y,\mfs,c,\Gamma)=\bigoplus_{[\mfa]}\mathbb Z\Lambda[\mfa]\otimes \Gamma[\mfa], \] where $[\mfa]$ denotes the irreducible critical point.We define the partial\[\partial^o_o =\sum_{[\mfa]}\sum_{[\mfb]}\sum_z\sum_{[\gamma]}\in \check M_z([\mfa],[\mfb])\epsilon[\gamma]\otimes \Gamma(z),\]where the sum is over all the moduli space $M_z([\mfa],[\mfb])$ with dimension $1$ and $[\mfb]$ denotes the irreducible critical point. The contribution for a given pair of critical points takes the form\[\sum_zn_z\Gamma(z).\]We set the support as\[supp(n)=\{z| ~n_z\neq0\}.\]By the definition of c-complete, we havethat\[supp(n)\cap \{z|\mathcal E^{top}_{\omega,\mathfrak q}(z)\leq C\}\] is finite, using the completeness of the local system, we have that the form\[\sum_zn_z\Gamma(z)\] is convergent. Similarly, one can verify that the maps $\check\partial$, $\hat\partial$ and $\bar\partial$ are well-defined.Combining with the equicontinuous property of the local system $\Gamma$, the proofs of $\check\partial^2$, $\hat\partial^2$ and $\bar\partial^2$ go through as the classical case.

\noindent We give an example of such a complete local system, i.e. a Novikov ring \cite{Novikov}. Let $I\subset\mathbb R$ be the set of the image of the homomorphism
\[\mathcal E^{top}:\pi_1(\mcB^\sigma(Y,\mfs))\to \mathbb R,~z\mapsto \mathcal E^{top}(z),\] where $\mathcal E(z)$ denotes the difference of the Chern-Simons-Dirac functional  between the different representatives of the quotient point. Set $\mathbb F=\mathbb Z_2$. We define $\mathbb F[I]$ by
 \[\mathbb F[I]=\{\sum _{i\in I} r_it^i|\mbox{ only for finitely many }i,~r_i\ne0\}.\]For $k\in\mathbb R$, let $U_{-k}$ be the $\mathbb F$-module spanned by the generators $t^i$, $i\in I$ satisfying $i\leq -k$. Using these as open neighborhoods of $0$, we form the completion $\bar{\mathbb F}[I]$, i.e.  each element is of the form $$\sum^C_{i= -\infty} r_it^i.$$  We define a local system by taking at $[\mfa_0]$ to be $\bar R[I]$, and specifying that for each closed loop $z$ based at $[\mfa_0]$, the automorphism $\Gamma(z)$ be the multiplication by   $t^{-\mathcal E^{top}(z)}$. This is a $c$-complete local system.

Similar to the foliation case of the previous section or to the non-exact perturbation on manifold case, we have the following theorem.
\begin{thm}
 Let $(Y,\mfs)$ be a $3$ orbifold defined as above with a \spinc structure $\mfs$. Then we have  the   monopole Floer homologies
  \[\widecheck{HM}_*(Y,\mfs,c;\Gamma), ~ \widehat{HM}_*(Y,\mfs ,c;\Gamma),~ \overline{HM}_*(Y,\mfs,c;\Gamma).\]
  where $\Gamma$ is a complete local system. Moreover, these homologies depend only on the isomorphism class of the \spinc structure $\mfs$, $c$ and $Y$,and are independent of the metrics or the  perturbations.
\end{thm}

  We give a necessary condition to avoid the complete local system.
% Before we proceed, we give a review of classical non-exact perturbation(see \cite[Section 29-Section 30]{KM}).Let $Y$ be a closed oriented three manifold and $\mfs$ be a \spinc structure. Choosing a class $c\in H^2(Y,\mathbb R)$, we consider the (non-exact)perturbed Chern-Simons-Dirac functional, \[\L_\omega(B,\Psi)=\L(B,\Psi)-\frac12\int_Y(B-B_0)\wedge\omega,\] where $\omega\in \frac{2\pi}ic$ and $(B,\Psi)\in \mcC(Y,\mfs)$. By the argument of \cite[Section 30]{KM}, we can define the basic monopole Floer homology groups with Novikov ring as coefficient.

% \begin{thm}[c.f. {\cite[Section 29.2]{KM}}]\label{thm-perturbation-no-need-Novikov}   Let $(Y,\mfs,c)$ be given as above. Suppose that there is a constant $t\in\mathbb R$ such that   \[2\pi^2c_1(\mfs)-c=t2\pi^2c_1(\mfs).\]   Then, the Novikov ring is not necessary, i.e. for any local coefficient, one can define the basic monopole Floer homologies. \end{thm}
% The proof is given in \cite[Section 29.2]{KM}. Similar to Theorem \ref{thm-perturbation-no-need-Novikov}, we have the following theorem.

 \begin{thm}
   Let $Y$ be a closed oriented $3$-orbifold defined as above, and $(\mfs,c)$ be as above. Suppose that there is a  constant $t$ such that $$-(c_1(\mfs)-c)+t(c_1(\mfs)-c(\mfs,Y))=0.$$ Then, for any local system the monopole Floer homologies are well-defined.
 \end{thm}

\subsection{Suspension}

Another way to construct the foliation is by suspension, here we give two   references of this subsection, see \cite[Chapter 3.8]{Molino} and \cite{Richard}. Let $(Y,g)$ be a closed oriented $3$ Riemannian manifold. Suppose that a compact Lie group $G$ actions on $(Y,g)$ isometrically and preserving the orientation of $Y$, and we have a representation
\[f: \pi_1(X)\to G\]
such that the closure of $Im(f)$ is $G$, where
  $X$ is a closed oriented manifold with fundamental group $\pi_1(X)$.
We set $M=\tilde X\times Y/f$, where $\tilde X$ denotes the universal covering of $X$ and $(x,y)\sim (x[\gamma]^{-1},f([\gamma])y)$ for $[\gamma]\in \pi_1(X)$. Fixing  a point $p= [y_0,x_0]\in M$, its leaf is defined by the set of the form\[\mathcal F_{p}=\{[x,y_0]\big| ~x\in \tilde X\}. \]

\noindent
Since one can find a $G$-invariant volume form over $Y$, lifting back on $M$ it holds that $H^3_b(M, F)\neq0$,  which implies  that the foliation is taut by Proposition \ref{prop-taut}. Before preceding, we have the following lemma(see \cite{Lin}).
\begin{lemma}
  Let $(M,F)$ be defined as above. Then,  we have an identification
   \[\pi_0(Map^G(Y,S^1))\cong   H^1(M,\mathbb Z)\cap H^1_b(M),\]
   where $Map^G(Y,S^1)$ denotes the space of $G$-invariant $S^1$-valued functions.
\end{lemma}
%\begin{pf}  It is known that $H^1(M,\mathbb Z)\cong \pi_0(Map(M,S^1))$, which means that  for each element $[w]\in H^1(M,\mathbb Z)$ we have a  representation $u:M\to S^1$ of this homotopy class $[u]_{ht}\in\pi_0(Map(M,S^1))$  satisfying the condition:   \[ [\frac1{2\pi i}u^{-1}du]=[w].\] Therefore, any element $[w]\in H^1(M,\mathbb Z) \cap H^1_b(M)$ corresponds to a representation $u$ of the homotopy class $[u]_{ht}\in\pi_0( Map(M,S^1))$, such that $[\frac1{2\pi i}u^{-1}du]\in H^1_b(M)$. This implies that there is $f\in i\Omega^0(M)$ such that $u^{-1}du+df\in \Omega^1_b(M)$. Setting $u'=e^{f}u$, we have that $L_\xi u'\equiv0$ for any $\xi\in \Gamma(F)$. This implies that   \[H^1(M,\mathbb Z)\cap H^1_b(M)\cong\pi_0(Map_b(M,S^1)), \]   where $Map_b(M,S^1)=\{u|u\in Map(M,S^1)~L_\xi u\equiv0,\mbox{ for any }\xi\in \Gamma(F)\}$. By the above argument, we have that $u$ corresponds to a $G$-invariant function on $Y$, we still use the same notation $u$ to express this $G$-invariant function of $Map(Y,S^1)$. Therefore, Thus, a $\mathbb Z$-module subset $\Gamma$ of $H^1_b(M)$ is a lattice of   $H^1_b(M)$ if and only if $rank(\Gamma)=b^G_1(Y)$.\end{pf}

Suppose there is a $G$-equivariant \spinc structure.  Given a $G$-equivariant spinor bundle $$S'\to Y,$$ we construct a foliated spinor bundle $S=\tilde X\times S'/f$, where the action of $[\gamma]\in \pi_1(X)$ is defined by $[\gamma](x,s_p)=(x[\gamma]^{-1},f[\gamma]s_p)$.
By \cite{Richard}, it is known that there is  an identification\[\Gamma^G(Y,S')\cong\Gamma_b(M,S).\]  %We define the groupoid $\mathcal G_F$ of  $(M,F)$ by $$\mathcal G_F=\{(x_1,x_2,[\gamma])\},$$ where $x_1$ and $x_2$ are points of a leaf $\mathcal F_y$, and $[\gamma]$  is an equivalence class of piecewise smooth paths in $\mathcal F_y$ from $x_1$ to $x_2$, two paths $a$ and $b$ are equivalent if and only if $a^{-1}b$ has a trivial holonomy. The multiplication is defined by $(x_1,x_2,[\gamma])(x_2,x_3,[\gamma'])=(x_1,x_3,[\gamma\gamma'])$. Let $p_1=[x_1,y_1]$ and $p_2=[x_2,y_2]$ be two points of a leaf of $M$, and $g=(p_1,p_2,[\gamma])\in\mathcal G_F$. Since $p_1$ and $p_2$ stay in the same leaf, we may choose a different representation of $p_2$ such that $y_2=y_1$. We define the action $ g_*: S_{p_1}\to  S_{p_2}$, by $[x_1,v_{y_1}]\mapsto [x_2,v_{y_1}]$, where $v_{y_1} \in S'$ denotes an element of the fiber over $y_1$. One can verify that $g_*g'_*=g'g_*$.  Let  $s:Y\to S'$ be a $G$-invariant section.  We set  $Ks:M\to S$ by $Ks([x,y])=[x,s(y)]$. The identity\[[x\gamma^{-1},s(f(\gamma)y)]=[x\gamma^{-1}, L_{f(\gamma)}s(y)]\] implies that $K$ is well-defined. For a $G$-invariant section $s:Y\to S'$, we have that\[g_*Ks([x_1,y])=g_*([x_1,s(y)])=[x_2,s(y)]=Ks([x_2,y]),\]for $g=([x_1,y],[x_2,y],[\gamma])\in\mathcal G_F$. This implies that $K:\Gamma^G(Y,S')\to \Gamma_b(M,S)$. Conversely, let  $s:M\to S$ be a basic section. Fixing $x\in\tilde X$ and $y\in Y$, the formula $s([x,y])=[xv_y]$ uniquely determines a vector $v_y\in S_y$. Since $s$ is basic, we have that $v_y$ is independent of the choice of this fixed $x$, i.e. $s$ defines a section $s:Y\to S'$.  Again, by using the formula $[x\gamma^{-1},s(f(\gamma)y)]=[x\gamma^{-1}, L_{f(\gamma)}s(y)]$, we have that $s$ is $Im(f)$-invariant, hence $G$-invariant.For a $G$-equivariant Dirac operator $D$, we have a basic Dirac operator constructed by the above identification, which is still denoted by $D: \Gamma_b(M,S)\to \Gamma_b(M,S)$.
 Summarizing the above arguments, we have the following proposition(see \cite{Lin}).
\begin{prop}
   Let  $(M,F)$ be a manifold with  a foliation constructed  as above,  and let $Y$ admit a $G$-equivariant spinor bundle. Suppose  it holds that $rank(\pi_0(Map^G(Y,S^1)))=b^G_1(Y)$, where $Map^G(Y,S^1)$ denotes the set of  $G$-invariant $S^1$-valued functions and $b^G_1$ denotes the dimension of the first cohomology for the $G$-invariant deRham complex. Then $(M,F)$ satisfies the Assumption \ref{assum-main-1}.  % We can define the  Suppose that $H^1_b(M)\cap H^1(M,)$
\end{prop}

\noindent
{\bf Remark}:

%\begin{itemize}\item When $G$ is connected, we have that $$H^{i,G}_{dR}(Y)\cong H^{i}_{dR}(Y).$$ It suffices for show that $G\sigma=\sigma$ for each homology class $\sigma$. Let $K$ be a complex containing $\sigma$ and $K'$ be a complex containing $K$ in its interior. At the identity $e\in G$, we can find a small neighborhood $U\subset G$ such that $g:K\to K'$ is homotopic to $e:K\to K'$ for each $g\in U$. By the connectedness of $G$, it is proved. \end{itemize}
\begin{enumerate}
  \item The condition that $rank(\pi_0(Map^G(Y,S^1)))=b^G_1(Y)$ is necessary. For example, let $Y=T^3=(S^1)^3$, $G=S^1$ action canonically on the first slot of $Y$, and $X=S^1$ with $f:\pi_1(X)\to S^1$ by sending the generator element of $\pi_1(X)$ to a dense element of $S^1$, e.g. $1\mapsto e^{i2\pi\theta}$ for some $\theta\notin\mathbb Q$. We have that $\dim H^1_b(M)=\dim H^{1,G}_{dR}(Y)=3$, and $H^1_b(M)\cap H^1(M,\mathbb Z)\cong\pi_0(\{u:M\to S^1|\mbox{ such that} L_\xi u\equiv0,\mbox{ for any }\xi\in\Gamma(F)\}))\cong\mathbb Z^2$.
  \item When $G$ is connected, we have that $b^{G}_1=b_1$, since  any homology cycle $\sigma$ is homotopic to $g_*\sigma$, for any $g\in G$. When $G$ actions freely, we have that $\pi_0(Map^G(Y,S^1))\cong H^1(Y/G,\mathbb Z)$, and $b_1(Y/G)=b^G_1(Y)$.
\end{enumerate}

At the end of this section, we give an explicit example. Let $Y=SO(3)$, and $T_1,~T_2$ and $T_3$ be three maximal tori(circles), such that their Lie algebras span the Lie algebra of $Y$, i.e. $so(3)$. We choose a closed oriented manifold $X$ whose fundamental groups is isomorphic to $\mathbb Z*\mathbb Z*\mathbb Z$, e.g. $X=\sharp_3 S^1\times S^k$ with $k\geq2$. We can consider a family of representations
\[f_t: \pi_1(X)\to T_1, T_2, T_3\]
such that the first component $(1,0,0)$ sends to an element $\left(\begin{array}{ccc}
     \cos 2\pi t& -\sin 2\pi t &\\
     \sin2\pi  t&\cos 2\pi t &\\
     &&1\\
\end{array}\right)$ of $T_1$, the second component $(0,1,0)$ sends to an element $\left(\begin{array}{ccc}
     \cos 2\pi t& & \sin 2\pi t\\
      & 1&\\
     -\sin 2\pi  t&&\cos2\pi  t\\
\end{array}\right)$  of $T_2$ and the third component $(0,0,1)$ sends to an  element   $\left(\begin{array}{ccc}
1&&\\
     &\cos2\pi  t& -\sin2\pi  t \\
     &\sin2\pi  t&\cos2\pi  t \\
\end{array}\right)$  of $T_3$.  We set $M_t=Y\times \tilde X/f_t$, the codimension $3$ foliation $F_t$ on $M_t$ is defined by letting the leaves be
\[\mathcal F_{t,y}=\{[x,y]\big| ~x\in \tilde X\}. \] We choose a trivial $SO(3)$-equivariant spin structure of $Y$.
%Given a leaf $\mathcal F_{t,y}$, we have that $[g\cdot x,y]\in \mathcal F_{t,y}$ for any $g\in \pi_1(X)$. This is equivalent to say that\[[x,yf_t(g)]\in \mathcal F_{t,y}\]for any $g\in \pi_1(X)$. We have that$$\mathcal F_{t,y}\subset \mathcal F_{t,C_y}\subset \overline{\mathcal F}_{t,y}$$where $C_y$ denotes the closure of the orbit of $y$.Suppose that $S=\{[x_i,yf_t(g_i)]\}$ is a sequence of points of $ \mathcal F_{t,C_y}$. Let $\pi:\tilde X\to X$ be the covering map, $\{\pi(x_i)\}$ has convergent subsequence with limit $ \bar x\in X$. Without loss of generality, we let $\{\pi(x_i)\}$ be a convergent sequence. We may change the trivializing the cover at $\bar x$ such that $x_i\to x \in X$. Let $G_t$ be the closure of $Im(f_t)$ in $Y$. As $G_t$ is compact, we may assume that   $S=\{[x_i,yf_t(g_i)]\}$ has  $x_i\to x$ and $f_t(g_i)\to g_t$, i.e. $\{[x_i,yf_t(g_i)]\}$ converges to $[x,yg_t]$. Hence\[ \mathcal F_{t,C_y}= \overline{\mathcal F}_{t,y}.\]We also established a bijection\[Y/G_t\to M_t/\bar{F_t}. \]We want to show that it is also an isometric. Let $C_1$ and $C_2$ be two orbit in $G_t$, and $\ell$ be a minimal geodesic connecting them. For any $x\in \tilde X$, the curve $\ell'=[x,\ell]$ connecting $\mathcal F_{C_1}$ and $\mathcal F_{C_2}$. Since we are using the product metric, the length of any curve connecting these two can not be shorter than the length of $\ell'$. As $H^3_b(M_t, F_t)\neq0$. We have that $(M_t, F_t)$ admits a taut bundle-like metric.
\begin{itemize}
  \item When  the group $G_t$ is a finite group of $Y$. Since  $Y$ admits a metric of positive scalar curvature, by Proposition \ref{prop-fibration} and \cite[Proposition 36.1.3]{KM}, one  deduces that     \[\overline{HM}_*(M_t,F_t)\cong \mathbb F[U,U^{-1}],~      \widecheck{HM}_*(M_t,F_t)\cong \mathbb F[U,U^{-1}]/\mathbb F[U],~      \widehat{HM}_*(M_t,F_t)\cong \mathbb F[U].\]
  \item When $G_t$ is dense in $M_t$. Since the transverse \spinc structure $\mathfrak s$ is  trivial, then $\Gamma_b(S)=\mathbb C^2$.   By the argument at the beginning of this subsection, we have that  $\Omega^1_b(M_t)\cong \Omega^{1,G}(Y)\cong\mathbb R^3$. Let $g_Y$ be a bi-invariant metric of $Y$ with positive scalar curvature, then  the associated $\Dirac^{g_Y}$ is an $SO(3)$ equivariant, which corresponds to a basic Dirac operator $\Dirac_0$ with spin connection $A_0$. Their spectrums have a one-to-one corresponding. Therefore, the solutions of the basic Seiberg-Witten equations \eqref{eqn-SW-equation-1} corresponds to the solutions of $SO(3)$-invariant Seiberg-Witten equations, i.e.
       \[\begin{cases}
         \Dirac^A_0\Psi=0\\
         \frac12*_{g_Y}F_A=q(\Psi),
       \end{cases}\] where $(A,\Psi)\in\mathcal C^{SO(3)}(Y)=\{A_0+i\Omega^{1,SO(3)}(Y)\times \mathbb C^2\}$. Since $g_Y$ has a positive scalar curvature, it is known that $\Psi\equiv0$, i.e. there is no irreducible solution to  the basic Seiberg-Witten equations. Consider the reducible solutions, we have that $dA=0$ and $\Psi$ is an eigenvector of $\Dirac^A_0$. Since $A=A_0+a$ and $d A_0=0$, this implies that $a$ is closed. Combining with $H^{1,G}(Y)=0$, we have that $a=df$ for a $SO(3)$-invariant function $f$, which implies that $f$ is constant and $a=0$. Thus, all reducible solutions are eigenvector of $\Dirac_0$. Recall that $\Dirac_0=\sum_i e^i\nabla'_{e_i}$, where $\{e_i\}$ is an orthonormal frame of $TY$; and
        \[\nabla'=d+\frac12\sum_{i<j}\omega_{ij}e^ie^j\]
        where $d$ is the flat connection of $S'$ and $\omega_{ij}$ is the Levi-Civita connection associated to $g_Y$. Since $\Dirac_0$ is independent of the choice of the orthonormal frame, we choose a frame $\{e_1,e_2,e_3\}$  such that they are generated by the left action of a frame  $\{L_x,L_y,L_z\}$ of $T_1Y$, where $L_x=\left(\begin{array}{ccc}
          &&\\
          &&-1\\
          &1&\\
        \end{array}\right)$, $L_y=\left(\begin{array}{ccc}
          &&1\\
          &&\\
          -1&&\\
        \end{array}\right)$ and $L_z=\left(\begin{array}{ccc}
          &-1&\\
          1&&\\
          &&\\
        \end{array}\right)$. Since $g_Y$ is bi-invariant $\{e_1,e_2,e_3\}$ are left-invariant, we have that
        \[\omega_{12}(e_3)=\frac12g_Y([e_3,e_1],e_2)=\frac12g_Y(e_2,e_2)=\frac12,\]
        \[\omega_{13}(e_2)=\frac12g_Y([e_2,e_1],e_3)=-\frac12g_Y(e_3,e_3)=-\frac12,\]
        \[\omega_{23}(e_1)=\frac12g_Y([e_1,e_2],e_3)=\frac12g_Y(e_2,e_2)=\frac12.\]
        We inherit the convection $e^1\cdot e^2\cdot e^3=Id$ from the book \cite{KM}, by  assigning
        \[e^1\mapsto\left(\begin{array}{cc}
          i&\\
          &-i\\
        \end{array}\right),~e^2\mapsto\left(\begin{array}{cc}
          &-1\\
          1&\\
        \end{array}\right),~ e^3\mapsto\left(\begin{array}{cc}
          &i\\
          i&\\
        \end{array}\right).\]
        For any $\Psi\in\Gamma^{SO(3)}(Y,S')\cong\mathbb C^2$, we have that 
        \[\Dirac_0\Psi=\sum_ke^kd\Psi+\sum_{i<j,k}e^k\frac12\omega_{ij}(e_k)e^ie^j\Psi=
        \frac34\Psi.\]
        Hence, $\Dirac_0$ acts as a diagonal matrix with eigenvalues $(\frac34,\frac34)$. Summarizing the above arguments, we have that
        %Since  a section $\psi\in \Gamma(M,S)$ corresponds to an $SU(2)$-equivariant map $s:P\to \mathbb C^2$, where $P$ is the principal foliated spin(\spinc) bundle, via the commutative diagram       \[\xymatrix{P\ar[d]\ar[r]&P\times \mathbb C^2\ar[d]\\       M\ar[r]&S=P\times \mathbb C^2/SU(2).}\]By the construction $P=\tilde X\times Y\times SU(2)/f_t=M\times SU(2)$. A basic section $\psi\in \Gamma_b(S)$ corresponds to a constant map $s:P\to \mathbb C^2$. Fixing a basic spin connection $A$, for any $\xi\in\Gamma(TM)$, the derivative $\nabla^A_\xi\psi$ corresponds to $\xi^Ps$, where $\xi^P$ is a lift vector field of $\xi$ under the connection $\nabla^A$. Since $s$ is constant, $\xi^Ps=0$, which implies that he basic Dirac operator becomes a zero map on $\Gamma_b(S)$. Therefore, one   has
      \[\overline{HM}_*(M_t,F_t,\mfs)\cong \mathbb F\oplus \mathbb F  ,~ \widecheck{HM}_*(M_t,F_t,\mfs)\cong \mathbb F\oplus \mathbb F, \widehat{HM}_*(M_t,F_t,\mfs)\cong 0.\]
\end{itemize}
%It is known that $\mathcal C^\sigma(M_t,F_t,\mathfrak s_t)$ is homotopic to $\mathbb C^2\setminus\{0\}/\mathbb R^+$, after taking the gauge action, we have that $\mathcal B^\sigma(M_t,F_t,\mathfrak s_t)$ is homotopic to $S^2$, hence $H^2(\mathcal B^\sigma(M_t,\mathfrak s_t),\mathbb F)\cong F[U]/(U^2=0)$. We can choose a perturbation such that the operator $\Dirac_0+\mathfrak q(A,\Psi)$ has two distinct positive eigenvalues, for bigger one we assign its degree by two and for lower one we assign its degree by $0$.

\vspace{3mm}

{\bf Remark:}
{  Note} that  the closure of the image of $f_t$ can not be of one dimensional. Otherwise, let $H$ be this one-dimensional closed subgroup in $SO(3)$. It is well known that $H$ preserves a vector in $S^2\subset \mathbb R^3$, say $(x,y,z)^t\in S^2$.  We have that the group generated by $\left(\begin{array}{ccc}
1&&\\
     &\cos2\pi  t& -\sin2\pi  t \\
     &\sin2\pi  t&\cos2\pi  t \\
\end{array}\right)$ preserves $(x,y,z)^t$, which implies that either $y=z=0$ or $t=0,1$. We can apply the same arguments for the other two subgroups. In conclusion we have that either $x=y=z=0$ or $t=0,1$, which contradicts to our assumption.

\vspace{3mm}

College of Mathematics and Statistics, Chongqing University,
Huxi Campus, Chongqing, 401331, P. R. China

Chongqing Key Laboratory of Analytic Mathematics and Applications, Chongqing University, Huxi Campus, Chongqing, 401331, P. R.
China

 E-mail: lindexie@126.com

\end{document}